\documentclass{jams-l}

\RequirePackage{enumitem}
\RequirePackage{graphicx}

\newtheorem{theorem}{Theorem}[section]
\newtheorem{lemma}[theorem]{Lemma}
\newtheorem{proposition}[theorem]{Proposition}
\newtheorem{corollary}[theorem]{Corollary}

\theoremstyle{definition}
\newtheorem{definition}[theorem]{Definition}
\newtheorem{example}[theorem]{Example}

\theoremstyle{remark}
\newtheorem{remark}[theorem]{Remark}
\newtheorem{assumption}[theorem]{Assumption}

\numberwithin{equation}{section}

\providecommand{\abs}[1]{\lvert#1\rvert}
\providecommand{\norm}[1]{\lVert#1\rVert}

\begin{document}

\title[Localisation for constrained transports I]{Localisation for constrained transports I: theory}


\author{Krzysztof J. Ciosmak}
\address{Department of Mathematics, University of Toronto, Fields Institute for Research in Mathematical Sciences}
\curraddr{}
\email{k.ciosmak@utoronto.ca}
\thanks{The  author wishes to express his thanks to Robert J. McCann for comments that contributed to the improvement of the manuscript.\\
 Part of the research presented here was performed when the author was a post-doctoral fellow at the University of Oxford, supported by the European Research Council Starting Grant CURVATURE, grant agreement No. 802689. 
 Another part was performed when the author was a doctoral student at the University of Oxford, supported by St John’s College, Clarendon Fund and Engineering and Physical Sciences Research Council.
}


\date{}

\dedicatory{}

\begin{abstract}
We investigate an analogue of the irreducible convex paving in the context of generalised convexity.
Consider two Radon probability measures $\mu,\nu$ ordered with respect to a cone $\mathcal{F}$ of functions on $\Omega$ stable under maxima. Under the assumption that any $\mathcal{F}$-transport between $\mu$ and $\nu$ is local, we establish the existence of the finest partitioning of $\Omega$, depending only on $\mu,\nu$ and the cone $\mathcal{F}$, into $\mathcal{F}$-convex sets, called irreducible components, such that any $\mathcal{F}$-transport between $\mu$ and $\nu$ must adhere to this partitioning. \\
Furthermore, we demonstrate that a set, whose sections are contained in the corresponding irreducible components, is a polar set with respect to all $\mathcal{F}$-transports between $\mu$ and $\nu$ if and only if it is a polar set with respect to all transports.
This provides an affirmative answer to a generalisation of a conjecture proposed by Ob\l{}\'oj and Siorpaes regarding polar sets in the martingale transport setting. \\
Among our contributions is also a generalisation of the Strassen's theorem to the setting of generalised convexity.
\end{abstract}

\maketitle
\tableofcontents

\section{Introduction}\label{s:intro}

We develop a theory of transport of measures with additional constraints given by a convex cone stable under maxima. In particular, we generalise the results concerning the irreducible convex paving of the martingale transport theory. The irreducible convex paving was discovered in the one-dimensional setting by Beiglb\"ock, Nutz and Touzi in \cite{Beiglbock2017} and developed in the finite-dimensional setting by De March and Touzi in \cite{Touzi2019}.  
The results say that, roughly speaking, any martingale transport between two given probability measures $\mu,\nu$ is constrained by the finest collection of closed, convex sets, whose relative interia, called irreducible components, partition the underlying space up to a set of $\mu$-measure zero. That is, if we are to transport mass from one irreducible component, then we are only allowed to do so within its closure. Moreover, within the components, the mass can be transported with no further such restrictions. 

Our theory extends these results to the setting of transports between two given probability measures $\mu,\nu$ on a set $\Omega$ that obey constraints imposed by an arbitrary convex cone of functions on $\Omega$ that is stable under taking finite maxima. The elements of the paving turn out to be convex in generalised sense.

As we show in a companion paper \cite{Ciosmak2024}, our developments allow to carry over the results concerning the irreducible convex paving to the infinite-dimensional setting, as well as to the infinite-dimensional submartingale transports, generalising results of Nutz and Stebegg \cite{Nutz2018}. We also are able to include the localisation of the Monge--Kantorovich problem, studied by Klartag \cite{Klartag2017}, Cavalletti and Mondino \cite{Cavalletti2017,Mondino2017}. Novel applications include pavings for submartingale transports in the finite and infinite-dimensional setting, pavings for harmonic transports and pavings for the solutions of general martingale problems.

Our contributions are several-fold. Firstly, we provide a general, abstract framework, which allows for an extension of the theory of irreducible convex paving and martingale transports. 
The advantage of our approach, apart from its breadth, lies in its transparency, as compared to the approaches in \cite{Touzi2019} and in \cite{Obloj2017}. Secondly, the novel methods of proofs allow us also to give a description of the fine structure of the irreducible components. 
Thirdly, we also give a proof of the characterisation of polar sets in our abstract setting. Not only is the characterisation more general as compared to the one in \cite{Touzi2019}, but also it introduces novel approach, whose main idea is more intuitive and elementary that the one in \cite{Touzi2019}.


\subsection{General setting}

Let $\mathcal{G}$ be a cone of functions on $\Omega$,  and $\mathcal{F}$ the lattice cone generated by $\mathcal{G}$, meaning 
the set of functions realised as maxima of functions from $\mathcal{G}$, and let $\mathcal{H}$ be the complete lattice cone generated by $\mathcal{G}$, i.e., the set of functions realised as suprema of functions from $\mathcal{G}$.\footnote{See Section \ref{s:prelim} for the definitions of lattice cone and of complete lattice cone, Definition \ref{def:prelim}, as well as for the definitions of  other basic notions, if not defined before.}  For example,  $\mathcal{G}$ could consist of non-decreasing affine functions of one variable,  in which case $\mathcal{H}$ would be the set of all convex non-decreasing functions on the real line, or $\mathcal{G}$ could consist of all affine functions, and then $\mathcal{H}$ would be the set of all convex functions. These two basic examples we encourage the reader to keep in mind. 

We equip $\Omega$ with the weak topology $\tau(\mathcal{G})$ generated  by $\mathcal{G}$.

Let $p\in\mathcal{H}$ be a real-valued proper function, i.e., with compact preimages of compact subsets of the real line, bounded from below by a positive constant, and such that any function in $\mathcal{G}$ is of $p$-growth.

Let $\mu,\nu$ be two Radon probability measures on $\Omega$, with respect to which $p$ is integrable.
We shall say that $\mu$ and $\nu$ are in $\mathcal{F}$-order, and write $\mu\prec_{\mathcal{F}}\nu$, provided that for any $f\in\mathcal{F}$

\begin{equation*}
\int_{\Omega}f\, d\mu\leq\int_{\Omega}f\, d\nu.
\end{equation*}

We shall study $\mathcal{F}$-transports between $\mu$ and $\nu$, i.e., Radon probability measures $\pi$ on $\Omega\times\Omega$ such that the respective marginals of $\pi$ are $\mu$ and $\nu$ and for all $f\in\mathcal{F}$ 
and all measurable non-negative $h$ on $\Omega$ we have
\begin{equation*}
\int_{\Omega}h(\omega_1)f(\omega_1)\, d\mu(\omega_1)\leq \int_{\Omega\times\Omega}h(\omega_1)f(\omega_2)\, d\pi(\omega_1,\omega_2).
\end{equation*}
Equivalently, $\pi$ is an $\mathcal{F}$-transport between $\mu$ and $\nu$ if and only if, it is a distribution a pair of random variables $(X,Y)$, with $X\sim \mu$ and $Y\sim \nu$, such that for all $f\in\mathcal{F}$, $(f(X),f(Y))$ is a one-step submartingale. 

Let us remark that, for given two Radon measures $\mu,\nu$, the sets of $\mathcal{H}$-transports and $\mathcal{F}$-transports coincide and $\mu\prec_{\mathcal{F}}\nu$ if and only if $\mu\prec_{\mathcal{H}}\nu$.

The main results that of the paper are concerned with the structure of $\mathcal{F}$-transports.
Given two measures $\mu\prec_{\mathcal{F}}\nu$ we shall define $\mathcal{A}$-convex sets called \emph{irreducible components}, where $\mathcal{A}$ is the linear span of $\mathcal{G}$. These sets depend only on $\mu,\nu$, $\mathcal{F}$ and $\mathcal{A}$. Under the assumption that all $\mathcal{F}$-transports between $\mu$ and $\nu$ are \emph{local}, meaning that they enjoy a certain finite-dimensionality property, we shall demonstrate that the irreducible components are pairwise disjoint.\footnote{Let us remark that the locality assumption is not overly strict -- it is easy to provide examples of pairs of measures in the convex order on an infinite-dimensional linear space, for which any martingale transport is local.}
Moreover, any $\mathcal{F}$-transport between $\mu$ and $\nu$ has to occur within closures of these sets.

The partitioning shall be called \emph{irreducible} $\mathcal{F}$-\emph{convex paving}.

We shall also consider \emph{polar sets} with respect to all $\mathcal{F}$-transports between $\mu$ and $\nu$, i.e.,  Borel sets $U\subset\Omega\times\Omega$ which are of zero measure with respect to every $\mathcal{F}$-transport. Under the assumption of locality, we shall show that if for $\mu$-almost every $\omega_1\in\Omega$ the section
\begin{equation*}
U_{\omega_1}=\{\omega_2\in\Omega\mid (\omega_1,\omega_2)\in U\}
\end{equation*}
is contained in the corresponding irreducible component, then $U$ is a  polar set if and only if it  projects to sets of $\mu$-measure zero and of $\nu$-measure zero, respectively. This provides an affirmative answer to generalisation of a conjecture of Ob\l\'oj and Siorpaes \cite[Conjecture 1.3.]{Obloj2017}.

Among our results is also a generalisation of the Strassen's theorem to the setting of generalised convexity, which includes also the case of infinite-dimensional dual Banach spaces. 

Let us recall that the classical theorem of Strassen, see \cite{Strassen1965}, tells that for two Borel probability measures $\mu,\nu$ on a finite-dimensional linear space there exists an $\mathcal{F}$-transport between $\mu$ and $\nu$, where $\mathcal{F}$ is the cone of convex, lower semi-continuous functions, if and only if $\mu\prec_{\mathcal{F}}\nu$.
Our result extends this characterisation to infinite-dimensional spaces.

\subsection{Martingale transport in finite-dimensional spaces}

An example of the above-described situation is the case of martingale transport, where one takes $\mathcal{F}$ to be the cone of all convex functions on a finite-dimensional linear space, $\mathcal{G}$ to be the linear space of all affine functions and $p$ to be a norm on the space. In this setting, in the works of Ghoussoub, Kim, Lim \cite{Ghoussoub2019}, De March, Touzi \cite{Touzi2019} and Ob\l\'oj and Siorpaes \cite{Obloj2017}, it has been shown that given two Borel probability measures in convex order there exists the finest partitioning which constrains all martingale transports between the considered measures. The elements of the partitioning  have been called irreducible components and the partitioning itself has been called an irreducible convex paving.

It has been conjectured in  \cite[Conjecture 1.3.]{Obloj2017} that this partitioning allows for a characterisation of the polar sets with  respect to all martingale transports between the considered measures.  An answer to the conjecture has already been provided by \cite[Theorem 2.5.]{Touzi2019}.  This work provides an  affirmative resolution of a generalisation of this conjecture, which includes the case of local martingale transports in infinite-dimensional spaces. 

Let us also mention that the studies of the martingale transport problem in  the one-dimensional discrete setting have been conceived by Henry-Labord\`ere, Beiglb\"ock and Penkner \cite{Beiglbock2013}. The martingale transport problem served there as the dual of the problem of robust, model-free, superhedging of exotic derivatives in financial mathematics. 
It has been further developed by Beiglb\"{o}ck, Juillet \cite{Juillet2016}, who studied the supports of martingale optimal transports, in a similar way  to the work of Gangbo and McCann \cite{McCann1996}.

Further extensions of the Kantorovich duality were studied in \cite{Zaev2015}, where, see\cite[Theorem 4.3.]{Zaev2015}, the author reproves the result of \cite{Juillet2016}.
Let us also mention the related work of Beiglb\"ock, Nutz and Touzi \cite{Nutz2017}, where
a quasi-sure formulation of the dual problem is introduced. It allowed the authors to obtain a general duality result, with no duality gap, and prove the existence of dual optimisers. 

An extension of the martingale transport problem to the continuous-time setting has been investigated by Gallichon, Henry-Labord\`ere and Touzi in  \cite{Galichon2014}.

Further studies that investigated the quasi-sure duality and local structure of martingale optimal transports were performed in \cite{DeMarch2018} and in  \cite{DeMarch20182}.

An abstract version of the duality theorem has been proven by Ekren and Soner in \cite{Ekren2018}.
Works concerning duality in the setting of martingale transports in the Skorokhod space include \cite{Dolinsky2015} and \cite{Cheridito2021}. 

The duality in the martingale optimal transport problem remains an active area of research. In \cite{Backhoff2020} the authors develop a displacement interpolation for the martingale optimal transport, that would be an analogue of the McCann's interpolation \cite{McCann1995}, \cite{McCann1997} and the optimal Brenier map \cite{Brenier1991}. Further studies \cite{Backhoff2023}, \cite{Backhoff20232} in this direction investigate displacement convexity of Bass functional, along generalised geodesics. 

\subsection{Generalised convexity}\label{s:convex}

As we have mentioned for any $\mu\prec_{\mathcal{F}}\nu$, there exists the finest partitioning that constrains any $\mathcal{F}$-transport between $\mu$ and $\nu$. The elements of the partitioning are $\mathcal{A}$-convex sets, where $\mathcal{A}$ is the linear span of $\mathcal{G}$. 

Let $\mathcal{F}$ be a convex cone of functions on a topological space $\Omega$.
A  closed set $K\subset\Omega$ is said to be $\mathcal{F}$-convex  whenever it equals to its $\mathcal{F}$-convex hull, that  is
\begin{equation*}
\mathrm{clConv}_{\mathcal{F}}K=\{\omega\in\Omega\mid f(\omega)\leq \sup f(K)\text{ for all }f\in\mathcal{F}\}.
\end{equation*}
Let us observe that when $\mathcal{F}$ is the complete lattice cone generated by a linear space $\mathcal{A}$ then $K$ is $\mathcal{F}$-convex if and only if it is $\mathcal{A}$-convex. Moreover, if $\Omega$ is a locally convex topological vector space and $\mathcal{A}$ is the space of all continuous linear functionals on $\Omega$,  then the above-defined notion of convexity coincides with the usual notion, as follows from the Hahn--Banach theorem.

The notion of generalised convexity with respect to an arbitrary family of functions was introduced by Fan in \cite{Fan1963}. It allowed for a generalisation of the Krein--Milman theorem.
The notion has been developed and further generalised, see \cite{Dolecki1978}, \cite{Rubinov2000} and \cite{Singer1997}. These developments mainly were concerned with various dualities,  notions of subdifferentials and Fenchel transforms. These have found its application in the monopolist's problem by Figalli, Kim  and McCann \cite{Figalli2011} and by McCann and Zhang in \cite{McCann2019}.

We refer the reader to the book of H\"ormander \cite{Hormander2007} for an account comprising matter on notions of convexity with respect to: subharmonic functions, \cite[Chapter III]{Hormander2007}, plurisubharmonic  functions, \cite[Chapter  IV]{Hormander2007}. We note  also that the latter is equivalent to the notion of  holomorphic convexity of complex analysis, see e.g. \cite[Chapter 3]{Krantz2001}.

We refer the reader also to a paper \cite{Ciosmak2023} for a detailed discussion of generalised convexity that is relevant to the setting of our current developments. The paper is concerned with a generalisation of the Levi problem \cite{Levi1911} and the Cartan--Thullen theorem \cite{Cartan1932}. In particular, theorems in \cite{Ciosmak2023} characterise topological spaces $\Omega$, equipped with a linear space of continuous functions $\mathcal{A}$ on $\Omega$, such that there exists a proper function $p$ that belongs to the complete lattice cone $\mathcal{H}$ generated by $\mathcal{A}$. As it turns out, these spaces are precisely spaces that are complete with respect to $\mathcal{A}$. 

The notion of $\mathcal{F}$-convex sets is discussed in detail in Section  \ref{s:fconvex}.

\subsection{Assumptions}

Before we state our main results, let us describe our setting and state the standing assumptions.

\begin{assumption}\label{as:ass}
Let $\mathcal{G}$ be a convex cone of functions on a set $\Omega$ that contains constants.
Let $\mathcal{F}$ be the lattice cone generated by $\mathcal{G}$. We shall assume that there exists $p$ in the complete lattice cone $\mathcal{H}$ generated by $\mathcal{F}$, that is non-negative and proper with respect to $\tau(\mathcal{G})$. 
Moreover, we assume also that $\mathcal{G}$ separates points of $\Omega$, consists of functions of $(p+1)$-growth. 
We shall also denote by $\mathcal{A}$ the linear span of $\mathcal{G}$, i.e., $\mathcal{A}=\mathcal{G}-\mathcal{G}$.
\end{assumption}

 Let us remark that $\tau(\mathcal{G})$ denotes the weak topology on $\Omega$ generated by $\mathcal{G}$. We stress that we do not assume that $\Omega$ is metrisable, so that we will be able to include in our considerations dual infinite-dimensional spaces equipped with weak* topologies.
 
\begin{example}
    Let $\Omega$ be a separable dual Banach space and let $\mathcal{G}$ consist of weakly* continuous affine functionals. Then the assumptions are satisfied, with $p$ being the norm on $\Omega$, and $\mathcal{F}$ being the cone of lower semi-continuous, with respect to the weak* topology, convex functions on $\Omega$. The fact that $p\in\mathcal{F}$ follows by the Hahn--Banach theorem and the fact that $p$ is proper follows by the Banach--Alaoglu theorem.
\end{example}

\begin{assumption}\label{as:sep}
    We shall usually assume, except for the proof of the generalisation of the Strassen theorem, Theorem \ref{thm:mart}, that $\Omega$ is separable in $\tau(\mathcal{G})$.
\end{assumption}

In the above example, if $\Omega$ is separable, as a Banach space, then it is also separable in $\tau(\mathcal{G})=\tau(\mathcal{A})$.

\subsection{Strassen theorem}

The first result that we prove is a generalisation of the Strassen theorem \cite{Strassen1965} to our general setting, in which we do not assume that the underlying space is locally compact, unlike the previously known versions do, cf. \cite{Ciosmak20232}. 

\begin{theorem}\label{thm:strassenINTRO}
Suppose that $\mu,\nu$ are Radon, probability measures on $\sigma(\tau(\mathcal{G}))$ such that 
\begin{equation*}
\int_{\Omega}p\, d\mu<\infty\text{ and }\int_{\Omega}p\, d\nu<\infty
\end{equation*}
and that
\begin{equation}\label{eqn:majintro}
\int_{\Omega}f\,d\mu\leq\int_{\Omega} f\,d\nu
\end{equation}
for all $f\in\mathcal{F}$.

Then there exists a Radon measure $\pi$ on $\Omega\times\Omega$ with marginals $\mu$ and $\nu$ and such that for all $f\in\mathcal{F}$,
and non-negative, bounded, $\sigma(\tau(\mathcal{G}))$-measurable functions $h$ on $\Omega$
\begin{equation}\label{eqn:thefintro}
\int_{\Omega}h(\omega_1)f(\omega_1)\, d\mu(\omega_1)\leq\int_{\Omega\times\Omega}h(\omega_1)f(\omega_2)\, d\pi(\omega_1,\omega_2).
\end{equation}
Conversely, for a Radon measure $\pi$  on $\Omega\times\Omega$ for which (\ref{eqn:thefintro}) is satisfied, its marginals $\mu$, $\nu$ satisfy (\ref{eqn:majintro}).
\end{theorem}

Let us note that, see Proposition \ref{pro:submarti}, the condition (\ref{eqn:thefintro}) is equivalent to: if $(X,Y)$ is a random variable distributed according to $\pi$, then for $f\in\mathcal{F}$ the pair $(f(X),f(Y))$ is a one-step submartingale. In particular, for $f\in\mathcal{F}\cap(-\mathcal{F})$, the pair $(f(X),f(Y))$ is a one-step martingale, i.e.,
\begin{equation*}
\mathbb{E}(f(Y)\vert \sigma(X))=f(X).
\end{equation*}

All Radon measures $\pi$ that satisfy  the conclusion of Theorem \ref{thm:strassenINTRO} we  shall call $\mathcal{F}$-transports between $\mu$ and $\nu$. The set of such measures we shall denote by $\Gamma_{\mathcal{F}}(\mu,\nu)$.

The results discussed above are presented and proven in Section \ref{s:varstrassen} and Section \ref{s:existence}.

\subsection{Gelfand transform}

Before we will be able to state further results, we shall consider the Gelfand transform $\Phi\colon\Omega\to\mathcal{A}^*$, 
given by the formula
\begin{equation*}
\Phi(\omega)(a)=a(\omega)\text{  for all }\omega\in\Omega\text{  and }a\in\mathcal{A}.
\end{equation*}
Here $\mathcal{A}^*$ is the dual space to $\mathcal{A}$, equipped with the norm
\begin{equation*}
    \norm{a}_{\mathcal{A}}=\sup\frac{\abs{a}}{p}(\Omega)\text{ for }a\in\mathcal{A}.
\end{equation*}

The name is inspired by the Gelfand transform in the theory of Banach algebras.
Although we develop our theory for general spaces $\Omega$ equipped with a convex cone of functions $\mathcal{G}$ on $\Omega$, all this can be inferred from a theory for dual spaces with the corresponding convex cone being a subset of weakly* continuous affine functionals, which corresponds to submartingale transports, see Theorem \ref{thm:embed} and Remark \ref{rem:eui}.

The point of employing the embedding map $\Phi\colon\Omega\to\mathcal{A}^*$ is to introduce an analogue of a linear structure on $\Omega$. 

The  details are presented in Section \ref{s:gelfand}.

\subsection{Maximal disintegrations}

Assume that $\Omega$ is separable in $\tau(\mathcal{G})$.
In Section \ref{s:joint} we show that there always exists a joint support for $\Gamma_{\mathcal{F}}(\mu,\nu)$. However, Example \ref{exa:jointbad} demonstrates that the joint support does not suffice for characterisation of polar sets with respect to $\mathcal{F}$-transports.
Any $\mathcal{F}$-transport between $\mu$ and $\nu$ admits a disintegration with respect to its first marginal $\mu$, see Section \ref{s:maxx}.  The set of such disintegrations we shall denote by $\Lambda_{\mathcal{F}}(\mu,\nu)$

Moreover, in Section \ref{s:maximaldis} we show that for any $\mu\prec_{\mathcal{F}}\nu$ there exists a disintegration $\lambda\in\Lambda_{\mathcal{F}}(\mu,\nu)$, which has maximal supports, i.e., for any $\lambda'\in\Lambda_{\mathcal{F}}(\mu,\nu)$
\begin{equation*}
    \mathrm{supp}\lambda'(\omega,\cdot)\subset\mathrm{supp}\lambda(\omega,\cdot)\text{ for }\mu\text{-almost every }\omega\in\Omega.
\end{equation*}
It  is immediate that the supports of maximal disintegration are unique.

\subsection{Gleason parts and Harnack inequality}

For a convex set $K$ and $k\in K$ the Gleason part $G(k,K)$ of $k$ in $K$ is the set of all $k'\in K$ such that the line segment $[k,k']$ extends in $K$ in both directions.
For example, in the finite-dimensional setting, $G(k,K)$ is the relative interior of a face of $K$ that contains $k$ in its relative interior.

We refer the reader to the book of Bear \cite{Bear1970}, for an introduction to the topic, and to the book of Alfsen \cite[Chapter II, \textsection 5, p. 122]{Alfsen1971}.

The notion of Gleason part is related to the Harnack inequality. 
The Harnack inequality, see e.g. \cite{Evans2010},  says that if $\Omega\subset \mathbb{R}^n$ is an open connected set, $\omega_1,\omega_2\in\Omega$, then there exist $C>0$ such that for all a non-negative, harmonic functions $h$ on $\Omega$ there is
\begin{equation*}
\frac1Ch(\omega_1)\leq h(\omega_2)\leq C h(\omega_1).
\end{equation*}
If $\mathcal{H}(\Omega)$ denotes the space of harmonic functions on $\Omega$ and $\Phi\colon\Omega\to\mathcal{H}(\Omega)^*$ is the Gelfand transform, then the Harnack inequality states that for any  two points $\omega_1,\omega_2\in\Omega$ belong to a single Gleason part of the convex set of non-negative linear functionals in $\mathcal{H}(\Omega)^*$. 

For a detailed discussion, we refer the reader to Section \ref{s:interia}.

\subsection{Irreducible components and irreducible $\mathcal{F}$-convex paving}

The second main result demonstrates that for any two Radon probability measures $\mu\prec_{\mathcal{F}}\nu$ there exists the finest partitioning into $\mathcal{A}$-convex sets that constrains any $\mathcal{F}$-transport between $\mu$ and $\nu$. We suppose that Assumption \ref{as:ass} holds true.

Let $\lambda\in\Lambda_{\mathcal{F}}(\mu,\nu)$ be a maximal disintegration of $\Gamma_{\mathcal{F}}(\mu,\nu)$. 
For $\omega\in\Omega$ we shall define the irreducible component $\mathrm{irc}_{\mathcal{A}}(\mu,\nu)$ to be the relative interior of $\mathrm{clConv}_H\Phi(\mathrm{supp}\lambda(\omega,\cdot))\subset \mathcal{A}^*$. 

We see immediately that  the irreducible components depend only on $\mu,\nu$, $\mathcal{F}$ and $\mathcal{A}$.

Let us observe that the irreducible components are $\mathcal{A}$-convex, see Section \ref{s:fconvex}. Remark \ref{rem:fgaconvexity} tells us that they are also $\mathcal{F}$-convex.

 We shall below consider \emph{local} $\mathcal{F}$-transports. These are transports such that 
 \begin{equation*}
    \mathrm{clConv}_H\Phi(\mathrm{supp}\lambda(\omega,\cdot))\subset \mathcal{A}^*\text{ is finite-dimensional for }\mu\text{-almost every }\omega\in\Omega.
 \end{equation*}

For a subspace $\mathcal{B}\subset\mathcal{A}$ we also define $\Phi_{\mathcal{B}}\colon\Omega\to\mathcal{B}^*$ by the formula $\Phi_{\mathcal{B}}(\omega)(b)=b(\omega)$ for $b\in\mathcal{B}$ and $\omega\in\Omega$. The map $R_{\mathcal{B}}\colon\mathcal{A}^*\to\mathcal{B}^*$ is the restriction operator, that assigns to a functional on $\mathcal{A}$ its restriction to $\mathcal{B}$.

\begin{theorem}\label{thm:partitionINTRO}
Suppose that $\mu\prec_{\mathcal{F}}\nu$, that $\Omega$ is separable in $\tau(\mathcal{G})$ and that any $\mathcal{F}$-transport between $\mu$ and $\nu$ is local.
Then there exists a Borel measurable set $B\subset\Omega$ with $\mu(B)=1$ such that whenever $\omega_1,\omega_2\in B$ then 
\begin{equation*}
\mathrm{irc}_{\mathcal{A}}(\mu,\nu)(\omega_1)\cap\mathrm{irc}_{\mathcal{A}}(\mu,\nu)(\omega_2)=\emptyset
\end{equation*}
or
\begin{equation*}
\mathrm{irc}_{\mathcal{A}}(\mu,\nu)(\omega_1)=\mathrm{irc}_{\mathcal{A}}(\mu,\nu)(\omega_2).
\end{equation*}
These sets are $\mathcal{A}$-convex.
Moreover, for any $\lambda\in\Lambda_{\mathcal{F}}(\mu,\nu)$
\begin{equation*}
    \mathrm{supp}\lambda(\omega,\cdot)\subset \Phi^{-1}\big(\mathrm{cl}\mathrm{irc}_{\mathcal{A}}(\mu,\nu)(\omega)\big)\text{ for }\mu\text{-almost every }\omega\in\Omega.
\end{equation*}
The sets $\mathrm{irc}_{\mathcal{A}}(\mu,\nu)(\cdot)$ are the smallest convex sets that satisfy  the above condition.
Furthermore, for $\mathcal{B}=\mathcal{G}\cap(-\mathcal{G})$,
\begin{equation*}
    \Phi_{\mathcal{B}}(\omega)\in R_{\mathcal{B}}(\mathrm{irc}_{\mathcal{A}}(\mu,\nu)(\omega))\text{, }\mu\text{-almost every }\omega\in\Omega.
\end{equation*}
In particular, if $\mathcal{G}$ is a linear subspace, then
\begin{equation*}
    \omega\in \Phi^{-1}\big(\mathrm{irc}_{\mathcal{A}}(\mu,\nu)(\omega)\big)\text{ for }\mu\text{-almost every }\omega\in\Omega.
\end{equation*}
\end{theorem}
The above  theorem is proven in Section \ref{s:irreducible}, see Theorem \ref{thm:partition}.
In Section \ref{s:dispartition}, Theorem \ref{thm:partdis}, we show that the equivalence relation defined by the irreducible components is measurable, so that wee may disintegrate the space $\Omega$ into its equivalence classes. 

If $\mathcal{G}$ is a linear subspace, then  $\Phi(\omega)$ belongs to $\mathrm{clConv}_H\Phi(\mathrm{supp}\lambda(\omega,\cdot))$. If all $\mathcal{F}$-transports are local then
\begin{equation*}
    \mathrm{irc}_{\mathcal{A}}(\mu,\nu)(\omega)=G\big(\Phi(\omega),\mathrm{clConv}_H\Phi(\mathrm{supp}\lambda(\omega,\cdot))\big)\subset\mathcal{A}^*\text{  for }\mu\text{-almost every }\omega\in\Omega.
\end{equation*}
Relevant description of the  Gleason parts is presented in Lemma \ref{lem:harnack}, which shows that in the setting of interest to our theory, the Gleason parts are also characterised by a Harnack-type inequality.

We could have defined the components by the formula
\begin{equation*}
  G\big(\omega,  \mathrm{clConv}_{\mathcal{F}}\mathrm{supp}\lambda(\omega,\cdot)\big)=\Phi^{-1}\Big(G\Big(\Phi(\omega),\mathrm{clConv}_{G}\Phi\big(\mathrm{supp}\lambda(\omega,\cdot)\big)\Big)\Big)\text{ for }\omega\in\Omega,
\end{equation*}
see Section \ref{s:fconvex}, Remark \ref{rem:gleasonf}.
However, the fact that such components are pairwise disjoint is strictly weaker than the corresponding disjointedness of the irreducible components,  therefore we prefer the previous definition.

Let us remark that one could prove the disjointness of the above defined irreducible components without assuming locality. However, without this assumption, these components may be empty. 

\subsection{Polar sets}

The third result characterises the polar sets with respect to all $\mathcal{F}$-transports between $\mu$ and $\nu$,  under an additional assumption that the section of the considered polar set are relatively open in the corresponding component.

Let us recall that a Borel set $U$ is polar with respect to all $\mathcal{F}$-transports between $\mu$ and $\nu$ if 
\begin{equation*}
    \pi(U)=0\text{  for all }\pi\in\Gamma_{\mathcal{F}}(\mu,\nu).
\end{equation*}

\begin{theorem}\label{thm:polarinto}
Suppose that $\mu\prec_{\mathcal{F}}\nu$. Suppose that $\Omega$ is separable in $\tau(\mathcal{G})$  and that any $\mathcal{F}$-transport between $\mu$  and $\nu$ is local.

Suppose that $U\subset\Omega\times\Omega$ is a Borel set such that for $\mu $-almost every $\omega\in \Omega$
\begin{equation*}
\Phi(U_{\omega})\subset\mathrm{irc}_{\mathcal{A}}(\mu,\nu)(\omega).
\end{equation*}
Then $U$ is $\Gamma_{\mathcal{F}}(\mu,\nu)$-polar set if and only if there exist Borel sets $N_1,N_2\subset\Omega$ with
\begin{equation*}
\mu(N_1)=0,\nu(N_2)=0
\end{equation*}
and 
\begin{equation*}
U\subset (N_1\times\Omega)\cup (\Omega\times N_2).
\end{equation*} 
\end{theorem}

The theorem is proven in Section  \ref{s:polarity}, Theorem \ref{thm:polarfin}.
In the case when $\Omega$ is a finite-dimensional linear space, $\mathcal{F}$ is the cone of convex functions and $\mathcal{A}$ is the space of affine functions on $\Omega$ the above theorem was conjectured by Ob\l\'oj and Siorpaes in \cite[Conjecture 1.3.]{Obloj2017}. Our theorem shows that the conjecture holds true in a much greater generality.

\subsection{Fine structure of the irreducible components}\label{s:finei}

Our approach allows us to provide a fine description of the intersections of the irreducible components. This description is the content of Theorem \ref{thm:finestrucfin}.
It says that there exists Borel $B\subset\Omega$ with $\mu(B)=1$ such that if for some $\omega_1,\omega_2\in B$ such that $\mathrm{clConv}_{H}\Phi(\mathrm{supp}\lambda(\omega_i,\cdot))$, $i=1,2,$ are finite-dimensional  there exists
\begin{equation*}
    a^*\in \mathrm{clConv}_{H}\Phi(\mathrm{supp}\lambda(\omega_1,\cdot))\cap  \mathrm{clConv}_{H}\Phi(\mathrm{supp}\lambda(\omega_2,\cdot))
\end{equation*}
then the Gleason parts of $a^*$ in the sets $ \mathrm{clConv}_{H}\Phi(\mathrm{supp}\lambda(\omega_1,\cdot))$ and $ \mathrm{clConv}_{H}\Phi(\mathrm{supp}\lambda(\omega_2,\cdot))$ are  equal, under the additional assumption that there exist probability measures $\eta_1,\eta_2$, $C_1,C_2>0$
such that 
\begin{equation*}
    \eta_i\leq C_i\lambda(\omega_i,\cdot)\text{ for }i=1,2,
\end{equation*}
where $\lambda\in\Lambda_{\mathcal{F}}(\mu,\nu)$ is a maximal disintegration, and
\begin{equation*}
    a^*(a)=\int_{\Omega}a\, d\eta_i\text{ for all }a\in\mathcal{A}\text{ and }i=1,2.
\end{equation*}

\subsection{Choice of cone $\mathcal{G}$}

\begin{remark}
The assumption that $\mathcal{G}$ is a convex cone serves us merely to facilitate the exposition. The theory works as well for a general set $\mathcal{G}$ of functions on $\Omega$ -- it suffices to consider the conical hull  $\mathrm{Cone}\mathcal{G}$ in place of a general set of functions.
\end{remark}

\begin{remark}
One can take $\mathcal{G}=\mathcal{F}$ and obtain a  partitioning into $(\mathcal{F}-\mathcal{F})$-convex sets. Note however, e.g., that if $\mathcal{F}$ is the cone of convex, non-decreasing functions, then $\mathcal{A}=\mathcal{F}-\mathcal{F}$ is dense in the space of all continuous functions, as follows by the Stone--Weierstrass theorem. Therefore any closed set is $\mathcal{A}$-convex and any Gleason part is a singleton. The difficulty is to find a subcone $\mathcal{G}$ large enough that it yields a fine partitioning, but small enough that the components are non-trivial.
\end{remark}

\begin{remark}
Let us note that, in particular, we can handle the following situation. Let $\mathcal{F}$ be a lattice cone and let $\mu\prec_{\mathcal{F}}\nu$ be two Radon probability measures.  Let $\mathcal{G}\subset\mathcal{F}$ be a convex cone containing constants, not necessarily generating $\mathcal{F}$. Since any $\mathcal{F}$-transport is also $\mathcal{G}$-transport, we see that there is exists a partitioning of $\Omega$ into $\mathcal{A}$-convex sets, which yields constraints for all $\mathcal{F}$-transports between $\mu$ and $\nu$.
\end{remark}

\subsection{Balayages of measures}\label{s:balayage}

Let $\mathcal{F}$ be a lattice cone of Borel measurable functions on a topological space $\Omega$. For two Radon probability measures $\mu,\nu$ we write $\mu\prec_{\mathcal{F}}\nu$ whenever for any $f\in\mathcal{F}$ 
\begin{equation*}
\int_{\Omega}f\,d\mu\leq\int_{\Omega}f\, d\nu,
\end{equation*}
provided that the integrals are well-defined.

Let us  stress that there is no loss of generality, if we assume that $\mathcal{F}$ is a complete lattice cone, as follows by Remark \ref{rem:finite}, see Lemma \ref{lem:taulem}.

The above order on measures, related to the notion of balayage, originates in the works on potential theory of Poincar\'e. It was studied by Mokobodzki in \cite{Mokobodzki1984}, where he defines $\nu$ to be balayage of $\mu$ relative to $\mathcal{F}$ if $\mu\prec_{\mathcal{F}}\nu$.
We refer the reader also to \cite{Bowles2019}, \cite{Ghoussoub2023},  \cite{Kim2024} and to \cite{Ciosmak20232} for recent developments and investigations of balayage in the context of optimal transport. We refer also to \cite[Chapter XI]{Meyer1966} for an accessible introduction to the theory of balayages and for a proof of Strassen's theorem \cite[T 51, p. 245]{Meyer1966}.

\begin{example}\label{exa:balayage}
Let us enumerate some interesting orders induced by convex cones.
\begin{enumerate}
\item Let $\mathcal{F}$ be the cone of subharmonic functions on $\mathbb{R}^n$. Then, by a result of Skorokhod, two non-negative measures $\mu,\nu\in\mathcal{P}(\mathbb{R}^n)$ are in $\mathcal{F}$-order, i.e., $\mu\prec_{\mathcal{F}}\nu$, if and only if there exists a stopping time $\tau$ and a Brownian motion $(B_t)_{t\in [0,\infty)}$ such that $B_0$ is distributed according to $\mu$ and $B_{\tau}$ is distributed according to $\nu$.  We refer the reader for a thorough survey on the  topic \cite{Obloj2004} and to \cite{Palmer2019}, \cite{Beiglbock2017}, \cite{Kim2020} for recent studies on the topic.
\item If $\mathcal{F}$ is the cone of convex functions on $\mathbb{R}^n$, then the $\mathcal{F}$-order is the convex order on  the space of probability measures. By Strassen's theorem two probability measures $\mu,\nu$ are in convex order if and only if they are  respective marginal distributions of a one-step martingale. 
\item\label{i:optimal} If $\mathcal{F}$ is the tangent cone at a $1$-Lipschitz $v$ to the set of $1$-Lipschitz functions on a metric space $\Omega$, then two probability measures $\mu,\nu$ are in $\mathcal{F}$-order if and only if $v$ is a Kantorovich potential for the Monge--Kantorovich mass transport problem from $\mu$ to $\nu$. Indeed, $\mu\prec_{\mathcal{F}}\nu$ if and only if for any $1$-Lipschitz $u$
\begin{equation*}
\int_{\Omega} (v-u)\,d\mu\leq\int_{\Omega} (v-u)\, d\nu\text{, that is }\int_{\Omega} u\,d(\nu-\mu)\leq\int_{\Omega} v\, d(\nu-\mu).
\end{equation*}
We refer the reader to the recent work \cite{Ciosmak20232} that explores this observation.
\item Let $\mathcal{F}$ be the tangent cone at $v$ to the set of $1$-Lipschitz and convex functions. Then two probability measures $\mu,\nu\in\mathcal{P}(\mathbb{R}^n)$ are in $\mathcal{F}$-order if and only if $v$ is an optimal function in the dual problem to the weak optimal transport problem between $\mu$ and $\nu$, introduced in \cite{Gozlan2017}, see also \cite{Gozlan2018}, \cite{Backhoff2019}. 
\item Let $\mathcal{F}$ be the cone of convex and increasing functions on $\mathbb{R}$. This ordering has been studied in \cite{Nutz2018}. In this setting a version of irreducible convex paving has also been obtained, \cite[Proposition 3.4]{Nutz2018}.
\end{enumerate}
\end{example}

\subsection{Examples}

Let us now shortly discuss two examples of our theory. We refer to the accompanying paper \cite{Ciosmak2024} for a thorough discussion of applications.

Let  us note  that  our developments  apply to any  pair of  measures $\mu,\nu$ that are in order with respect  to a lattice cone generated by a convex cone of functions containing constants.  Nevertheless, it is of interest to present some concrete examples.

\subsubsection{Martingale transports in infinite-dimensional spaces}

The basic example is the martingale transport. Our theory allows for a treatment of the infinite-dimensional case, which generalises results  of \cite{Juillet2016}, \cite{Nutz2017}, \cite{Touzi2019}, \cite{Obloj2017}. Specifically, suppose that $X^*$ is a dual, separable Banach space. Let $\mathcal{A}$ denote the space of all affine weakly* continuous functionals on $X^*$ and $\mathcal{F}$ be the corresponding lattice cone of generated by $\mathcal{A}$. It consists of convex   and continuous, with respect to the weak* topology on $X^*$, functions on  $X^*$.

Suppose we have two Radon probability measures $\mu,\nu$ on $X^*$ with finite first moments
such that $\mu\prec_{\mathcal{F}}\nu$.

 Theorem \ref{thm:strassenINTRO} shows that there exists a Radon probability measure $\pi$ on the space $X^*\times X^*$, with marginals $\mu$ and $\nu$, that is a distribution of a one-step martingale.

Suppose now that any $\mathcal{F}$-transport between $\mu$ and $\nu$ is local, see Definition \ref{def:local}.
Theorem \ref{thm:partitionINTRO} shows that there exists a partitioning of $X^*$, up to a set of  $\mu$-measure zero, into relatively open, convex and finite-dimensional sets $\mathrm{irc}_{\mathcal{A}}(\mu,\nu)(\cdot)$. Moreover, for any $\mathcal{F}$-transport and for any its disintegration $\lambda\in\Lambda_{\mathcal{F}}(\mu,\nu)$, 
\begin{equation*}
    \mathrm{supp}\lambda(x^*,\cdot)\subset \mathrm{cl}\,\mathrm{irc}_{\mathcal{A}}(\mu,\nu)(x^*)\text{ for }\mu\text{-almost every  }x^*\in X^*.
\end{equation*}
Moreover
\begin{equation*}
 x^*\in\mathrm{irc}_{\mathcal{A}}(\mu,\nu)(x^*)\text{ for }\mu\text{-almost every  }x^*\in X^*.
\end{equation*}
Furthermore, a description of fine structure of the intersections of the irreducible components is available, see Section \ref{s:finei}.

Theorem \ref{thm:polarinto} provides a characterisation of subsets of $X^*\times X^*$ that are polar with respect of all $\mathcal{F}$-transports. Suppose that all $\mathcal{F}$-transports between  $\mu$ and $\nu$ are local.  Let $U\subset X^*\times X^*$ be a Borel set such that 
\begin{equation*}
    U_{x^*}\subset \mathrm{irc}_{\mathcal{A}}(\mu,\nu)(x^*)\text{ for }\mu\text{-almost every }x^*\in X^*
\end{equation*}
Then $U$ is polar with respect to all $\mathcal{F}$-transports, if and only if there exist Borel sets $N_1,N_2\subset X^*$ such that 
\begin{equation*}
    \mu(N_1)=0,\nu(N_2)=0\text{  and }U\subset N_1\times X^*\cup X^*\times N_2.
\end{equation*}
In other words, $U$ is a polar set if and only if it is polar with respect to all transports between $\mu$ and $\nu$.
This provides an affirmative answer to a generalisation of a conjecture of Ob\l\'oj and Siorpaes of \cite{Obloj2017} to the infinite-dimensional setting.  

We shall present this example in more detail in \cite{Ciosmak2024}.

\subsubsection{Supermartingale transports}

Another example concerns  supermartingale transport. That is we let $\mathcal{H}$ be the cone  of lower  semi-continuous functions on $\mathbb{R}$ that are convex and non-increasing with respect to each co-ordinate. 

Two measures $\mu,\nu\in\mathcal{P}(\mathbb{R})$ with finite first moments are in $\mathcal{H}$-order if and only if they are distributions of the respective marginals of a one-step supermartingale, i.e., of a random variable $(X,Y)$ such that
\begin{equation*}
X\geq \mathbb{E}(Y\vert\sigma(X)).
\end{equation*}

In the one-dimensional case this problem has been studied in \cite{Nutz2018}. In that paper \cite[Proposition 3.4.]{Nutz2018}, it is proven that given two Borel probability measures $\mu\prec_{\mathcal{H}}\nu$ on $\mathbb{R}$, with finite first moments, there exists a partitioning of $\mathbb{R}$ into open intervals such that any $\mathcal{H}$-transport between $\mu$ and $\nu$ has to occur within their closures. Moreover, a characterisation of polar sets for all $\mathcal{H}$-transports has been provided. 

In  \cite{Ciosmak2024} we show how Theorem \ref{thm:partitionINTRO} and Theorem \ref{thm:polarinto}, can be applied to extend these results to the multi-dimensional setting.
Specifically, taking $\mathcal{G}\subset\mathcal{H}$ to be the cone of affine non-increasing functions will yield a partitioning into convex sets. Since the cone $\mathcal{G}$ generates $\mathcal{H}$ we also obtain a characterisation of polar sets.

More generally, we may consider a dual, separable Banach space $X^*$. Let $(x_i)_{i=1}^{\infty}$ be a linearly dense subset of $X$. Let $\mathcal{F}$ be the lattice cone of functions on $X^*$ generated by this set. If for a pair $\mu\prec_{\mathcal{F}}\nu$ any $\mathcal{F}$-transport is local, we can show existence of the finest partitioning  into convex sets and characterise polar sets, in an analogous way to the finite-dimensional setting.

A similar treatment can be applied to the submartingale transport.
The submartingale transport is relevant to the studies of the monopolist's problem investigated  by Rochet and Chon\'e \cite{Rochet1998} and, more recently, by Figalli, Kim, McCann \cite{Figalli2011} and McCann, Zhang \cite{McCann2019}. 

\subsubsection{Optimal transport localisation}

One of the aims of this paper is to develop a general localisation framework that would comprise both the irreducible convex paving pertaining to martingale transports and the localisation of the optimal transport theory, which we shall describe in detail in \cite{Ciosmak2024}.

The localisation in the context of Riemannian geometry has been developed by Klartag \cite{Klartag2017}, in the context of Finsler geometry by Ohta \cite{Ohta2018}, and by Cavalletti and Mondino \cite{Mondino2017}, \cite{Cavalletti2017} to essentially non-branching metric measure spaces that satisfy curvature-dimension conditions.

Let us note that it has been conjectured in \cite{Klartag2017} that the method can be generalised to multi-dimensional setting employing optimal transport of vector measures. The optimal transport of vector measures has been introduced in \cite{Ciosmak20211}, where one part of Klartag's conjecture has been refuted. Another version of optimal transport of vector measures has been proposed in \cite{Ciosmak20212}. The other part of the conjecture of Klartag has been partially resolved in the affirmative in \cite{Ciosmak2021}. The method that was proposed in \cite{Klartag2017} for a resolution has been shown in \cite{Ciosmak20213} to contain a gap.

\subsubsection{Martingale problems}

 Our framework includes also the solutions of martingale problems. Suppose that $L$ is a linear operator on the subspace $\mathcal{D}$ of the space of Borel functions, $\mu$ is a Borel probability measure and $(X_t)_{t\in [0,1]}$ is a solution of the martingale problem for $(L,\mu)$, i.e., for any $f\in\mathcal{D}$
\begin{equation*}
    f(X_t)-\int_0^t(Lf)(X_s)\,ds
\end{equation*}
is a well-defined martingale. Let $\mathcal{A}$ denote the kernel of $L$ and let $\nu$ be the distribution of $X_1$. Then there is a partitioning of the underlying space into $\mathcal{A}$-convex sets $\Phi^{-1}(\mathrm{apirc}_{\mathcal{G}_0}(\mu,\nu)(\cdot))$ such that for all $t\in [0,1]$
\begin{equation*}
    X_t\in \Phi^{-1}\big(\mathrm{cl}(\mathrm{apirc}_{\mathcal{G}_0}(\mu,\nu)(X_0))\big)\text{ almost surely.}
\end{equation*}
We can treat similarly e.g. holomorphic martingales.

We refer the reader to \cite{Ciosmak2024} for a detailed discussion of these ideas. It is shown there that the obtained partitioning can be non-trivial.

\subsection{Methods of proofs}

Let us now discuss the most interesting methods of the proofs of the main results. 

\subsubsection{Previous approaches}

The results concerning the irreducible convex paving that were obtained in \cite{Touzi2019}, \cite{Ghoussoub2019} and \cite{Obloj2017} relied on the dual point of view.
In \cite[Proposition 4.1., p. 20, Proposition 6.3., p. 27]{Touzi2019} the authors use a duality result to prove the existence of partitioning.

The basic observation related to this point of view is as follows. Suppose that $\mu,\nu$ are two Borel probability measures on $\mathbb{R}^n$ with finite first moments, in convex order.
That is, for any convex, lower semi-continuous function $f$ on $\mathbb{R}^n$, we have
\begin{equation}\label{eqn:ef}
\int_{\mathbb{R}^n}f\,d \mu\leq\int_{\mathbb{R}^n}f\, d\nu.
\end{equation}
If for some such $f$ we had an equality in the above inequality, then we could readily show that any martingale transport between $\mu$ and $\nu$ is constrained by the faces of $f$, i.e., maximal sets on which $f$ is affine, cf. Lemma \ref{lem:faces}.

Our method for the description of the irreducible $\mathcal{F}$-convex paving is, to an extent, elementary and relies merely on the primal problem concerning the martingale transports. 

\subsubsection{Irreducible components and irreducible $\mathcal{F}$-convex paving}\label{ss:components}

The novelty of our approach is that to prove that the irreducible components form a partitioning we use the idea presented below and the result of Kellerer \cite[Proposition 3.3., p. 424 and Proposition 3.5., p. 425]{Kellerer1984} that characterises polar sets with respect to all transports between any two given measures; see also \cite[Proposition 2.1.]{Beiglbock2009}. No other duality result in employed in our considerations.

The proof of Theorem \ref{thm:partitionINTRO}, presented in full in Section \ref{s:irreducible}, is based on the following observation.

Let $\mu\prec_{\mathcal{F}}\nu$. One can choose a maximal disintegration  $\lambda\in\Lambda_{\mathcal{F}}(\mu,\nu)$, i.e., it is that for any other $\lambda'\in\Lambda_{\mathcal{F}}(\mu,\nu)$
\begin{equation*}
\mathrm{supp}\lambda'(\omega,\cdot)\subset \mathrm{supp}\lambda(\omega,\cdot) \text{ for }\mu\text{-almost every }\omega\in\Omega.
\end{equation*}
We refer to Section \ref{s:maximaldis} for the proof.

As we have already specified, an irreducible component $\mathrm{irc}_{\mathcal{A}}(\mu,\nu)(\omega)$ for $\omega\in\Omega$ is the relative interior of the convex set $\mathrm{clConv}_{H}\Phi(\mathrm{supp}\lambda(\omega,\cdot))$.

If the irreducible component  is non-empty, then it equals to the  set of functionals $a^*\in\mathcal{A}^*$ for which there exists a function $h\in L^{\infty}(\Omega,\lambda(\omega,\cdot))$, bounded from below by a positive constant, such that for all $a\in\mathcal{A}$
\begin{equation*}
a^*(a)=\int_{\Omega}ah\, d\lambda_0(\omega,\cdot)
\end{equation*}
see Theorem \ref{thm:finite}.  Let us assume that $\mu$ is countably supported.
If two irreducible components of $\omega_1,\omega_2$ intersected, then we would be able to find  $h_1\in L^{\infty}(\Omega,\lambda(\omega_1,\cdot))$ and $h_2\in L^{\infty}(\Omega,\lambda_0(\omega,\cdot))$ such that 
\begin{equation*}
\int_{\Omega}ah_1\, d\lambda(\omega_1,\cdot)=\int_{\Omega}ah_2\, d\lambda(\omega_2,\cdot)\text{ for all }a\in\mathcal{A}.
\end{equation*}
That would allow for a modification of the original kernel $\lambda$ in the following way
\begin{equation*}
d\lambda'(\omega_1,\cdot)=(1-h_1)\,d\lambda(\omega_1,\cdot)+h_2\,d\lambda(\omega_2,\cdot),
\end{equation*}
and
\begin{equation*}
d\lambda'(\omega_2,\cdot)=(1-h_2)\,d\lambda(\omega_2,\cdot)+h_1\,d\lambda(\omega_1,\cdot),
\end{equation*}
where we assumed, without loss of generality, that $h_1,h_2$ are bounded by one.
The new kernel $\lambda'$ would give rise to an $\mathcal{F}$-transport with the same marginals as the original one. Therefore, as $h_1,h_2$ are positive
\begin{equation*}
\mathrm{supp}\lambda(\omega_1,\cdot)=\mathrm{supp}\lambda(\omega_2,\cdot).
\end{equation*}
Thus also 
\begin{equation*}
\mathrm{clConv}_{H}\Phi(\mathrm{supp}\lambda(\omega_1,\cdot))=\mathrm{clConv}_{H}\Phi(\mathrm{supp}\lambda(\omega_2,\cdot)).
\end{equation*}
It follows that the corresponding parts have to be equal as well.

This shows that if the components intersect, they are equal. This is to say, the components indeed partition the underlying space.

To treat the general case, we need to modify the kernel in an alternative way. The idea is to take all points whose components intersect the considered component and, using a symmetric transport between $\mu$ and $\mu$, obtain a modified kernel. In this way we infer that the set of pairs of points with non-equal and intersecting irreducible components is a polar set for all transports between $\mu$ and $\mu$.
The aforementioned result of Kellerer \cite{Kellerer1984} completes the proof.

\subsection{Polar sets}

Let us now discuss the proof of Theorem \ref{thm:polarinto}, presented in full in Section \ref{s:polarity}.
The idea is to show that if a Borel set $U\subset\Omega\times\Omega$, that is polar with respect to all $\mathcal{F}$-transports between $\mu$ and $\nu$, has sections in the preimages of the irreducible components, then it is polar with respect to all transports between $\mu$ and $\nu$. 

The main difficulty in extending the one-dimensional results of \cite{Nutz2017} to multi-dimensional setting in \cite{Touzi2019}, was that, as the authors of \cite{Touzi2019} write, that the convex functions in multi-dimensional setting do not have a simple generating family.

We circumvent these problems by considering bounded convex functions on polytopes, which enjoy certain compactness properties.

The main, and worth noticing, change in the approach, as compared with the one presented in \cite{Nutz2017}, is that we do not use the potential functions of $\mu$ and of $\nu$. Instead, we observe that if $S$ is a polytope in a finite-dimensional linear space, and if $V$ is its finite set of vertices, then the space $\mathcal{B}$ of convex, non-negative function $b$ such that $\max b(V)=1$ is uniformly Lipschitz on compact subsets of relative interia of faces of $S$. 
Moreover, a pointwise limit of elements of $\mathcal{B}$ belongs to $\mathcal{B}$, thanks to the finiteness of the set of vertices. If we moreover assume that any $b\in\mathcal{B}$ vanishes at some point of the relative interior of $S$, then no function in $\mathcal{B}$ is affine. 
This compactness serves us to obtain the necessary room in order to construct, for any given transport $\gamma$ between $\mu$ and $\nu$, an $\mathcal{F}$-transport between  $\mu$ and $\nu$ with dominates disintegrations of $\gamma$, in the sense of absolute continuity.
Theorem \ref{thm:mart} is used to finish off the proof. Note that we do not apply Theorem \ref{thm:mart} on  each irreducible component separately, since we do not know whether the constructed marginals, on each of the components, are non-negative measures.

We see therefore that our usage of duality is limited to Theorem  \ref{thm:mart} and the result of Kellerer \cite{Kellerer1984}, unlike in \cite{Touzi2019}.

\subsection{Issues in the infinite-dimensional case}

As we have seen, in  Theorem \ref{thm:partitionINTRO} and in Theorem \ref{thm:polarinto} we assume that all $\mathcal{F}$-transports between $\mu$ and $\nu$ are local.  

This assumption is used, for the sake of Theorem \ref{thm:partitionINTRO}, in two ways. The first is  that any finite-dimensional convex set has non-empty interior and we may therefore use Theorem \ref{thm:interior}. The second is \cite[Proposition  1.2.9, p. 18]{Zalinescu2002}.
Theorem \ref{thm:polarinto} needs the finite-dimensionality assumption for two additional reasons. The first is its use in \ref{lem:faces} and the second is the exhaustability of finite-dimensional convex sets by convex polygons, cf. proof of Theorem \ref{thm:polarfin}.

As shown in Example \ref{exa:interior}, the relative interior of a convex set is can be empty. If $\mathcal{G}$ is a linear subspace, a natural idea is to consider the Gleason part of $\Phi(\omega)$ in $\mathrm{clConv}_H\Phi(\mathrm{supp}\lambda(\omega,\cdot))$, $\omega\in\Omega$, instead of the relative interior. The Gleason parts have the advantage of being non-empty.

However, the characterisation of the Gleason parts, Theorem \ref{thm:finite}, does not carry  over to the infinite-dimensional case, which does not allow to an extension of our results to this case. This is discussed in Section \ref{s:issues}, Theorem \ref{thm:nonresult}.

A resolution of that obstacle can be achieved by the procedure of finite-dimensional approximations, see \cite{Ciosmak2024}. 

 \addtocontents{toc}{\setcounter{tocdepth}{10}}
 
\section{Preliminaries}\label{s:prelim}

Let us recall several definitions, that we shall need in our developments.

\begin{definition}\label{def:proper}
Let $\Omega$, $Z$ be topological spaces. A map $p\colon \Omega\to Z$ is said to be proper whenever preimages $p^{-1}(Z)$ of compact sets in $Z$ are compact sets in $\Omega$.
\end{definition}

\begin{definition}
Let $\Omega$ be a set and let $\mathcal{C}$ be a set of functions on $\Omega$. The coarsest topology on $\Omega$ with respect to which all functions in $\mathcal{C}$ are continuous we shall call the topology generated by $\mathcal{C}$ and denote by $\tau(\mathcal{C})$.
\end{definition}

\begin{definition}\label{def:growth}
Let $\Omega$ be a set and let $p$ be a function on $\Omega$ bounded from below by a positive constant. We shall say that a function $a$ on $\Omega$ is of $p$-growth if there exists $C>0$ such that for all $\omega\in\Omega$ there is
\begin{equation*}
\abs{a(\omega)}\leq Cp(\omega).
\end{equation*}
We denote by $\mathcal{D}_p(\Omega)$ the space of all functions of $p$-growth, with norm 
\begin{equation*}
\norm{a}_{\mathcal{D}_p(\Omega)}=\sup\bigg\{\frac{\abs{a(\omega)}}{p(\omega)}\mid \omega\in\Omega\bigg\}.
\end{equation*}
If $\Omega$ is a topological space, then we let $\mathcal{C}_p(\Omega)$ denote the space of all continuous functions of $p$-growth.
We let $\mathcal{C}(\Omega)$ denote the space of all continuous functions on $\Omega$ and $\mathcal{M}_p(\Omega)$ denote the space of signed Radon measures $\mu$ on $\Omega$ with 
\begin{equation*}
\int_{\Omega}p\, d\abs{\mu}<\infty.
\end{equation*}
For $\mu\in  \mathcal{M}_p(\Omega)$ we  set 
\begin{equation*}
\norm{\mu}_{\mathcal{M}_p(\Omega)}=\int_{\Omega}p\, d\abs{\mu}.
\end{equation*}
We let $\mathcal{P}_p(\Omega)$ denote the subset of $\mathcal{M}_p(\Omega)$ consisting of probability measures.
\end{definition}

\begin{definition}\label{def:baire}
Let $\Omega$ be a topological space, with topology $\tau$. Then:
\begin{enumerate}
\item $\sigma(\tau)$ is called the Borel $\sigma$-algebra on $\Omega$,
\item $\sigma(\mathcal{C}(\Omega))$ is called the Baire $\sigma$-algebra on $\Omega$,
\item a measure $\mu$ on $\sigma(\tau)$ is called a Radon measure if for any set $B\in\sigma(\tau)$ and $\epsilon>0$ there exists a compact set $K\subset B$ such that $\abs{\mu}(B\setminus K)<\epsilon$,
\item a measure $\mu$ on an algebra $\Theta$ is called regular, if for any $A\in\Theta$ and any $\epsilon>0$ there exists closed set $F\subset A$ such that $A\setminus F\in\Theta$ and  $\abs{\mu}(A\setminus F)<\epsilon$,
\item\label{i:tight} a measure $\mu$ on an algebra $\Theta$ is tight if for any $\epsilon>0$, there exists a compact set $K$ such that $\abs{\mu}(A)<\epsilon$ for any $A\in\Theta$ such that $A\cap K=\emptyset$.
\end{enumerate}
\end{definition}

\begin{definition}\label{def:prelim}
Let $\Omega$ be a set. Let $\mathcal{A}$ be a set of functions on $\Omega$ with values in $(-\infty,\infty]$. A set $\mathcal{K}$ of functions on $\Omega$ is said to be:
\begin{enumerate}
\item a convex cone whenever $\alpha f+\beta g \in\mathcal{K}$ for any $f,g\in\mathcal{K}$ and any numbers $\alpha,\beta\geq 0$;
\item stable under maxima provided that the maximum $f\vee g$ of any two functions   $f,g\in\mathcal{K}$ belongs to $\mathcal{K}$;
\item stable under suprema provided that the supremum $\sup\{f_i\mid i \in I\}$  of any family $(f_i)_{i\in I}\subset\mathcal{K}$ belongs to $\mathcal{K}$;
\item a vector lattice whenever it is a linear space such that is stable under maxima;
\item the vector lattice generated by a set $\mathcal{A}$ of functions on $\Omega$ if it is the smallest vector lattice containing $\mathcal{A}$;
\item a lattice cone whenever it is a convex cone stable under maxima;
\item the lattice cone generated by $\mathcal{A}$ whenever it is the smallest lattice cone of functions on $\Omega$ containing $\mathcal{A}$;
\item a complete lattice cone whenever it is a convex cone stable under suprema;
\item the complete lattice cone generated by $\mathcal{A}$ whenever it is the smallest complete lattice cone of functions on $\Omega$ containing $\mathcal{A}$;
\item separating points of $\Omega$ whenever for any two distinct $\omega_1,\omega_2\in\Omega$, there exist $f,g\in\mathcal{K}$ such that $f(\omega_1)\neq g(\omega_2)$.
\end{enumerate}
\end{definition}

Let us stress that we allow the functions to take value $+\infty$. Note however that this is not allowed for linear subspaces.

We refer to \cite{Lusky1981} for a study of the above notions in contexts related to this work.

\begin{lemma}\label{lem:supgen}
Suppose that $\mathcal{G}$ is a convex cone  of functions on $\Omega$. Let $\mathcal{H}$ denote the complete lattice cone generated by $\mathcal{G}$. Then for any $f\in\mathcal{H}$ there exists a family $(g_i)_{i\in I}\subset\mathcal{G}$ such that
\begin{equation}\label{eqn:forma}
f(\omega)=\sup\{g_i(\omega)\mid i\in I\}\text{ for all }\omega\in \Omega.
\end{equation}
\end{lemma}
\begin{proof}
Let $\mathcal{L}$ be the set of functions of the form (\ref{eqn:forma}).  We shall show that $\mathcal{L}$ is a convex cone. Let $\alpha,\beta\geq 0$, $f,g,\in\mathcal{L}$. Let $(g_i)_{i\in I}, (g_j)_{j\in J}\subset\mathcal{A}$ be families of functions in $\mathcal{G}$ corresponding to $f$ and $g$ respectively.
Then
\begin{equation*}
\alpha f+\beta g=\sup\{\alpha g_i+\beta g_j\mid j\in J, i\in I\}\in\mathcal{L}.
\end{equation*}
Clearly, $\mathcal{H}\supset\mathcal{L}$. Since $\mathcal{L}$ is a complete lattice cone, $\mathcal{H}=\mathcal{L}$.
\end{proof}

\begin{lemma}\label{lem:supgenfin}
Suppose that $\mathcal{G}$ is a convex cone  of functions on $\Omega$. Let $\mathcal{F}$ denote the lattice cone generated by $\mathcal{G}$. Then for any $f\in\mathcal{F}$ there exists a finite family $(g_i)_{i=1}^{n}\subset\mathcal{G}$ such that
\begin{equation*}
f(\omega)=\sup\{g_i(\omega)\mid i=1,\dotsc,n\}\text{ for all }\omega\in \Omega.
\end{equation*}
\end{lemma}
\begin{proof}
The proof follows analogous lines to the lines of the proof of Lemma \ref{lem:supgen}.
\end{proof}

\begin{corollary}\label{col:integral}
    Suppose that any function in  $\mathcal{G}$ is of $p$-growth. Then so is any function in $\mathcal{F}$.
\end{corollary}

Let us recall the following theorem of Daniell, see \cite[Theorem 7.8.1., p. 99]{Bogachev20072}.

\begin{theorem}\label{thm:daniell}
Let $\mathcal{B}$ be a vector lattice of functions on a set $\Omega$ such that $1\in\mathcal{B}$. Suppose that $\lambda$ is a non-negative linear functional on $\mathcal{B}$ such that $\lambda(1)=1$ and such that for any sequence of functions $(f_n)_{n=1}^{\infty}\subset\mathcal{B}$ monotonically decreasing to zero, there is $\lim_{n\to\infty}\lambda(f_n)=0 $.
Then there exists a unique probability measure $\mu$ on $\sigma(\mathcal{B})$ such that each function in $\mathcal{B}$ is integrable with respect to $\mu$ and for all $b\in\mathcal{B}$
\begin{equation*}
\lambda(b)=\int_{\Omega}b\,d\mu.
\end{equation*}
\end{theorem}

\subsection{Ordering with respect to a convex cone}

Let is define the ordering with respect to a convex cone of functions.
We have already discussed this notion in Section \ref{s:balayage}.

 \begin{definition}
Suppose that $\Omega$ is a topological space. Let $\mathcal{F}$ be a convex cone of Borel measurable functions on $\Omega$. Let $\mu,\nu\in\mathcal{P}(\Omega)$ be two Borel probability measures on $\Omega$. We shall say that $\mu$ and $\nu$ are in $\mathcal{F}$-order whenever for any $f\in\mathcal{F}$ that is integrable with respect to $\mu$ and $\nu$ there is
\begin{equation*}
\int_{\Omega}f\,d\mu\leq\int_{\Omega}f\,d\nu.
\end{equation*}
If the above condition is satisfied, we shall write $\mu\prec_{\mathcal{F}}\nu$.
\end{definition}

\begin{remark}
Let $\Omega$ be a topological space. Given a family of pairs of Borel probability measures $\{(\delta_{\omega},\eta_{\omega})\mid \omega\in \Omega\}$ we may define a convex cone $\mathcal{F}$ of  measurable functions $f\colon\Omega\to\mathbb{R}$ such that 
\begin{equation*}
f(\omega)\leq \int_{\Omega} f\, d\eta_{\omega}\text{ for all }\omega\in\Omega.
\end{equation*}
It was shown by Choquet and Deny that any lattice cone of continuous functions that contains constants can be obtained in this way, see \cite[Theorem 47]{Meyer1966}.
\end{remark}

\subsection{De la Vall\'ee-Poussin theorem}\label{s:poussin}

It turns out that it is of vital importance for applications to include in our considerations cones consisting of unbounded functions. An instance of such a situation is the case of martingale transport, when one considers two Borel probability measures $\mu,\nu$ on a finite-dimensional linear space in convex order, that is $\mu\prec_{\mathcal{F}}\nu$, where $\mathcal{F}$ is the cone of lower semi-continuous, convex functions. Inspired by an idea of \cite{Nutz2017} we shall employ de la Vall\'ee-Poussin theorem in the following version.  We provide a short proof for completeness, based on \cite{Chafai2014}, but see also \cite[Theorem 4.5.9., p. 272]{Bogachev20071}.

\begin{theorem}\label{thm:poussin}
Suppose that $\Omega$ is a topological space, $p\colon \Omega\to\mathbb{R}$ is a non-negative measurable function and that $\mu,\nu\in\mathcal{M}(\Omega)$ are non-negative measures such that 
\begin{equation*}
\int_{\Omega}p\,d\mu<\infty\text{ and }\int_{\Omega}p\,d\nu<\infty.
\end{equation*}
 Then there exists a continuous, convex and increasing function $\xi\colon \mathbb{R} \to[1,\infty)$ with superlinear growth, i.e.,
\begin{equation*}
\lim_{t\to\infty}\frac{\xi(t)}{t}=\infty
\end{equation*}
such that
\begin{equation*}
\int_{\Omega}\xi(p)\,d\mu<\infty\text{ and }\int_{\Omega}\xi(p)\,d\nu<\infty,
\end{equation*}
and $\xi(t)\geq t$ for all $t\geq 0$.
\end{theorem}
\begin{proof}
For $m\in\mathbb{N}$ pick $t_m\in\mathbb{R}$ such that 
\begin{equation*}
\int_{\Omega}p\mathbf{1}_{p^{-1}([t_m,\infty))}\,d\mu<\frac1{2^m}\text{ and }\int_{\Omega}p\mathbf{1}_{p^{-1}([t_m,\infty))}\,d\nu<\frac1{2^m}.
\end{equation*}
Define for $t\geq 0$
\begin{equation*}
\xi(t)=\sum_{m=1}^{\infty}(t-t_m)_+.
\end{equation*}
Then $\xi$ is continuous, convex, increasing and
\begin{equation*}
\lim_{t\to\infty}\frac{\xi(t)}{t}=\lim_{t\to\infty}\sum_{m=1}^{\infty}\bigg(1-\frac{t_m}{t}\bigg)_+=\infty.
\end{equation*}
Moreover, by construction,
\begin{equation*}
\int_{\Omega}\xi(p)\,d\mu\leq \sum_{i=1}^{\infty}\int_{\Omega}p\mathbf{1}_{p^{-1}([t_m,\infty))}\, d\mu\leq 1.
\end{equation*}
Similarly,
\begin{equation*}
\int_{\Omega}\xi(p)\,d\nu\leq \sum_{i=1}^{\infty}\int_{\Omega}p\mathbf{1}_{p^{-1}([t_m,\infty))}\, d\nu\leq 1.
\end{equation*}
Adding a function $t\mapsto t+1$ to $\xi$ yields a desired function.
\end{proof}

\section{Strassen's theorem for general spaces}\label{s:varstrassen}

In this section, we shall present a general version of Strassen's theorem, without any assumptions on the underlying spaces.

We shall first prove the following theorem, which is a  generalisation of  the Strassen's theorem (cf. \cite[Theorem T.51, p. 244]{Meyer1966}, see also \cite[Theorem 2.1]{Ciosmak20232}). The need for the generalisation is caused by our need to work with Banach spaces.

We refer the reader to \cite{Lusky1981} for a thorough study of various generalisations of the Strassen's theorem and related concepts concerning convex cones and balayages.

\begin{theorem}\label{thm:varstrassen}
Let $X$ be a Banach space, let $(\Omega,\Sigma)$ be a measurable space and let $\mu$ be a probability measure on $\Sigma$. Let $p$ be a non-negative continuous function integrable with respect to $\mu$. Let $\omega\mapsto h_{\omega}$ be a map from $\Omega$ to continuous, convex functions on $X$, which is weakly measurable, that is, for every $x\in X$ the map $\omega\mapsto h_{\omega}(x)$ is $\Sigma$-measurable, and such that there a exists non-negative number $c$ such that
\begin{equation}\label{eqn:ccc}
\abs{h_{\omega}(x)}\leq c \norm{x}p(\omega)\text{ for all }x\in X, \omega\in\Omega.
\end{equation}
Set 
\begin{equation*}
h(x)=\int_{\Omega}h_{\omega}(x)\,d\mu(\omega).
\end{equation*}
Suppose that $x^*\in X^*$. Then the following conditions are equivalent:
\begin{enumerate}
\item\label{i:boundee} $x^*\leq h$,
\item\label{i:dis} there exists a map $\omega\mapsto x_{\omega}^*$ from $\Omega$ to $X^*$ which is weakly measurable, in the sense that $\omega\mapsto x_{\omega}^*(x)$ is measurable for any $x\in X$, and such that for all $x\in X$ 
\begin{equation*}
x^*(x)=\int_{\Omega}x_{\omega}^*(x)\,d\mu(\omega)
\end{equation*}
and such that for all measurable sets $A\subset \Omega$ and $x\in X$
\begin{equation*}
\int_A x_{\omega}^*(x)\, d\mu(\omega)\leq \int_{A}h_{\omega}(x)\, d\mu(\omega).
\end{equation*}
\end{enumerate}
Moreover, if the above conditions are satisfied and $Z\subset X$ is a countable set, then for $\mu$-almost every $\omega\in \Omega$ there is 
\begin{equation*}
x_{\omega}^*(z)\leq h_{\omega}(z)\text{ for all }z\in Z.
\end{equation*}
\end{theorem}
\begin{proof}
Clearly, \ref{i:dis} implies \ref{i:boundee}. Suppose now that \ref{i:boundee} holds true.
Consider the Banach space $L^1(\Omega,X,\mu)$ of equivalence classes of Borel measurable, Bochner-integrable maps $f\colon \Omega\to X$ with finite norm
\begin{equation*}
\norm{f}_1=\int_{\Omega}\norm{f(\omega)}p(\omega)\,d\mu(\omega).
\end{equation*}
Define a functional $\lambda$ on the subspace of $L^1(\Omega,X,\mu)$ consisting of constant functions by the formula
\begin{equation*}
\lambda(x)=x^*(x).
\end{equation*}
By the assumption, $\lambda\leq h$, on that subspace.
By the Hahn--Banach theorem for convex
majorants we may extend it to a functional $\Lambda$ defined on $L^1(\Omega, X,\mu)$ such that it satisfies
\begin{equation}\label{eqn:compare}
\Lambda(g)\leq \int_{\Omega}h_{\omega}(g(\omega))\,d\mu(\omega)\text{ for all }g\in L^1(\Omega,X,\mu).
\end{equation}
Note that the right-hand side of the above inequality is well-defined by integrability of $g$ and the assumption (\ref{eqn:ccc}). 
Then there exists a weakly measurable map $x_{\omega}^*$ such that for all $g\in L^1(\Omega,X,\mu)$
\begin{equation*}
\Lambda(g)=\int_{\Omega}x_{\omega}^*(g(\omega))\,d\mu(\omega),
\end{equation*} 
cf. \cite[Theorem T51, p. 245]{Meyer1966}.
For $z\in X$ and a measurable set $A\subset\Omega$ consider map $z\mathbf{1}_A\in L^1(\Omega,X,\mu)$. Then it follows by (\ref{eqn:compare}) that
\begin{equation}\label{eqn:A}
\int_A x_{\omega}^*(z)\,d\mu(\omega)\leq \int_A h_{\omega}(z)\,d\mu(\omega).
\end{equation}
Now, for the proof of the last part of the theorem, observe that for any $z\in Z$ it follows by (\ref{eqn:A}) that for $\mu$-almost every $\omega\in\Omega$ there is
\begin{equation*}
x_{\omega}^*(z)\leq h_{\omega}(z).
\end{equation*}
Employing the countable additivity of $\mu$ and countability of $Z$ yields the assertion.
\end{proof}

\begin{remark}\label{rem:el}
Let us elaborate on the relation of $L^1(\Omega, X,\mu)^*$ with $L^{\infty}(\Omega,X^*,\mu)$. In \cite[Chapter IV]{Diestel1977} it is proven that these spaces are isomorphic if and only if $X^*$ enjoys Radon--Nikodym property with respect to the given measure $\mu$. However, the proof of this equivalence, \cite[Chapter IV, Theorem 1]{Diestel1977}, fails if we merely require weak measurability of functions in $L^{\infty}(\Omega,X^*,\mu)$. Hence, the assertion in the above proof of the existence of a weakly measurable map $x_{\omega}^*$ that represents functional $\Lambda$, cf.  \cite[Theorem T51, p. 244]{Meyer1966} and \cite[10.5.4., p. 374]{Bogachev20072}, does not contradict results of \cite{Diestel1977}. For a discussion of measurability of vector-valued maps we refer to \cite[Chapter II]{Diestel1977}.
\end{remark}

\section{Existence of transports}\label{s:existence}

Usually Strassen's theorem, cf. Theorem \ref{thm:varstrassen}, is used to infer the existence of a Markov kernel, which then is used to construct an appropriate transport. In this section, we follow a different strategy of employing Theorem \ref{thm:varstrassen} to directly infer the existence of a transport.

\begin{definition}
Let $\Omega$ be a topological space. Let $\pi$ be a Radon measure on $\Omega\times\Omega$. Let $p_i\colon \Omega\times\Omega\to\Omega$ denote the projections on the first and on the second co-ordinate respectively. Then the $i^{\text{th}}$ marginal of $\pi$ is the push-forward measure ${p_i}_{\#}\pi$ of $\pi$ via $p_i$, $i=1,2$.
\end{definition}

\begin{theorem}\label{thm:mart}
Let $\Omega$ be a set.
Let $\mathcal{G}$ be a convex cone of functions on $\Omega$ that contains constants.
Let $\mathcal{F}$ be the lattice cone generated by $\mathcal{G}$.
Let $p$ be a real-valued, non-negative function in $\mathcal{H}$, proper with respect to $\tau(\mathcal{G})$.

Suppose that $\mathcal{G}$ separates points of $\Omega$ and consists of functions of $(p+1)$-growth.

Suppose that $\mu,\nu$ are Radon, probability measures on $\sigma(\tau(\mathcal{G}))$ such that 
\begin{equation}\label{eqn:finitep}
\int_{\Omega}p\, d\mu<\infty\text{ and }\int_{\Omega}p\, d\nu<\infty
\end{equation}
and that
\begin{equation}\label{eqn:maj}
\int_{\Omega}f\,d\mu\leq\int_{\Omega} f\,d\nu
\end{equation}
for all $f\in\mathcal{F}$.

Then there exists a Radon measure $\pi$ on $\Omega\times\Omega$ with marginals $\mu$ and $\nu$ and such that for all $f\in\mathcal{F}$,
and non-negative, bounded, $\sigma(\tau(\mathcal{A}))$-measurable functions $h$ on $\Omega$
\begin{equation}\label{eqn:thef}
\int_{\Omega}h(\omega_1)f(\omega_1)\, d\mu(\omega_1)\leq\int_{\Omega\times\Omega}h(\omega_1)f(\omega_2)\, d\pi(\omega_1,\omega_2).
\end{equation}
In particular, for all $f\in\mathcal{F}\cap(-\mathcal{F})$,
and $h$ as above 
\begin{equation}\label{eqn:tha}
\int_{\Omega}h(\omega_1)f(\omega_1)\, d\mu(\omega_1)=\int_{\Omega\times\Omega}h(\omega_1)f(\omega_2)\, d\pi(\omega_1,\omega_2).
\end{equation}
Conversely, for a Radon measure $\pi$  on $\Omega\times\Omega$ for which (\ref{eqn:thef}) is satisfied, its marginals $\mu$, $\nu$ satisfy (\ref{eqn:maj}).
\end{theorem}

\begin{proposition}\label{pro:submarti}
Let $\pi$ be as in the above theorem. Let $(X,Y)$ be a random variable with values in $\Omega\times\Omega$ distributed according to $\pi$.
Then for any $f\in\mathcal{F}$ 
the pair $(f(X),f(Y))$ is a one-step submartingale with respect to filtration $(\sigma(X),\sigma(Y))$, i.e.,
\begin{equation*}
f(X)\leq \mathbb{E}(f(Y)\vert \sigma(X))\text{ almost surely.}
\end{equation*}
In particular, for all $f\in\mathcal{F}\cap(-\mathcal{F})$
$(f(X),f(Y))$ is a one-step martingale with respect to filtration $(\sigma(X),\sigma(Y))$.
\end{proposition}
\begin{proof}
Note first that Corollary \ref{col:integral} shows that any $f\in\mathcal{F}$ is integrable with respect to $\mu$ and $\nu$. Let $(X,Y)$ be distributed according to $\pi$. Then, (\ref{eqn:thef}) may rewritten as
\begin{equation*}
\mathbb{E}f(X)\mathbf{1}_A(X)\leq \mathbb{E}f(Y)\mathbf{1}_A(X)
\end{equation*}
for all $A\in\sigma(\tau(\mathcal{A}))$.
This is to say, for all $f\in\mathcal{F}$
\begin{equation*}
f(X)\leq \mathbb{E}(f(Y)\vert \sigma(X)).
\end{equation*}
That is, $(f(X),f(Y))$ is a one-step submartingale.
\end{proof}

\subsection{Lemmata}

\begin{lemma}\label{lem:poussin}
Let $\mathcal{H}$ be the complete lattice cone of functions on $\Omega$ generated by a convex cone $\mathcal{G}$, $1\in \mathcal{G}$. Let $p\in\mathcal{H}$. Let $\xi\colon \mathbb{R}\to \mathbb{R}$ convex and non-decreasing. Then the composition $\xi(p)$ belongs to $\mathcal{H}$.
\end{lemma}
\begin{proof}
By Lemma \ref{lem:supgen} there exists a family $(g_i)_{i\in I}\subset \mathcal{G}$ of functions such that for every $\omega\in\Omega$ there is
\begin{equation*}
p(\omega)=\sup\{g_i(\omega)\mid i\in I\}.
\end{equation*}
Now, as $\xi$ is convex and non-decreasing, there exists a family of non-decreasing affine functions $(\lambda_j)_{j\in J}$ on $\mathbb{R}$ such that for all $t\geq 0$ there is
\begin{equation*}
\xi(t)=\sup\{\lambda_j(t)\mid j\in J\}.
\end{equation*}
Then 
\begin{equation*}
\xi(p)(\omega)=\sup\{\lambda_j(g_i)(\omega)\mid i\in I,j\in J\}.
\end{equation*}
Indeed, for any $j\in J$, there exist $\alpha_j\geq 0$ and $\beta_j\in\mathbb{R}$ such  that $\lambda_j(t)=\alpha_j t+\beta_j$ for all $t\in\mathbb{R}$, so 
\begin{equation*}
\lambda_j(\sup\{g_i(\omega)\mid i\in I\})=\sup\{\lambda_j(g_i)(\omega)\mid i\in I\}.
\end{equation*}
Now, for any $i\in I,j\in J$, function $\lambda_j(g_i)\in\mathcal{G}$. The proof is complete.
\end{proof}

\begin{lemma}\label{lem:sequence}
Let $p$ be a non-negative function on a set $\Omega$. Then 
there exists a monotonically increasing sequence of $[0,1]$-valued functions $(\phi_n)_{n=1}^{\infty}$ on $\Omega$ such that:
\begin{enumerate}
\item $\phi_n=0$ on $p^{-1}([n+1,\infty))$,
\item $\phi_n=1$ on $p^{-1}([0,n])$.
\end{enumerate}
\end{lemma}
\begin{proof}
Let for $n=1,2,\dotsc$
\begin{equation*}
\phi_n=\min\{1,\max\{n+1-p,0\}\}.
\end{equation*}
The required properties of the sequence are readily verified.
\end{proof}

\begin{lemma}\label{lem:depe}
Let $\xi\colon \mathbb{R}\to [1,\infty)$ be a function of superlinear growth. Let $p\colon\Omega\to[1,\infty)$ be a function on $\Omega$. For $f\in \mathcal{D}_p(\Omega)$
\begin{equation}\label{eqn:conver}
\lim_{n\to\infty} f(1-\phi_n) =0\text{ in }\mathcal{D}_{\xi(p)}(\Omega).
\end{equation}
\end{lemma}
\begin{proof}
Let $\epsilon>0$ and pick $n_0$ large enough, so that for all $t\geq n_0$ there is
\begin{equation*}
\frac{t}{\xi(t)}\leq \epsilon.
\end{equation*}
Then for $\omega\in p^{-1}([n_0,\infty))$ we have
\begin{equation*}
\frac{p(\omega)}{\xi(p(\omega))}\leq \epsilon.
\end{equation*}
Therefore for $n\geq n_0$
\begin{equation}\label{eqn:factor}
\frac{\abs{f(\omega)(1-\phi_n)(\omega)}}{\xi(p(\omega))}=\frac{\abs{f(\omega)}}{p(\omega)}\frac{p(\omega)(1-\phi_n)(\omega)}{\xi(p(\omega))}\leq \epsilon \norm{f}_{\mathcal{D}_p(\Omega)}.
\end{equation}
The proof is complete.
\end{proof}

\begin{lemma}\label{lem:envelope}
Let $\Omega$ be a set. Let $\mathcal{G}$ be a convex cone of continuous functions on $\Omega$ and let $\mathcal{H}$ be the complete lattice cone generated by $\mathcal{G}$. Let $\xi(p)\in\mathcal{H}$, $\xi(p)\geq 1$. For $\omega\in\Omega$ and a function $x\in\mathcal{D}_{\xi(p)}(\Omega)$ define
\begin{equation*}
h_{\omega}(x)=\inf\Big\{-y(\omega)\mid  y\in\mathcal{H}, -y\geq x\Big\}.
\end{equation*}
Then:
\begin{enumerate}
\item\label{i:ef} for each $x\in \mathcal{D}_{\xi(p)}(\Omega)$, $\omega\mapsto h_{\omega}(x)$ belongs to $-\mathcal{H}$,
\item\label{i:convex} for each $\omega$, $x\mapsto h_{\omega}(x)$ is a continuous, subadditive and positively homogeneous function on $\mathcal{D}_{\xi(p)}(\Omega)$,
\item\label{i:bound} $\abs{h_{\omega}(x)}\leq \norm{x}_{\mathcal{D}_{\xi(p)}(\Omega)}\xi(p(\omega))$, for $\omega\in\Omega$ and $x\in \mathcal{D}_{\xi(p)}(\Omega)$,
\item\label{i:a} for $f\in(-\mathcal{H})\cap \mathcal{D}_{\xi(p)}(\Omega)$, $h_{\omega}(f)=f(\omega)$; thus for all $f\in\mathcal{H}\cap(-\mathcal{H})$, $h_{\omega}(f)=f(\omega)$,
\item\label{i:positive} if $x\leq 0$, then also $h_{\omega}(x)\leq 0$,
\item\label{i:lowerb} for each $x\in  \mathcal{D}_{\xi(p)}(\Omega)$,  $x\leq h_{\omega}(x)$.
\end{enumerate}
\end{lemma}
\begin{proof}
Note first that $h_{\omega}(x)$ is well-defined, as for any $x\in\mathcal{D}_{\xi(p)}(\Omega)$ the set of all  $y$ such that $y\in -\mathcal{H}$ such that $-y\geq x$ is non-empty, as it contains $-\xi(p)\norm{x}_{\mathcal{D}_{\xi(p)}(\Omega)}$. The bound \ref{i:bound} follows immediately. 
In particular for each $\omega\in\Omega$, $x\mapsto h_{\omega}(x)$ is a continuous function.
As $\mathcal{H}$ is a convex cone, so $x\mapsto h_{\omega}(x)$ is subadditive and positively homogeonous.

For any $x\in\mathcal{D}_{\xi(p))}(\Omega)$ the function $\omega\mapsto h_{\omega}(x)$ belongs to $-\mathcal{H}$, as $\mathcal{H}$ is stable under suprema. 

The last three assertions are trivial.
\end{proof}

\begin{lemma}\label{lem:net}
Let $\mathcal{G}$ be a convex cone of functions on $\Omega$. Let $\mathcal{H}$ be the complete lattice cone generated by $\mathcal{G}$. Then for any $f\in\mathcal{H}$ there exists an increasing net $(f_{\alpha})_{\alpha\in A}$ of functions in the lattice cone generated by $\mathcal{G}$ that converges to $f$.
\end{lemma}
\begin{proof}
By Lemma \ref{lem:supgen}, for $f\in\mathcal{H}$ there exists a family $(g_i)_{i\in I}$ of functions in $\mathcal{G}$ such that
\begin{equation*}
f(\omega)=\sup\{g_i(\omega)\mid i\in I\}\text{ for all }\omega\in\Omega.
\end{equation*}
Let $A$ be the set of finite subsets of $I$. For $\alpha\in A$, we define 
\begin{equation*}
f_{\alpha}(\omega)=\max\{g_i(\omega)\mid i\in \alpha\}\text{ for }\omega\in\Omega.
\end{equation*}
Clearly, $A$ is a net, with order given by inclusion, and thus $(f_{\alpha})_{\alpha\in A}$ is a net.
We claim that for all $\omega\in\Omega$
\begin{equation*}
\lim_{\alpha\in A}f_{\alpha}(\omega)=f(\omega).
\end{equation*}
Let $\omega\in \Omega$ and $\epsilon>0$. We need to show that there is $\alpha_0\in A$ such that for $\alpha\supset \alpha_0$ we have
\begin{equation*}
f_{\alpha}(\omega)>f(\omega)-\epsilon.
\end{equation*}
This is however trivial, it suffices to take $\alpha_0=\{i\}$, where $i\in I$ is such that
\begin{equation*}
f_i(\omega)>f(\omega)-\epsilon.
\end{equation*}
\end{proof}

Let us bring to our attention the following lemma, a modification of \cite[Lemma 7.2.6., p. 75]{Bogachev20072}.
\begin{lemma}\label{lem:taulem}
Let $\mu$ be a finite Radon measure on a topological space $\Omega$ and let $(f_{\alpha})_{\alpha\in A}$ be an increasing net of lower semi-continuous, functions on $\Omega$ converging to $f$. Assume that for some $\alpha_0\in A$, $f_{\alpha_0}$ is integrable and that $f$ is integrable with respect to $\mu$. Then
\begin{equation*}
\lim_{\alpha\in A}\int_{\Omega}f_{\alpha}\, d\mu=\int_{\Omega}f\,d\mu.
\end{equation*}
\end{lemma}
\begin{proof}
For natural $m,n$ and $\alpha\in A$ define
\begin{equation*}
f_{\alpha,n,m}=-m+\sum_{l=-m}^m\sum_{k=1}^n\frac1n\mathbf{1}_{f_{\alpha}^{-1}\big((l+\frac{k-1}{n},\infty)\big)}
\end{equation*}
and 
\begin{equation*}
f_{n,m}=-m+\sum_{l=-m}^m\sum_{k=1}^n\frac1n\mathbf{1}_{f^{-1}\big((l+\frac{k-1}{n},\infty)\big)}.
\end{equation*}
By lower semi-continuity of $f_{\alpha}$, we see that the sets $f_{\alpha}^{-1}\big((l+\frac{k-1}{n},\infty)\big)$ are open and form a net increasing to the open set $f^{-1}\big((l+\frac{k-1}{n},\infty)\big)$, for each $m,n$. As $\mu$ is a Radon measure, it is $\tau$-additive (see \cite[Definition 7.2.1., p. 73 and Proposition 7.2.2., p. 74]{Bogachev20072}). Therefore for every $n,m$
\begin{equation*}
\lim_{\alpha\in A}\int_{\Omega}f_{\alpha,n,m}\, d\mu=\int_{\Omega}f_{n,m}\, d\mu.
\end{equation*}
Observe that  
\begin{equation*}
\text{if }-m<f_{\alpha}(\omega)\leq m\text{, then }\abs{f_{\alpha,n,m}-f_{\alpha}}(\omega)\leq\frac1n.
\end{equation*} 
Similarly, 
\begin{equation*}
\text{if }-m<f(\omega)\leq m\text{, then }\abs{f_{n,m}-f}(\omega)\leq\frac1n.
\end{equation*}
Now, let $\alpha_0\in A$ be such that $f_{\alpha_0}$ is integrable with respect to $\mu$. We may suppose that 
\begin{equation*}
f_{\alpha_0}\leq f_{\alpha}\leq f\text{ for all }\alpha\in A.
\end{equation*}
Take $\epsilon>0$. Let $m_0$ be such that 
\begin{equation*}
\int_{\Omega}\abs{f} \mathbf{1}_{f^{-1}((-m_0,m_0])}\, d\mu<\epsilon\text{ and }\int_{\Omega}\abs{f_{\alpha_0}} \mathbf{1}_{f^{-1}((-m_0,m_0])}\, d\mu<\epsilon.
\end{equation*}
Let $n_0$ be some number such that $n_0>\frac{1}{\epsilon}$. Let $\alpha_1\in A$ be such that for all $\alpha\geq \alpha_1$ we have  
\begin{equation*}
\int_{\Omega}f_{\alpha,n_0,m_0}\, d\mu> \int_{\Omega}f_{n,m}\, d\mu-\epsilon.
\end{equation*}
Then 
\begin{equation*}
\int_{\Omega}f\, d\mu\geq \int_{\Omega}f_{\alpha}\, d\mu> \int_{\Omega}f\, d\mu-5\epsilon.
\end{equation*}
\end{proof}

\begin{remark}\label{rem:finite}
Suppose that $\mathcal{F}$ is the lattice cone of functions $\Omega$ generated by a convex cone $\mathcal{G}$. Let $\mathcal{H}$ be the complete lattice cone generated by $\mathcal{G}$. Lemma \ref{lem:supgen} shows that any element of $\mathcal{H}$ is a supremum of elements of $\mathcal{G}$. Then the above lemma shows that for two Radon probability measures $\mu,\nu\in\mathcal{P}(\Omega)$ 
\begin{equation*}
\mu\prec_{\mathcal{F}}\nu\text{ is equivalent to }\mu\prec_{\mathcal{H}}\nu.
\end{equation*}
Therefore, there is no loss of generality in assuming that $\mathcal{F}$ is a complete lattice cone.
\end{remark}

\subsection{Proof}

\begin{proof}[Proof of Theorem \ref{thm:mart}]
Remark \ref{rem:finite} shows that if $\mu\prec_{\mathcal{F}}\nu$, then $\mu\prec_{\mathcal{H}}\nu$ as well.
Let us first note that all functions in $\mathcal{H}$ are $\sigma(\tau(\mathcal{G}))$-measurable, as they are lower semi-continuous with respect to $\tau(\mathcal{G})$, by Lemma \ref{lem:supgen}.


By Theorem \ref{thm:poussin} there exists a continuous, convex, increasing function $\xi\colon \mathbb{R}\to [1,\infty)$ of superlinear growth such that 
\begin{equation}\label{eqn:ppp}
\int_{\Omega}\xi(p)\, d\mu<\infty, \int_{\Omega}\xi(p)\, d\nu<\infty.
\end{equation}

Consider the Banach space $X=\mathcal{C}_{\xi(p)}(\Omega)$.
Let $x^*$ be an element of $X^*$ represented by a measure $\nu\in\mathcal{P}(\Omega)$. 

Let $h$ be defined as in Lemma \ref{lem:envelope}; note that $\xi(p)\in\mathcal{H}$, by Lemma \ref{lem:poussin}.
Pick some $x\in X$. We shall show that 
\begin{equation}\label{eqn:major}
x^*(x)\leq \int_{\Omega}h_{\omega}(x)\,d\mu(\omega).
\end{equation} 

By Lemma \ref{lem:envelope}, \ref{i:ef}, $\omega\mapsto h_{\omega}(x)$ belongs to $-\mathcal{H}$, so it is $\sigma(\tau(\mathcal{G}))$-measurable.
Moreover, by Lemma \ref{lem:envelope}, \ref{i:bound}, and (\ref{eqn:ppp}) it is integrable with respect to $\mu$ and to $\nu$.
By the definition of $h$, \ref{i:lowerb}, and the assumption (\ref{eqn:maj}) we have
\begin{equation*}
x^*(x)=\int_{\Omega} x\,d\nu \leq \int_{\Omega} h_{\omega}(x)\,d\nu(\omega)\leq\int_{\Omega}h_{\omega}(x)\,d\mu(\omega).
\end{equation*}
Theorem \ref{thm:varstrassen} now implies that there exists a weakly measurable function $\omega\mapsto x_{\omega}^*$ with values in $X^*$ such that 
\begin{equation}\label{eqn:repr}
x^*(x)=\int_{\Omega}x_{\omega}^*(x)\,d\mu(x)\text{ for all }x\in X
\end{equation}
and
\begin{equation}\label{eqn:esti}
\int_{A}x_{\omega}^*(x)\, d\mu(\omega)\leq\int_A h_{\omega}(x)\, d\mu(\omega)\text{ for all }x\in X\text{ and }A\in\sigma(\tau(\mathcal{G})).
\end{equation}
For any set $A\in\sigma(\tau(\mathcal{G}))$ we consider a functional $x^*_A$ on $X$ given by the formula
\begin{equation*}
x_A^*(x)=\int_Ax_{\omega}^*(x)\, d\mu(\omega).
\end{equation*}
By Lemma \ref{lem:envelope}, \ref{i:a} and \ref{i:positive}, and by (\ref{eqn:esti}), for all $A\in\sigma(\tau(\mathcal{G}))$,
\begin{equation}\label{eqn:necessities}
 x_A^*(1)=\mu(A)\text{ and }x_A^*(x)\geq 0\text{ if }x\geq 0.
\end{equation}
We claim that it also satisfies the convergence condition of Theorem \ref{thm:daniell} on the space $\mathcal{C}_p(\Omega)$. Let us fix a set $A\in\sigma(\tau(\mathcal{G}))$. Let $(f_n)_{n=1}^{\infty}\subset\mathcal{C}_p(\Omega)$ be a sequence on non-negative functions converging monotonically to zero. We shall show that
\begin{equation}\label{eqn:claim}
\lim_{n\to\infty}x_A^*(f_n)=0.
\end{equation}
By (\ref{eqn:necessities}) functional $x_A^*$ is non-negative. Therefore the above limit is non-negative.
Let $\epsilon>0$. Let $(\phi_m)_{m=1}^{\infty}$ be a sequence constructed in Lemma \ref{lem:sequence}. By Lemma \ref{lem:depe}, there exists $m_0$ such that for all $m\geq m_0$
\begin{equation}\label{eqn:epsilonik}
\norm{f_1(1-\phi_m)}_{D_{\xi(p)}(\Omega)}<\epsilon.
\end{equation}
Now, (\ref{eqn:esti}) and Lemma \ref{lem:envelope}, \ref{i:convex}, \ref{i:positive} imply that 
\begin{equation*}
x_A^*(f_n)\leq \int_A \Big(h_{\omega}(f_n\phi_m)+h_{\omega}(f_1(1-\phi_m))\Big)\,d\mu(\omega).
\end{equation*}
Thanks to Lemma \ref{lem:envelope}, \ref{i:bound}, and to (\ref{eqn:epsilonik}), we get that
\begin{equation}\label{eqn:xbiga}
x_A^*(f_n)\leq (\norm{f_n\phi_{m_0}}_{\mathcal{D}_{\xi(p)}(\Omega)}+\epsilon)\int_A \xi(p)\,d\mu.
\end{equation}
Now, by Lemma \ref{lem:sequence}, $\phi_{m_0}=0$ on $p^{-1}([0,m_0+1])$. By the assumption, the set $p^{-1}([0,m_0+1])$ is compact. By the Dini's lemma, the sequence $(f_n)_{n=1}^{\infty}$ converges uniformly on $p^{-1}([0,m_0+1])$ to zero. As $\xi(p)$ is bounded from below, we see that there exists $n_0$ such that for $n\geq n_0$ we have
\begin{equation*}
\norm{f_n\phi_{m_0}}_{\mathcal{D}_{\xi(p)}(\Omega)}< \epsilon.
\end{equation*}
This, together with (\ref{eqn:xbiga}) and $\mu$-integrability of $\xi(p)$ completes the proof of (\ref{eqn:claim}).

Now, (\ref{eqn:necessities}), (\ref{eqn:claim}) imply, thanks to Theorem \ref{thm:daniell}, that for any $A\in\sigma(\tau(\mathcal{G}))$ there exists a unique non-negative measure $P(A,\cdot)$ on $\sigma(\mathcal{C}_p(\Omega))$ such that for all $x\in\mathcal{C}_p(\Omega)$
\begin{equation}\label{eqn:pe}
\int_{\Omega}x\, dP(A,\cdot)= \int_A x_{\omega}^*(x)\,d\mu(\omega)
\end{equation}
and such that 
\begin{equation}\label{eqn:mu}
P(A,\Omega)=\mu(A).
\end{equation}
Moreover, any $f\in \mathcal{C}_p(\Omega)$ is integrable with respect to $P(A,\cdot)$. 

By uniqueness, we see that the map
\begin{equation}\label{eqn:addit}
\sigma(\tau(\mathcal{A}))\ni A\mapsto P(A,\cdot)\in\mathcal{M}(\Omega)
\end{equation}
is additive. 

Now, $P(A,\cdot)$ is tight, cf. Definition \ref{def:baire}, \ref{i:tight}. Indeed, let $\epsilon>0$ and let $K\subset\Omega$ be a compact set such that $\nu(\Omega\setminus K)<\epsilon$. If $R\in\sigma(\mathcal{C}_p(\Omega))$ is such that $K\cap R=\emptyset$ then, by (\ref{eqn:addit}),
\begin{equation*}
P(A,R)\leq P(\Omega,R)=\nu(R)\leq \nu(\Omega\setminus K)<\epsilon.
\end{equation*}
As $\mathcal{G}$ separates points of $\Omega$, $\tau(\mathcal{G})$ is Hausdorff. By \cite[Corollary 7.1.8., p. 70]{Bogachev20072} and \cite[Theorem 7.3.2.,(i), p. 78]{Bogachev20072}, any tight Baire measure on $\Omega$ admits a unique extension to a Radon measure. 
Therefore $P(A,\cdot)$ uniquely extends to a Radon measure on $\sigma(\tau(\mathcal{G}))$.
We shall denote this extension still by $P$.

Let $\Theta$ denote the algebra of subsets of $\Omega\times\Omega$ that are of the form
\begin{equation*}
\bigcup_{i=1}^m A_i\times B_i\text{ with }A_i,B_i\in\sigma(\tau(\mathcal{G})).
\end{equation*}
We define $\pi$ on $\Theta$ by the formula
\begin{equation}\label{eqn:pi}
\pi\bigg(\bigcup_{i=1}^m A_i\times B_i\bigg)=\sum_{i=1}^m P(A_i,B_i),
\end{equation}
whenever $(A_i\times B_i)_{i=1}^m$ is a pairwise disjoint family of sets, $A_i,B_i\in\sigma(\tau(\mathcal{G}))$.
We claim that $\pi$ is regular, additive, of bounded variation and tight. 

Clearly, it is additive by (\ref{eqn:addit}). It is non-negative, so its total variation is at most one. 

Note also that by (\ref{eqn:repr}) and (\ref{eqn:pe}) and uniqueness,
\begin{equation}\label{eqn:nu}
P(\Omega,\cdot)=\nu.
\end{equation}

Now, $\mu,\nu$ are Radon measures on $\Omega$. Let $\epsilon>0$. Pick compact sets $K_1,K_2\subset\Omega$ for which
\begin{equation*}
\mu(\Omega\setminus K_1)<\epsilon\text{ and }\nu(\Omega\setminus K_2)<\epsilon.
\end{equation*}
Then $K_1\times K_2$ is compact in $\Omega\times\Omega$ and
\begin{equation*}
\pi(\Omega\setminus K_1\times K_2)\leq P(\Omega,\Omega\setminus K_1)+P(\Omega\setminus K_2,\Omega)<2\epsilon.
\end{equation*}
Thus, $\pi$ is tight.

Similarly, we can show that $\pi$ is regular. Indeed, let $A,B\in\sigma(\tau(\mathcal{G}))$, $\epsilon>0$. Take a closed set $E\subset\Omega$ such that $\mu(A\setminus E)<\epsilon$. Let $F\subset B$ be a closed set such that $\nu(B\setminus F)<\epsilon$. Then $E\times F\subset A\times B$ is a closed set such that
\begin{equation*}
\pi(A\times B\setminus E\times F)\leq P(A\setminus E, \Omega)+P(\Omega, B\setminus F)=\mu(A\setminus E)+\nu(B\setminus F)<2\epsilon.
\end{equation*}

Now, \cite[Theorem 7.3.2.(i), p. 78]{Bogachev20072} allows us to infer that $\pi$ extends uniquely to a Radon measure on $\Omega\times\Omega$.

By the definition of $\pi$ and by (\ref{eqn:mu}) and (\ref{eqn:nu}) we have that for all $A\in\sigma(\tau(\mathcal{G}))$
\begin{equation*}
\pi(A\times \Omega)=\mu(A)\text{ and }\pi(\Omega\times A)=\nu(A).
\end{equation*}
That is $\pi$ has the required marginals.

We see  that (\ref{eqn:esti}), (\ref{eqn:pe}) and Lemma \ref{lem:envelope}, \ref{i:a}, imply that for all $f\in\mathcal{H}\cap\mathcal{C}_p(\Omega)$ we have
\begin{equation}\label{eqn:continu}
\int_A f\, d\mu\leq\int_{\Omega} f\, d\pi(A,\cdot).
\end{equation}
Now,  $\mathcal{G}\subset\mathcal{C}_p(\Omega)$, (\ref{eqn:continu}), Lemma \ref{lem:net} and Lemma \ref{lem:taulem} imply (\ref{eqn:thef}).

Then (\ref{eqn:tha}) follows immediately, as $\mathcal{F}\cap(-\mathcal{F})$ is a linear subspace.

The converse claim is trivially satisfied, by taking $h=1$.
\end{proof}

\subsection{Remarks}

\begin{remark}
Suppose that $p\in\mathcal{H}$ is some real-valued function such that $p_+$ is proper, $\mathcal{G}$ consists of functions of $(p_++1)$-growth and (\ref{eqn:finitep}) is satisfied -- with  negative infinity values allowed. Then taking $p_+$ in place of $p$, we can conclude the assertions of the theorem.
\end{remark}

\begin{remark}
An important feature of Theorem \ref{thm:mart} is that we do not assume that $\Omega$ is locally compact in a topology with respect to which $\mathcal{G}$ consists of continuous functions. That will allow us to employ the theorem to study, e.g., martingale transports on infinite-dimensional Banach spaces.
\end{remark}

\begin{remark}
Suppose that $\Omega$ is a locally compact Hausdorff space, and that $\mathcal{G}$ consists of continuous functions. In this case the proof of the above theorem is significantly simpler. This is due to separability of the space of compactly supported functions and thanks to the Riesz' representation theorem.
\end{remark}

\begin{remark}
Let us note an application of Strassen's theorem to optimal transport problems presented in  \cite[1.8.3., p. 58]{Lusky1981}.
 \end{remark}

 \begin{remark}
     A generalisation of the theorem that would be interesting to investigate concerns Radon probability measures $\mu,\nu$ for which (\ref{eqn:maj}) holds true, with all $f\in\mathcal{F}$ integrable with respect to $\mu$ and $\nu$, but with no assumption (\ref{eqn:finitep}). Tightly related to this context is also Remark \ref{rem:moments}.
 \end{remark}

\section{Gelfand transform}\label{s:gelfand}

\begin{theorem}\label{thm:banach}
Suppose that $X$ is a normed space.  Let $X^*$ be its dual Banach space. Let $G\subset X$ be a linearly dense, convex cone. 
Let $H$ denote the lattice cone of functions on $X^*$ generated by $G$ and constant functions.
Suppose that $\mu,\nu$ are two Radon probability measures on $\sigma(\tau(X))$ such that
\begin{equation*}
\int_{X^*}\norm{\cdot}\, d\mu<\infty, \int_{X^*}\norm{\cdot}\, d\nu<\infty, 
\end{equation*}
and that for any function $f\in H$
there is
\begin{equation*}
\int_{X^*} f\,d\mu\leq\int_{X^*} f\, d\nu.
\end{equation*}
Then there exists a Radon measure $\pi$ on $X^*\times X^*$ with marginals $\mu$ and $\nu$ such that for any $\sigma(\tau(X))$-measurable bounded, non-negative function $h$ and any $f\in H$ 
\begin{equation*}
\int_{Y}h(y_1)f(y_1)\, d\mu(y_1)\leq \int_{Y\times Y}h(y_1)f(y_2)\, d\pi(y_1,y_2).
\end{equation*}
\end{theorem}
\begin{proof}
We shall apply Theorem \ref{thm:mart} with $\Omega=X^*$, $\mathcal{G}=G+\mathbb{R}$ and $p=\norm{\cdot}$.
Let $\mathcal{H}$ be the complete lattice cone generated by $\mathcal{G}$.
The fact that $p$ belongs to $\mathcal{H}$ follows by the Hahn--Banach theorem and the linear density of $G$ in $X$.

Then the fact that sub-level sets of $p$ are compact in the topology generated by $\mathcal{G}$ follows by the Banach--Alaoglu theorem. 

The fact that $\mathcal{G}$ separates points of $X^*$ follows by the Hahn--Banach theorem and the  density of  $G-G$ in $X$.

The fact that all all functions in $\mathcal{G}$ are of $p$-growth follows by their continuity.

Thus, the assumptions of Theorem \ref{thm:mart} are satisfied. Hence, its conclusion holds true. This is however the assertion of the currently proven theorem.
\end{proof}

\begin{remark}\label{rem:moments}
The assumption that the moments $\int_{X^*}\norm{\cdot}\, d\mu,\int_{X^*}\norm{\cdot}\, d\nu$ are finite is rather strong, see e.g., \cite[Theorem 7.14.45., p. 145]{Bogachev20072}. 
Does the above theorem hold true if we merely assume that each element of $H$ is integrable with respect to $\mu$ and $\nu$? Simple examples show that these two properties are not equivalent for infinite-dimensional spaces.
\end{remark}

\begin{remark}
Let us note that martingales in Banach spaces is a widely studied topic with several applications. We refer the reader to \cite{Pisier2016}, \cite{Cox2017} and to \cite[Section 7.14(xii).]{Bogachev20072} for related studies.
\end{remark}

Let $X$ be a normed space. Below we shall consider bounded weak* topology on a dual space $X^*$. A set $U\subset X^*$ is open in bounded weak* topology provided that the intersection of $U$ with any closed ball in $X^*$ is relatively weakly* open. Let us recall the following theorem of Krein and \v{S}mulian, \cite{Krein1940}. 

\begin{theorem}\label{thm:krein-smulian}
Let $X$ be a Banach space. Then a convex set $Z\subset X^*$ is closed in the weak* topology if and only if it is closed in the bounded weak* topology.
\end{theorem}

\begin{theorem}\label{thm:embed}
Let $T_1$ be the class of triples $(\Omega,\mathcal{G},p)$ such that
\begin{enumerate}
\item $\mathcal{G}$ is a convex cone of functions on a set $\Omega$ that contains constants and separates points of $\Omega$,
\item $p$ is a real-valued, non-negative function in the complete lattice cone $\mathcal{H}$ generated by $\mathcal{G}$, proper with respect to $\tau(\mathcal{G})$,
\item $\mathcal{G}\subset\mathcal{C}_{p+1}(\Omega)$.
\end{enumerate}
Let $T_2$ be the class of triples $(Z,G,\norm{\cdot}_{X^*})$ such that 
\begin{enumerate}
\item $G$ is a linearly dense convex cone in a normed space $X$,
\item  $Z$ is a subset of $X^*$, boundedly weakly* closed, linearly dense in weak* topology.
\end{enumerate}

To any $(\Omega,\mathcal{G},p)$ in the class $T_1$ we may assign $(Z,G,\norm{\cdot}_{H^*})$ in the class $T_2$ in the following way. Let
\begin{enumerate}
\item\label{i:hzero} $G=\mathcal{G}$, normed by $\norm{\cdot}_{\mathcal{D}_{p+1}(\Omega)}$, $X=G-G$,
\item\label{i:phi} $Z=\Phi(\Omega)$, where $\Phi\colon \Omega\to X^*$ is given by the formula
\begin{equation*}
\Phi(\omega)(g_1-g_2)=g_1(\omega)-g_2(\omega)\text{ for }g_1,g_2\in G,\omega\in\Omega.
\end{equation*}
\end{enumerate}
Then
\begin{enumerate}
\item\label{i:homeo} $\Phi$ is a homeomorphism of $(\Omega,\tau(\mathcal{G}))$ and $(\Phi(\Omega),\tau(G))$,
\item\label{i:norm} $\norm{\Phi(\omega)}_{X^*}=p(\omega)+1$, for $\omega\in\Omega$,
\item\label{i:closed} $\Phi(\Omega)$ is a closed subset of $X^*$ in the bounded weak* topology, linearly dense in the weak* topology,
\item\label{i:aextend} for any $g\in \mathcal{G}$, $g(\Phi^{-1})$ extends to a weakly* continuous linear map on $X^*$,  
\item\label{i:arestrict} for any weakly* continuous linear map $h$ on $X^*$, the function 
\begin{equation*}
\omega\mapsto h(\Phi(\omega))
\end{equation*}
belongs to $\mathcal{G}-\mathcal{G}$,
\item\label{i:fextend} for any $f$ in the complete lattice cone $\mathcal{H}$ generated  by $\mathcal{G}$, $f(\Phi^{-1})$ extends to a map in the complete lattice cone $H$ generated by $G$,
\item\label{i:frestrict} for any map $g\in H$, the function 
\begin{equation*}
\omega\mapsto g(\Phi(\omega))
\end{equation*}
belongs to $\mathcal{H}$.
\end{enumerate}
\end{theorem}
\begin{proof}
Let $(\Omega,\mathcal{G},p)$ be in the class $T_1$. 
The fact that  $\Phi$ is homeomorphism onto $\Phi(\Omega)$, i.e., \ref{i:homeo}, follows directly from the definitions of the considered weak topologies.

Note that, by \ref{i:hzero} and \ref{i:phi},
\begin{equation*}
\norm{\Phi(\omega)}_{X^*}=\sup\{\abs{(g_1-g_2)(\omega)}\mid \norm{g_1-g_2}_{\mathcal{D}_{p+1}(\Omega)}\leq 1\}\leq (p(\omega)+1).
\end{equation*}
However, by Lemma \ref{lem:supgen}, $(p+1)\in\mathcal{H}$ is a supremum of functions in $\mathcal{G}$, we get an equality in the above inequality, i.e., \ref{i:norm} is proven.

Let us check that $Z=\Phi(\Omega)$ is boundedly weakly* closed. It suffices to show that the intersection of $Z$ with any closed ball is weakly* closed. Let $(\omega_{\alpha})_{\alpha\in A}$ be a net in $\Omega$ such that $(\Phi(\omega_{\alpha}))_{\alpha\in A}$ is bounded and converges to some $\xi\in X^*$ in the weak* topology. By \ref{i:norm} we see that $(p(\omega_{\alpha}))_{\alpha\in A}$ is bounded. As $p$ is proper with respect to $\tau(\mathcal{G})$, the net $(\omega_{\alpha})_{\alpha\in A}$ is a subset of a compact set. As $\Phi$ is a homeomorphism, its image is compact. As $G$ separates points of $X^*$, the image is closed in the weak* topology. Thus $\xi\in\Phi(\Omega)$.

We shall prove that $\Phi(\Omega)$ is linearly dense in $X^*$ in the weak* topology. Suppose this is not the case. Then there exist some $v\in X^*$ and a weakly* continuous functional on $X^*$ separating the weak* closure of $\Phi(\Omega)$ and $v$, i.e., there is some $g_1,g_2\in\mathcal{G}$ such that $\Phi(\Omega) $ annihilates $g_1-g_2$ and $v(g_1-g_2)>0$. It follows that $g_1-g_2=0$, which contradicts the above strict inequality. Thus \ref{i:closed} is proven.

It follows also that if $(Z,G,\norm{\cdot}_{X^*})$ is defined as in the formulation of the theorem, then it is in the class $T_2$. 

Let $g\in\mathcal{G}$. Extension in \ref{i:aextend} is given by the formula
\begin{equation*}
X^*\ni x^*\mapsto x^*(g)\in\mathbb{R}.
\end{equation*}
We shall now define an extension for $f\in\mathcal{H}$. By Lemma \ref{lem:supgen} for $f\in\mathcal{H}$ there exists a family $(g_i)_{i\in I}$ of elements of $\mathcal{G}$ such that
\begin{equation}\label{eqn:repef}
f=\sup\{g_i\mid i\in I\}.
\end{equation}
We define a map
\begin{equation*}
X^*\ni x^*\mapsto \sup\{x^*(g_i)\mid i\in I\}\in (-\infty,\infty].
\end{equation*}
By (\ref{eqn:repef}), this is an extension. As functionals $x^*\mapsto x^*(g)$ are weakly* continuous, for $g\in G$, this extension belongs to the complete lattice cone $H$ generated by $G$. Thus \ref{i:fextend} is proven.

Note that any weakly* continuous linear functional $h$ on $X^*$ is given by some element $h_0\in G-G$. Therefore for $\omega\in\Omega$, 
\begin{equation*}
 h(\Phi(\omega))=h_0(\omega).
\end{equation*}
This proves \ref{i:arestrict}. Point \ref{i:frestrict} is proven similarly.
\end{proof}

\begin{remark}\label{rem:complete}
If $\mathcal{G}-\mathcal{G}$ equipped with $\norm{\cdot}_{\mathcal{D}_{p+1}(\Omega)}$ is complete, then using Theorem \ref{thm:krein-smulian}, one may replace the bounded weak* topology in the above theorem by the weak* topology.
\end{remark}

\begin{corollary}
Theorem \ref{thm:mart} and Theorem \ref{thm:banach} are equivalent.
\end{corollary}
\begin{proof}
Employing the results of Theorem \ref{thm:embed}, we see that measures $\mu$, $\nu$ are Radon probability measures on $\sigma(\tau(\mathcal{G}))$ with
\begin{equation*}
\int_{\Omega}p\, d\mu<\infty,\int_{\Omega}p\, d\nu<\infty
\end{equation*}
and
\begin{equation*}
\int_{\Omega}f\, d\mu\leq\int_{\Omega}f\, d\nu\text{ for all }f\in\mathcal{F}
\end{equation*}
if and only if 
$\Phi_{\#}\mu$, $\Phi_{\#}\nu$ are Radon probability measures on $\sigma(\tau(G))$ with
\begin{equation*}
\int_Z\norm{\cdot}_{X^*}\, d\Phi_{\#}\mu<\infty,\int_Z\norm{\cdot}_{X^*}\, d\Phi_{\#}\nu<\infty 
\end{equation*}
and
\begin{equation*}
\int_Zh\, d\Phi_{\#}\mu\leq\int_Zh\, d\Phi_{\#}\nu\text{ for all }h\in H.
\end{equation*} 
Similarly, a Radon measure $\pi$ on $\Omega\times\Omega$ with marginals $\mu$ and $\nu$ is such that 
\begin{equation*}
\int_{\Omega}e(\omega_1)f(\omega_1)\, d\mu(\omega_1)\leq\int_{\Omega\times\Omega}e(\omega_1)f(\omega_2)\,d\pi(\omega_1,\omega_2)
\end{equation*}
for all bounded, measurable, non-negative $e$ and all $f\in\mathcal{F}$
if and only if $(\Phi,\Phi)_{\#}\pi$ is a Radon measure on $X^*\times X^*$ with marginals $\Phi_{\#}\mu$ and $\Phi_{\#}\nu$ such that 
\begin{equation*}
\int_{X^*}s(h_1)h(g_1)\, d\Phi_{\#}\mu(g_1)\leq\int_{X^*\times X^*}s(g_1)h(g_2)\,d(\Phi,\Phi)_{\#}\pi(g_1,g_2)
\end{equation*}
for all bounded measurable, non-negative $s$ and $h\in H$.

It suffices now to observe that if $\rho$ is a Radon probability measure on $X^*\times X^*$ with marginals $\Phi_{\#}\mu,\Phi_{\#}\nu$, then there exists a Radon probability measure $\pi$ on $\Omega\times\Omega$ such that $\rho=(\Phi,\Phi)_{\#}\pi$. Indeed, the condition on marginals implies that $\rho$ is concentrated on $\Phi(\Omega)\times\Phi(\Omega)$ and thus 
\begin{equation*}
\pi=(\Phi^{-1},\Phi^{-1})_{\#}\rho
\end{equation*}
satisfies our requirements.
\end{proof}

\begin{remark}\label{rem:eui}
The theorem says that, up to a homeomorphism, the assumptions in Theorem \ref{thm:mart}, with  $\mathcal{G}$ being a linearly dense convex cone, imply that $\Omega$ is a boundedly weakly* closed, linearly dense, subset of a dual Banach space, that $p+1$ extends to the norm on that Banach space and that the convex cone $\mathcal{G}$ consists of weakly* continuous linear functionals. 
\end{remark}

\section{Topological interlude}

\begin{lemma}\label{lem:regular}
Let $\Omega$ be a set and let $\mathcal{G}$ be a convex cone of functions on $\Omega$ containing constants. Then the topology generated by $\mathcal{G}$ is completely regular.
Moreover, if $F\subset\Omega$ is closed and $\omega\notin F$, there exists a non-negative function $f$ in the vector lattice generated by $\mathcal{G}$, such that 
\begin{equation*}
f(\omega)>0\text{ and }f(F)=\{0\}.
\end{equation*}
\end{lemma}
\begin{proof}
It suffices to show that if $F\subset\Omega$ is a closed set and $\omega\notin F$, then there exists a non-negative continuous function $f$ on $\Omega$ such that $f=0$ on $F$ and $f(\omega)>0$.
As $\omega\notin F$, there exists an open set $U$ containing $\omega$ and disjoint from $F$. We may assume that there exist $a_1,\dotsc,a_k\in\mathcal{A}$ such that 
\begin{equation*}
U=\{\omega\in\Omega\mid a_i(\omega)>0\text{ for }i=1,\dotsc,k\}.
\end{equation*}
Then $f=\sum_{i=1}^k {a_i}_{+}$ satisfies the requirements.
\end{proof}

Let $p\colon\Omega\to\mathbb{R}$ be a proper, non-negative function on $\Omega$. Let $\xi\colon [0,\infty)\to [1,\infty)$ be a convex, non-decreasing function.

Below we shall consider the space $\mathcal{P}_{\xi(p)}(\Omega)$ of all Radon probability measures $\mu$ on $\Omega$ with 
\begin{equation*}
\int_{\Omega}\xi(p)\, d\mu<\infty,
\end{equation*}
equipped with the weak topology induced by the space $\mathcal{C}_{\xi(p),0}(\Omega)$ of continuous functions $h$ such that for each $\epsilon>0$ there is a compact set $K\subset\Omega$ such that
\begin{equation*}
\sup\Big\{\frac{\abs{h(\omega)}}{\xi(p(\omega))}\mid \omega\notin K\Big\}<\epsilon.
\end{equation*}

\begin{lemma}\label{lem:polish}
Suppose that $\Omega$ is separable in $\tau(\mathcal{G})$ and $p\colon\Omega\to \mathbb{R}$ is proper. Then $\mathcal{P}_{\xi(p)}(\Omega)$ is a Polish space.
\end{lemma}
\begin{proof}
The countable dense set in $\mathcal{P}_{\xi(p)}(\Omega)$ is the rational  convex hull of the Dirac measures at the points of the countable dense set of $\Omega$, as follows by the  Hahn--Banach theorem. 

Let us now show that $\mathcal{P}_{\xi(p)}(\Omega)$ is metrisable and complete in its metric.
As $p$ is a proper function, its sublevel sets are compact. By the Stone--Weierstrass theorem the space of continuous functions on a compact, separable and completely regular space -- cf. Lemma \ref{lem:regular} -- is separable. Let 
\begin{equation*}
K_k=p^{-1}([0,k])\text{ for }k=1,2,\dotsc,
\end{equation*} 
and pick a countable set $(\phi^k_n)_{n=1}^{\infty}$, dense in $\mathcal{C}(K_k)$. Lemma \ref{lem:regular} tells us that these functions can be taken from the vector lattice generated by $\mathcal{G}$. In particular, we may assume that they are defined on $\Omega$ and belong to $\mathcal{C}_{p+1}(\Omega)$. For $\eta_1,\eta_2\in\mathcal{P}_{\xi(p)}(\Omega)$ define
\begin{equation*}
d(\eta_1,\eta_2)=\sum_{k,n=1}^{\infty}\frac{1}{2^{k+n}}\frac{\Big\lvert\int_{\Omega}\phi^k_n\, d\eta_1-\int_{\Omega}\phi^k_n\, d\eta_2\Big\rvert}{1+\Big\lvert\int_{\Omega}\phi^k_n\, d\eta_1-\int_{\Omega}\phi^k_n\, d\eta_2\Big\rvert}.
\end{equation*}
Then the density on appropriate compactae and  Lemma \ref{lem:depe} show that this metric gives the weak topology induced by $\mathcal{C}_{\xi(p),0}(\Omega)$. 
To prove completeness, let us take a Cauchy sequence 
\begin{equation*}
(\eta_k)_{k=1}^{\infty}\subset\mathcal{P}_{\xi(p)}(\Omega)
\end{equation*}
of Radon probability measures on $\Omega$. Each of the elements of the sequence induces a continuous functional on the Banach space $\mathcal{C}_{\xi(p),0}(\Omega)$. By the Banach--Steinhaus uniform boundedness principle it follows that 
\begin{equation*}
\sup\Big\{\int_{\Omega}\xi(p)\, d\eta_k\mid k=1,2,\dotsc\Big\}<\infty.
\end{equation*}
It follows that the sequence $(\eta_k)_{k=1}^{\infty}$ is uniformly tight. Applying the Banach--Alaoglu theorem and \cite[Theorem 7.10.6., p. 112]{Bogachev20072}, similarly to the application in the proof of \cite[Theorem 8.6.7., p. 206]{Bogachev20072}, yields that the sequence converges to some Radon probability measure $\eta\in\mathcal{P}_{\xi(p)}(\Omega)$. 
\end{proof}

\begin{remark}\label{rem:standard}
Since $\mathcal{G}\subset\mathcal{C}_{\xi(p),0}(\Omega)$, we see, by the above proof, that separability and $\sigma$-compactness of $\Omega$ imply separability of $\mathcal{G}$.

 It follows that $\tau(\mathcal{G})$ is countably generated. 
Therefore, by \cite[Proposition 12.1., p. 73]{Kechris1995}, $\Omega$ is a standard Borel space, i.e., it is isomorphic to a subset of a metrisable and separable space, \cite[Definition 12.5., p. 74]{Kechris1995}, \cite[Theorem 13.1., p. 82]{Kechris1995}.
\end{remark}

\section{$\mathcal{F}$-transports and their modifications}\label{s:modifications}

Below we suppose that $\Omega$ is a set. Let $\mathcal{G}$ be a convex cone of functions on $\Omega$, generating a lattice cone $\mathcal{F}$, and $p$ be a proper, non-negative function that belongs to the complete lattice cone $\mathcal{H}$ generated by $\mathcal{G}$. We assume that functions in $\mathcal{G}$ are of $(p+1)$-growth.

\begin{definition}
Suppose that $\mu,\nu\in\mathcal{P}_{\xi(p)}(\Omega)$ are two Radon probability measures.
Any Radon probability measure $\pi$ on $\Omega\times\Omega$ with marginals $\mu,\nu$ and such that for all $f\in\mathcal{F}$ 
and all bounded, non-negative measurable $h$ we have
\begin{equation*}
\int_{\Omega}f(\omega_1)h(\omega_1)\, d\mu(\omega_1)\leq\int_{\Omega\times\Omega}f(\omega_2)h(\omega_1)\, d\pi(\omega_1,\omega_2)
\end{equation*}
we shall call an $\mathcal{F}$-transport between $\mu$ and $\nu$. 
We shall denote the set of all such measures by $\Gamma_{\mathcal{F}}(\mu,\nu)$. 
\end{definition}

\begin{remark}
 If $\mu,\nu\in\mathcal{P}_p(\Omega)$,  then Theorem \ref{thm:poussin} shows that $\mu,\nu\in\mathcal{P}_{\xi(p)}(\Omega)$ for some increasing, convex $\xi\colon\mathbb{R}\to [1,\infty)$ with superlinear growth, so that $\xi(p)\in\mathcal{H}$ is proper. We shall be applying  
Theorem \ref{thm:mart} with the function $\xi(p)$ in lieu of $p$. Theorem \ref{thm:mart} tells us that $\Gamma_{\mathcal{F}}(\mu,\nu)\neq\emptyset$ if and only if $\mu\prec_{\mathcal{F}}\nu$. 
\end{remark}

\subsection{Modifications of $\mathcal{F}$-transports}

Let us recall that $\mathcal{A}=\mathcal{G}-\mathcal{G}$.

\begin{lemma}\label{lem:equivalence}
 Let $\rho\in\mathcal{P}(\Omega\times\Omega)$ be a Radon probability measure with marginals $\mu,\nu$. Suppose that 
\begin{equation*}
    \int_{\Omega\times\Omega} \big(\xi(p)(\omega_1)+\xi(p)(\omega_2)\big)\, d\rho(\omega_1,\omega_2)<\infty.
\end{equation*}
Assume that
\begin{enumerate}
\item\label{i:equationa} for all $g\in\mathcal{G}$ and all bounded Borel measurable functions $h$
\begin{equation*}
\int_{\Omega\times\Omega} h(\omega_1)g(\omega_1)\, d\rho(\omega_1,\omega_2)=\int_{\Omega\times\Omega} h(\omega_1)g(\omega_2)\, d\rho(\omega_1,\omega_2).
\end{equation*}
\end{enumerate}
Then
\begin{enumerate}
\item[(ii)]\label{i:ineqef} for all $f\in\mathcal{F}$ and all bounded Borel measurable, non-negative functions $h$
\begin{equation*}
\int_{\Omega\times\Omega} h(\omega_1)f(\omega_1)\, d\rho(\omega_1,\omega_2)\leq\int_{\Omega\times\Omega} h(\omega_1)f(\omega_2)\, d\rho(\omega_1,\omega_2).
\end{equation*}
\end{enumerate}
Conversely, if \ref{i:ineqef} holds true and $\mathcal{A}\subset\mathcal{F}$, then \ref{i:equationa} holds true as well.
\end{lemma}
\begin{proof}
Suppose that \ref{i:equationa} holds true. Note that it suffices to prove \ref{i:ineqef} for $h=\mathbf{1}_A$, where $A\subset\Omega$ is a Borel set. Observe that the restriction of a Radon measure to a Borel set is a Radon measure. Any $f\in\mathcal{F}$ is a maximum of a finite family of functions in $\mathcal{G}$, by Lemma \ref{lem:supgenfin}. Let 
\begin{equation*}
f=\max\{g_i\mid i=1,2\dotsc,k\}.
\end{equation*}
For $i=1,2,\dotsc,k$ define inductively
\begin{equation*}
A_i=\{\omega\in A\mid f(\omega)=g_i(\omega)\}\setminus \bigcup_{j=1}^{i-1}A_j.
\end{equation*}
Then $(A_i)_{i=1}^k$ is a pairwise disjoint family of Borel sets, whose union is $A$. Applying \ref{i:equationa} on each of the sets $A_i$ yields
\begin{equation*}
\int_{A_i\times\Omega}g_i(\omega_1)\, d\rho(\omega_1,\omega_2)=\int_{A_i\times\Omega}g_i(\omega_2)\, d\rho(\omega_1,\omega_2).
\end{equation*}
Since $f=g_i$ on $A_i$ and $g_i\leq f$ on $\Omega$,
\begin{equation*}
\int_{A_i\times\Omega}f(\omega_1)\, d\rho(\omega_1,\omega_2)\leq \int_{A_i\times\Omega}f(\omega_2)\, d\rho(\omega_1,\omega_2).
\end{equation*}
Adding up the resulting inequalities yields the assertion.

The converse implication is trivial.
\end{proof}

\begin{proposition}\label{pro:mod}
Let $\mu,\nu$ be two Radon probability measures on $\Omega$. Suppose that $\pi\in\Gamma_{\mathcal{F}}(\mu,\nu)$ is an $\mathcal{F}$-transport between $\mu$ and $\nu$.  Let $\rho\in\mathcal{M}(\Omega\times\Omega)$ be a finite signed Radon measure such that
\begin{enumerate}
\item $ \int_{\Omega\times\Omega} \big(\xi(p)(\omega_1)+\xi(p)(\omega_2)\big)\, d\abs{\rho}(\omega_1,\omega_2)<\infty$,
\item the marginals ${p_1}_{\#}\rho$, ${p_2}_{\#}\rho$ of $\rho$ are zero measures,
\item there exists $\epsilon>0$ such that $ 0\leq\epsilon \rho+\pi$,
\item $\int_{\Omega\times\Omega} h(\omega_1)a(\omega_2)\, d\rho(\omega_1,\omega_2)=0$ for all $a\in\mathcal{A}$ and all bounded measurable $h$.
\end{enumerate}
Then $\pi+\epsilon\rho\in\Gamma_{\mathcal{F}}(\mu,\nu)$.
\end{proposition}
\begin{proof}
Note that since ${p_1}_{\#}\rho=0$, then
\begin{equation*}
\int_{\Omega\times\Omega} h(\omega_1)a(\omega_1)\, d\rho(\omega_1,\omega_2)=0\text{ for all }a\in\mathcal{A}\text{ and bounded measurable }h.
\end{equation*}
The conclusion follows now readily from Lemma \ref{lem:equivalence}, as $\mathcal{A}=\mathcal{G}-\mathcal{G}$.
\end{proof}

\section{Maximal disintegrations}\label{s:maxx}

\subsection{Joint support}\label{s:joint}

In the classical optimal transport theory, it can be shown that for any two Radon probability measures there exists a transport between them whose support contains the supports of all other transports between the two measures. 

A similar condition was studied in the martingale optimal transport theory, see \cite[Theorem 2.1]{Touzi2019}. In there, disintegrations with maximal closed convex hulls of their supports are studied.

\begin{proposition}\label{pro:support}
Let $\mu\prec_{\mathcal{F}}\nu$ be two Radon probability measures. Then there exists the minimal closed set $S\subset\Omega\times\Omega$ such that $\pi(S)=1$ for every $\pi\in\Gamma_{\mathcal{F}}(\mu,\nu)$.
\end{proposition}
\begin{proof}
Let $\mathcal{S}$ be the family of all closed sets $S$ such that $\pi(S)=1$ for every $\pi\in\Gamma_{\mathcal{F}}(\mu,\nu)$. Clearly, $\Omega\times\Omega\in\mathcal{S}$, so it is non-empty.
Let $S_0$ be the intersection of all elements of $\mathcal{S}$. It is a closed set. We shall show that it belongs to $\mathcal{S}$.

Any measure in $\Gamma_{\mathcal{F}}(\mu,\nu)$ is a Radon measure. Any Radon measure is $\tau$-additive, see \cite[Definition 7.2.1., p. 73 and Proposition 7.2.2., p. 74]{Bogachev20072}.

For any $\pi\in\Gamma_{\mathcal{F}}(\Omega\times\Omega)$, the union of the family $(\Omega\times\Omega\setminus S)_{S\in\mathcal{S}}$ of open sets of $\pi$-measure zero has $\pi$-measure zero.
The proof is complete.
\end{proof}

\begin{definition}\label{def:joint}
The minimal closed set $S\subset\Omega\times\Omega$ such that $\pi(S)=1$ for all $\pi\in\Gamma_{\mathcal{F}}(\mu,\nu)$ we shall call the joint support of $\Gamma_{\mathcal{F}}(\mu,\nu)$.
\end{definition}

We shall see now that the joint support of $\Gamma_{\mathcal{F}}(\mu,\nu)$ allows us to characterise polar sets under additional assumption.

\begin{remark}
As we shall see later in Section \ref{s:maximaldis}, the joint support is not sufficient for our purposes and therefore we shall employ another set with certain minimality properties.
\end{remark}

\begin{proposition}\label{pro:supported}
Suppose that $\Omega$ is separable. Let $\mu\prec_{\mathcal{F}}\nu$ be two Radon measures in $\mathcal{F}$-order. Let $S\subset\Omega\times\Omega$ be the joint support of $\Gamma_{\mathcal{F}}(\Omega\times\Omega)$. There exists $\pi\in \Gamma_{\mathcal{F}}(\mu,\nu)$ such that $S$ is the support of $\pi$.
\end{proposition}
\begin{proof}
Since $\Omega$ is separable, so is $\Gamma_{\mathcal{F}}(\mu,\nu)$ considered in the weak* topology. Let $(\pi_i)_{i=1}^{\infty}\subset \Gamma_{\mathcal{F}}(\mu,\nu)$ be a dense set. 

Define 
\begin{equation*}
\pi_0=\sum_{i=1}^{\infty}\frac1{2^i}\pi_i.
\end{equation*}
Then $\pi_0\in\Gamma_{\mathcal{F}}(\mu,\nu)$ and $S$ is the support of $\pi_0$. 

Indeed, suppose on the contrary, that there is an open set $U\subset\Omega\times\Omega$ such that $U\cap S\neq\emptyset$ and such that $\pi_0(U)=0$. We may assume that
\begin{equation*}
U=\{(\omega_1,\omega_2)\in\Omega\times\Omega\mid a_i(\omega_1)>0,b_i(\omega_2)>0\text{ for }i=1,2,\dotsc,k\}
\end{equation*}
for some functions $a_1,\dotsc,a_k,b_1,\dotsc,b_k\in\mathcal{A}$. Let 
\begin{equation*}
f=\min\{\prod_{i=1}^k{a_i}_+{b_i}_+,1\}.
\end{equation*}
Then $f$ is continuous, bounded, non-negative, positive on $U$, zero on the complement of $U$.
Since $\pi_0(U)=0$, also $\int_Uf\, d\pi_i=0$ for $i=1,2,\dotsc$.
The map 
\begin{equation*}
\Gamma_{\mathcal{F}}(\mu,\nu)\ni \pi\mapsto \int_{\Omega\times\Omega}f\, d\pi\in\mathbb{R}
\end{equation*}
is continuous. Since $(\pi_i)_{i=1}^{\infty}$ is dense, it follows that $\int_U f\, d\pi=0$ for all $\pi\in\Gamma_{\mathcal{F}}(\mu,\nu)$.
This stands in contradiction with the definition of $S$.
\end{proof}

\subsection{Disintegrations}\label{s:disinteg}

So far we have managed to not employ disintegrations of measures in our considerations. Let us recall that Souslin spaces are continuous images of Polish spaces. 
We shall employ regular conditional measures, see \cite[Definition 10.4.2., p. 358]{Bogachev20072}. In the following theorem, we combine \cite[Example 10.4.11. p. 363]{Bogachev20072} and \cite[Lemma 10.4.3, p. 358]{Bogachev20072}

\begin{theorem}\label{thm:dis}
Let $\mu$ be a probability measure on a Souslin space $X$ equipped with the Borel $\sigma$-algebra $\sigma(\tau)$. Let $r$ be a  $\mu$-measurable map to a Souslin space $Y$, equipped with the Borel $\sigma$-algebra $\mathcal{E}$. Then there exists a map 
\begin{equation*}
Y\times\sigma(\tau)\ni (y,B)\mapsto \mu(y,B)\in\mathbb{R}
\end{equation*}
such that
\begin{enumerate}
\item\label{i:radon} for every $y\in Y$,$ \sigma(\tau)\ni B\mapsto \mu(y,B)\in\mathbb{R}$ is a Borel probability measure on $X$,
\item for every $B\in\sigma(\tau)$ the map $Y\ni y\mapsto \mu(y,B)\in\mathbb{R}$ is Borel measurable,
\item\label{i:integral} for all $B\in\sigma(\tau)$ and all $E\in\mathcal{E}$, 
\begin{equation*}
\mu\big(B\cap r^{-1}(E)\big)=\int_E\mu(y,B)\, dr_{\#}\mu(y),
\end{equation*}
\item for $r_{\#}\mu$-almost every $y\in Y$, $\mu(y,\cdot)$ is concentrated on $r^{-1}(\{y\})$,
\item if $\sigma(\tau)$ is countably generated, then the considered map is essentially unique: if $\mu_1,\mu_2\colon Y\times \sigma(\tau)\to\mathbb{R}$ are two maps satisfying conditions \ref{i:radon}-\ref{i:integral}, then there exists a Borel set $Z\subset Y$ with $\mu(Z)=0$ such that for all $y\in Y\setminus Z$ and all Borel sets $B\in\sigma(\tau)$ 
\begin{equation*}
\mu_1(y,B)=\mu_2(y,B).
\end{equation*}
\end{enumerate}
\end{theorem}

We shall now assume, as we do in Theorem \ref{thm:mart}, that there exists a non-negative, proper function $p\in\mathcal{H}$ such that $\mathcal{G}\subset\mathcal{D}_{p+1}(\Omega)$ and that
\begin{equation*}
\int_{\Omega}p\, d\mu<\infty, \int_{\Omega}p\, d\nu<\infty.
\end{equation*}
By Theorem \ref{thm:poussin}, there exists a continuous, increasing, convex function $\xi\colon\mathbb{R}\to [1,\infty)$ with superlinear growth such that
\begin{equation*}
\int_{\Omega}\xi(p)\, d\mu<\infty, \int_{\Omega}\xi(p)\, d\nu<\infty.
\end{equation*}

\begin{remark}\label{rem:integral}
Suppose that $\pi\in\Gamma_{\mathcal{F}}(\mu,\nu)$ and that $p$ is a Borel measurable function integrable with respect to $\nu$. We shall be assuming that $\Omega$ is separable and $\sigma$-compact, so that $\mathcal{P}_{\xi(p)}(\Omega)$ is a Polish space, by Lemma \ref{lem:polish}. Theorem \ref{thm:dis} shows that there exists a disintegration of $\pi$ with respect to the projection $p_1$ on the first co-ordinate and any such disintegration $\lambda$ of $\pi$ is concentrated on the fibres, i.e.,
\begin{equation*}
\lambda(\omega,\{\omega\}^c\times\Omega)=0\text{ for }\mu\text{-almost every }\omega\in\Omega.
\end{equation*}
We shall in what follows identify $\lambda(\omega,\cdot)$ with the measure on the appropriate fibre given by
\begin{equation*}
\sigma(\tau(\mathcal{A}))\ni B\mapsto \lambda(\omega,\{\omega\}\times B)\in\mathbb{R}.
\end{equation*}
Note that if  $\lambda$ is a disintegration of $\pi$ with respect to the projection $p_1$ on the first co-ordinate. Then, for $\mu$-almost every $\omega\in\Omega$, $\xi(p)$ is integrable with respect to $\lambda(\omega,\cdot)$, so that $\lambda(\omega,\cdot)\in\mathcal{P}_{\xi(p)}(\Omega)$.
\end{remark}

\begin{definition}
Let $\pi\in\Gamma_{\mathcal{F}}(\mu,\nu)$. Let $T(\pi)$ be the set of all disintegrations
\begin{equation*}
\lambda(\cdot,\cdot)\colon  \Omega\times \sigma(\tau(\mathcal{A}))\to\mathbb{R}
\end{equation*}
of $\pi$ with respect to the projection $p_1$ on the first co-ordinate. That is,
\begin{enumerate}
\item for every $\omega\in \Omega$, $ \sigma(\tau(\mathcal{A}))\ni B\mapsto \lambda(\omega,B)\in\mathbb{R}$ is a Radon probability measure on $\Omega$ in $\mathcal{P}_{\xi(p)}(\Omega)$,
\item for every $B\in\sigma(\tau(\mathcal{A}))$ the map $\Omega\ni \omega\mapsto\lambda(\omega,B)\in\mathbb{R}$ is Borel measurable,
\item for all Borel sets $B\subset\Omega\times\Omega$ and $E\subset\Omega$, 
\begin{equation*}
\pi\big(B\cap r^{-1}(E)\big)=\int_E\lambda(\omega,B_{\omega})\, d\mu(\omega),
\end{equation*}
where 
\begin{equation*}
B_{\omega}=\{\omega'\in\Omega\mid (\omega,\omega')\in\Omega\}.
\end{equation*}

\end{enumerate}
\end{definition}

\begin{remark}
Theorem \ref{thm:dis} and Remark \ref{rem:integral} tell us that for any $\pi\in\Gamma_{\mathcal{F}}(\mu,\nu)$, $T(\pi)\neq\emptyset$. Indeed, by \cite[Theorem 6.6.6., p. 21]{Bogachev20072} and Remark \ref{rem:standard}, $\Omega$ is a Souslin space, and also $\Omega\times\Omega$ is a Souslin space. Moreover, $\lambda(\omega,\cdot)$ for $\omega\in\Omega$ can be taken to be Radon measures, as $\Omega$ is a countable union of compact metrisable sets, and any Borel probability measure on a Polish space is a Radon measure.
\end{remark}

\begin{definition}
We shall denote by $\Lambda_{\mathcal{F}}(\mu,\nu)$ the set of all maps $\lambda\in T(\pi)$ for some $\pi\in \Gamma_{\mathcal{F}}(\mu,\nu)$. 
\end{definition}

\subsection{Maximal disintegrations}\label{s:maximaldis}

We shall assume that $\Omega$ is separable with respect to $\tau(\mathcal{G})$. We recall that the existence of a proper function $p\in\mathcal{H}$ implies that $\Omega$ is $\sigma$-compact. The aim of this section is to introduce a notion of maximal disintegration of $\mathcal{F}$-transports. 

\begin{definition}\label{def:maxdis}
A map $\lambda\in \Lambda_{\mathcal{F}}(\mu,\nu)$ such that for any $\lambda'\in \Lambda_{\mathcal{F}}(\mu,\nu)$ there is
\begin{equation*}
\mathrm{supp}\lambda'(\omega,\cdot)\subset\mathrm{supp}\lambda(\omega,\cdot)\text{ for }\mu\text{-almost every }\omega\in \Omega,
\end{equation*}
we shall call a maximal disintegration of $\Gamma_{\mathcal{F}}(\mu,\nu)$ with respect to $p_1$, or, briefly, a maximal disintegration of  $\Gamma_{\mathcal{F}}(\mu,\nu)$. Note here that we consider the supports of measures with respect to the topology $\tau(\mathcal{G})$.
\end{definition}

\begin{proposition}\label{pro:maxsup}
Suppose that $\Omega$ is separable with respect to $\tau(\mathcal{G})$. Then there exists a maximal disintegration of $\Gamma_{\mathcal{F}}(\mu,\nu)$ with respect to $p_1$.
\end{proposition}
\begin{proof}
Since $\Omega$ is separable and $\sigma$-compact, $\mathcal{C}_0(\Omega)$ is separable as well, cf. Remark \ref{rem:standard},  and thus also $L^1(\Omega,\mathcal{C}_0(\Omega),\mu)$ is separable. 
Recall that Remark \ref{rem:el} tells us that the space 
\begin{equation*}
L^{\infty}(\Omega,\mathcal{C}_0(\Omega)^*,\mu)
\end{equation*}
of all weakly* measurable essentially bounded functions is dual to 
\begin{equation*}
L^1(\Omega,\mathcal{C}_0(\Omega),\mu).
\end{equation*}
By the Banach--Alaoglu theorem the unit ball of the former space is weakly* compact. The separability of the latter space implies that the unit ball of the former is metrisable. Thus it is also separable, in the weak* topology.
Now, as $\Lambda_{\mathcal{F}}(\mu,\nu)\subset L^{\infty}(\Omega,\mathcal{C}_0(\Omega)^*,\mu)$ is contained in the unit ball, it is separable.

Let $(\lambda_i)_{i=1}^{\infty}\subset \Lambda_{\mathcal{F}}(\mu,\nu)$ be a dense set. 
For each $i$ let $\lambda_i$ is a disintegration of some $\pi_i\in \Gamma_{\mathcal{F}}(\mu,\nu)$, $\lambda_i\in T(\pi_i)$. Set
\begin{equation*}
\lambda=\sum_{i=1}^{\infty}\frac1{2^i}\lambda_i.
\end{equation*}
Then $\lambda\in T(\pi)$, where
\begin{equation*}
\pi=\sum_{i=1}^{\infty}\frac1{2^i}\pi_i.
\end{equation*}
Let now $U\subset\Omega$ be open. 
As in Lemma \ref{lem:regular}, we may assume that 
\begin{equation*}
U=\{\omega\in\Omega\mid a_i(\omega)>0\text{ for }i=1,2,\dotsc,k\}
\end{equation*}
for some $a_1,\dotsc,a_k\in\mathcal{A}$. Let $f_U=\min\{\max\{\prod_{i=1}^k {a_i},0\},1\}$. Then $f_U\geq 0$ and is positive exactly on $U$.
Note that for $\omega\in \Omega$, 
\begin{equation*}
\mathrm{supp}\lambda(\omega,\cdot)\cap U=\emptyset\text{ if and only if }\lambda(\omega,U)=0.
\end{equation*} 
Therefore the set 
\begin{equation*}
B_U=\{\omega\in\Omega\mid \mathrm{supp}\lambda(\omega,\cdot)\cap U=\emptyset\}
\end{equation*}
is Borel.
Now, by construction,
\begin{equation*}
\int_{\Omega}f_U\, d\lambda(\omega,\cdot)=0\text{ for }\mu\text{-almost every }\omega\in B_U.
\end{equation*}
This implies that the same equality holds true with $\lambda_i$ in place of $\lambda$, for each number $i=1,2,\dotsc$. By density of $(\lambda_i)_{i=1}^{\infty}$ for any $\lambda'\in \Lambda_{\mathcal{F}}(\mu,\nu)$ there is  
\begin{equation*}
\int_{\Omega}f_U\, d\lambda'(\omega,\cdot)=0\text{ for }\mu\text{-almost every }\omega\in B_U.
\end{equation*}
It follows that 
\begin{equation}\label{eqn:inclusion}
\mathrm{supp}\lambda'(\omega,\cdot)\subset \Omega\setminus U \text{ for }\mu\text{-almost every }\omega\in B_U
\end{equation}
For each open set $U$, let $C_U$ denote the set of $\omega\in\Omega$ such that (\ref{eqn:inclusion}) holds true or $\omega\notin B_U$. Then $\mu(C_U)=1$ and
\begin{equation*}
\mathrm{supp}\lambda(\omega,\cdot)\subset \Omega\setminus U\text{ implies that }\mathrm{supp}\lambda'(\omega,\cdot)\subset \Omega\setminus U \text{ for }\mu\text{-almost every }\omega\in C_U.
\end{equation*}
We shall now employ separability and $\sigma$-compactness of $\Omega$. By Remark \ref{rem:standard}, $\Omega$ is metrisable on compactae. Hence, it possesses a countable basis of neighbourhoods on each of its compact sets.
Let $C$ be the intersection of the sets $C_U$ corresponding to the countable basis for one of the countable number of compact subsets of $\Omega$. Then $\mu(C)=1$ and $C$ is Borel. Moreover for $\omega\in C$ and for all open sets $U\subset\Omega$
\begin{equation*}
\mathrm{supp}\lambda(\omega,\cdot)\subset \Omega\setminus U\text{ implies that }\mathrm{supp}\lambda'(\omega,\cdot)\subset \Omega\setminus U.
\end{equation*}
Taking for $\omega\in \Omega$, $U=\mathrm{supp}\lambda(\omega,\cdot)^c$ completes the proof.
\end{proof}

\begin{remark}
It is readily visible that the supports of a maximal disintegration of $\Gamma_{\mathcal{F}}(\mu,\nu)$ are unique, up to a set of $\mu$-measure zero.
\end{remark}

\begin{remark}
The proof shows that one may pick a maximal disintegration $\lambda\in\Lambda_{\mathcal{F}}(\mu,\nu)$ in the quasi-relative interior of the set $\Lambda_{\mathcal{F}}(\mu,\nu)$, cf. \cite[Proposition  1.2.9, p. 18]{Zalinescu2002}. 
\end{remark}

\begin{example}\label{exa:jointbad}
Let us present an example of a measure $\pi\in \mathcal{P}([0,1]\times [0,1])$, $\pi\in\Gamma(\mu,\nu)$ for some $\mu,\nu\in\mathcal{P}([0,1])$, for which the sections of $\mathrm{supp}\pi$ are not equal $\mu$-almost everywhere to supports of a maximal disintegration.
Let 
\begin{equation*}
\mu=\frac12\delta_0+\sum_{i=1}^{\infty}\frac1{2^{i+1}}\delta_{1/i}
\text{ and }\pi=\frac12\delta_{(0,0)}+\sum_{i=1}^{\infty}\frac1{2^{i+1}}\delta_{1/i}\otimes \lambda,
\end{equation*}
where $\lambda$ is the Lebesgue measure.
Then
\begin{equation*}
\lambda(1/i,\cdot)=\lambda\text{ for }i=1,2,\dotsc\text{ and }\lambda(0,\cdot)=\delta_0
\end{equation*}
is a maximal disintegration of $\pi$.
However 
\begin{equation*}
\mathrm{supp}\pi=\bigcup_{i=1}^{\infty}\{1/i\}\times [0,1]\cup \{0\}\times [0,1].
\end{equation*}
Thus the $0$-section of the support does not coincide with the support of $\lambda(0,\cdot)$.
\end{example}

\section{Measurable selections}\label{s:selections}
\begin{definition}
For $\lambda\in \Lambda_{\mathcal{F}}(\mu,\nu)$ and $\omega\in\Omega$ we define $\Xi(\lambda,\omega)$ to be the set of all Radon probability measures $\eta\in\mathcal{P}_{\xi(p)}(\Omega)$ such that 
\begin{equation*}
C\lambda(\omega,\cdot)\geq \eta \geq 0 \text{ for some }C\in (0,\infty).
\end{equation*}
We define the set $\Theta(\lambda,\omega)$ to be the set of all continuous functionals $a^*\in\mathcal{A}^*$ such that
\begin{equation}\label{eqn:functional}
a^*(a)=\int_{\Omega}a\, d\eta\text{ for all }a\in\mathcal{A},
\end{equation}
where $\eta\in \Xi(\lambda,\omega)$.
We shall denote by $\zeta_{\lambda,\omega}$ the map $\zeta_{\lambda,\omega}\colon \Xi(\lambda,\omega)\to \Theta(\lambda,\omega)$ that to a measure $\eta$ associates the linear functional $a^*$ given by (\ref{eqn:functional}). .
\end{definition}

\begin{remark}
The Radon--Nikodym theorem, \cite[Theorem 3.2.2., p. 178]{Bogachev20071}, tells us that for each $\eta\in\Xi(\lambda,\omega)$ there exists non-negative $f\in L^{\infty}(\Omega,\lambda(\omega,\cdot))$ such that for each Borel measurable, bounded function $g$
\begin{equation*}
\int_{\Omega}g\, d\eta=\int_{\Omega} g f\, d\lambda(\omega,\cdot).
\end{equation*}
\end{remark}

\begin{definition}
We shall say that a Radon measure $\gamma$ on $\Omega\times\Omega$ is symmetric provided that
\begin{equation*}
\gamma=s_{\#}\gamma,
\end{equation*}
where $s\colon\Omega\times\Omega\to\Omega\times\Omega$ is given by the formula $s(\omega_1,\omega_2)=(\omega_2,\omega_1)$.
\end{definition}

A tool that will be useful is the Kuratowski--Ryll-Nardzewski measurable selection theorem \cite{Kuratowski1965}, see also \cite[Theorem 6.9.3., p. 36]{Bogachev20072} together with its corollary \cite[Corollary 6.9.4., p. 36]{Bogachev20072}. We cite them below together.

\begin{theorem}\label{thm:kuratowski}
Suppose that $Y$ is a Polish space and that $(X,\Sigma)$ is a measurable space. Let $S$ be a map on $X$ that takes values in the closed, non-empty subsets of $Y$. Suppose that for any open set $U\subset Y$ there is
\begin{equation*}
\{x\in X\mid S(x)\cap U\neq \emptyset\}\in \Sigma.
\end{equation*}
Then there exists a map $s\colon X\to Y$, measurable with respect to $\Sigma$ and Borel $\sigma$-algebra on $Y$, such that 
\begin{equation*}
s(x)\in S(x)\text{ for any }x\in X.
\end{equation*}
Moreover, one can find a sequence $(s_n)_{n=1}^{\infty}$ of measurable maps as above such that for each $x\in X$, the set $\{s_n(x)\mid n=1,2,\dotsc,\}$ is dense in $S(x)$.
\end{theorem}

In other words, the theorem tells that if a multifunction with values in closed sets is measurable, then it admits a measurable selection.

\begin{lemma}\label{lem:selection}
Suppose that $\Omega$ is separable in  $\tau(\mathcal{G})$. Let $\pi\in\Gamma_{\mathcal{F}}(\mu,\nu)$ and let $\lambda\in T(\pi)$. Let 
\begin{equation*}
\Lambda=\{(\omega_1,\omega_2)\in\Omega\times\Omega\mid \Theta(\lambda,\omega_1)\cap\Theta(\lambda,\omega_2)\neq \emptyset\} 
\end{equation*} 
and
\begin{equation*}
\Delta=\{(\omega_1,\omega_2)\in \Omega\times\Omega\mid \omega_1=\omega_2\}.
\end{equation*}
Then $\Lambda$ is Borel measurable. Moreover, there exists a sequence of Borel measurable functions
\begin{equation*}
\Lambda\setminus \Delta\ni (\omega_1,\omega_2)\mapsto \eta^k_{\omega_1,\omega_2}\in \mathcal{P}_{\xi(p)}(\Omega)
\end{equation*}
such that 
\begin{enumerate}
\item for $k=1,2,\dotsc$ there exists symmetric Borel measurable $\epsilon^k\colon\Lambda\setminus \Delta\to (0,\infty)$ such that
for all $(\omega_1,\omega_2)\in\Lambda\setminus \Delta$
\begin{equation*}
\eta^k_{\omega_1,\omega_2}\leq \epsilon^k(\omega_1,\omega_2)\lambda(\omega_1,\cdot),
\end{equation*}
\item for $k=1,2,\dotsc$ and for all $(\omega_1,\omega_2)\in \Lambda\setminus \Delta$ the functionals
\begin{equation*}
\mathcal{A}\ni a\mapsto\int_{\Omega}a\, d\eta^k_{\omega_1,\omega_2}\in\mathbb{R}\text{ and }\mathcal{A}\ni a\mapsto\int_{\Omega}a\, d\eta^k_{\omega_2,\omega_1}\in\mathbb{R}
\end{equation*}
are equal and  belong to $\Theta(\lambda,\omega_1)\cap\Theta(\lambda,\omega_2)$,

\item for each $(\omega_1,\omega_2)\in\Lambda\setminus\Delta$ the set
\begin{equation*}
\{(\eta^k_{\omega_1,\omega_2},\eta^k_{\omega_2,\omega_1})\in\mathcal{P}_{\xi(p)}(\Omega)\times \mathcal{P}_{\xi(p)}(\Omega)\mid k=1,2,\dotsc\}
\end{equation*}
is dense in the set of probabilities in $\Xi(\lambda,\omega_1)\times\Xi(\lambda,\omega_2)$ for which the respective functionals are equal and belong to $\Theta(\lambda,\omega_1)$ and $\Theta(\lambda,\omega_2)$ respectively.
\end{enumerate}
\end{lemma}
\begin{proof}
Since $\Omega$ is $\sigma$-compact, we may assume that it is compact and therefore, by Remark \ref{rem:standard}, it is also a Polish space.
For $n\in\mathbb{N}$ and $\omega\in\Omega$, let $\Theta_n(\lambda,\omega)$ denote the set of all functionals $a^*\in\mathcal{A}^*$ of the form 
\begin{equation*}
a^*(a)=\int_{\Omega}a\, d\eta,
\end{equation*}
where $\eta\in\mathcal{P}_{\xi(p)}(\Omega)$ is such that $\eta\leq n\lambda(\omega,\cdot)$.
If 
\begin{equation*}
\Theta(\lambda,\omega_1)\cap\Theta(\lambda,\omega_2)\neq\emptyset
\end{equation*}
then there exists $n\in\mathbb{N}$ such that 
\begin{equation*}
\Theta_n(\lambda,\omega_1)\cap\Theta_n(\lambda,\omega_2)\neq\emptyset.
\end{equation*}
This is to say,
\begin{equation*}
\Lambda=\bigcup_{n=1}^{\infty}\{(\omega_1,\omega_2)\in\Omega\times\Omega\mid \Theta_n(\lambda,\omega_1)\cap\Theta_n(\lambda,\omega_2)\neq \emptyset\}.
\end{equation*}
Thus, it suffices to prove the lemma with $\Theta_n(\lambda,\cdot)$ in place of $\Theta(\lambda,\cdot)$.
For each $(\omega_1,\omega_2)\in\Lambda$, the set $R(\omega_1,\omega_2)$ of all pairs
\begin{equation*}
(\eta_1,\eta_2)\in  \mathcal{P}_{\xi(p)}(\Omega)\times\mathcal{P}_{\xi(p)}(\Omega)
\end{equation*}
such that 
\begin{equation}\label{eqn:compact}
\eta_i\leq n\lambda(\omega_i,\cdot),i=1,2,
\end{equation} 
and for all $a\in\mathcal{A}$
\begin{equation*}
\int_{\Omega}a\, d\eta_1=\int_{\Omega}a\,d\eta_2
\end{equation*}
is closed in the weak topology. We have thus defined a map $R$ on $\Lambda$ with values in closed, convex, non-empty subsets of a Polish space $\mathcal{P}_{\xi(p)}(\Omega)\times\mathcal{P}_{\xi(p)}(\Omega)$. 
Moreover, condition (\ref{eqn:compact}) tells us that the measures in $R(\omega_1,\omega_2)$ are uniformly tight. Therefore, $R(\omega_1,\omega_2)$ is compact for each $(\omega_1,\omega_2)\in \Lambda$, see \cite[Theorem 8.6.7., p. 206]{Bogachev20072}.

Note now that the set
\begin{equation}\label{eqn:set}
\{(\omega,\eta)\in\Omega\times \mathcal{P}_{\xi(p)}(\Omega)\mid \eta\leq n\lambda(\omega,\cdot)\}
\end{equation}
is Borel. Indeed, inequality $\eta\leq n\lambda(\omega,\cdot)$ is equivalent to demanding that for a countable dense subset $\mathcal{G}_0\subset\mathcal{C}_{\xi(p),0}(\Omega)$ of non-negative functions, there is
\begin{equation*}
\int_{\Omega}g\, d\eta\leq n\int_{\Omega}g\, d\lambda(\omega,\cdot)\text{ for all }g\in\mathcal{G}_0.
\end{equation*}
The function
\begin{equation*}
(\omega,g)\mapsto\int_{\Omega}g\, d\lambda(\omega,\cdot)
\end{equation*}
is a Carath\'eodory function, whence it is jointly measurable by \cite[Lemma 4.51., p. 153]{Aliprantis2006} and Lemma \ref{lem:polish}. This shows Borel measurability of (\ref{eqn:set}).

Now \cite[Theorem 28.8, p. 220]{Kechris1995} and Remark \ref{rem:standard} tell us that $R$ is a Borel map provided that the set
\begin{equation*}
\Big\{(\omega_1,\omega_2,\eta_1,\eta_2)\in\Omega\times\Omega\times\mathcal{P}_{\xi(p)}(\Omega)\times\mathcal{P}_{\xi(p)}(\Omega)\mid (\eta_1,\eta_2)\in R(\omega_1,\omega_2)\Big\}
\end{equation*}
is Borel. This however follows readily by (\ref{eqn:set}). 

It follows immediately that $\Lambda$ is Borel measurable.

An application of Theorem \ref{thm:kuratowski} yields the existence of two Borel measurable sequences
\begin{equation*}
\Lambda\ni (\omega_1,\omega_2)\mapsto (\beta^k_{\omega_1,\omega_2},\zeta^k_{\omega_1,\omega_2})\in R(\omega_1,\omega_2)
\end{equation*}
such that for each $(\omega_1,\omega_2)\in \Lambda$, the set 
\begin{equation*}
\{(\beta^k_{\omega_1,\omega_2},\zeta^k_{\omega_1,\omega_2})\in R(\omega_1,\omega_2)\mid k=1,2,\dotsc\}
\end{equation*}
is dense in $R(\omega_1,\omega_2)$.

Now, we shall write $(\omega_1,\omega_2)\sim (\omega_1',\omega_2')$ if $\omega_1=\omega_2'$ and $\omega_2=\omega_1'$. 
Pick two Borel measurable, disjoint sets $E_1,E_2$ in $\Omega\times\Omega\setminus \Delta$ such that $E_1\cup E_2=\Omega\times\Omega\setminus \Delta$ and $E_i$ contains exactly one element of each equivalence class, $i=1,2$. The existence of such transverals follows e.g. by \cite[Theorem 12.16, p. 78]{Kechris1995}.
Define for $k=1,2,\dotsc$ and $(\omega_1,\omega_2)\in \Omega\times\Omega\setminus \Delta$
\begin{equation*}
\eta_{\omega_1,\omega_2}^k=\beta^k_{\omega_1,\omega_2}\mathbf{1}_{E_1}(\omega_1,\omega_2)+ \xi^k_{\omega_2,\omega_1}\mathbf{1}_{E_2}(\omega_1,\omega_2)
\end{equation*} 
Then the required conditions are satisfied.
\end{proof}

\section{Mixing}\label{s:mixing}

\begin{lemma}\label{lem:modification}
Suppose that $\Omega$ is separable in $\tau(\mathcal{G})$. Let $\pi\in\Gamma_{\mathcal{F}}(\mu,\nu)$ and let $\lambda\in T(\pi)$. Suppose that $\gamma\in\Gamma(\mu,\mu)$ is symmetric. Let 
\begin{equation*}
\Lambda=\{(\omega_1,\omega_2)\in\Omega\times\Omega\mid \Theta(\lambda,\omega_1)\cap\Theta(\lambda,\omega_2)\neq \emptyset\}.
\end{equation*} 
Suppose that
\begin{equation*}
\Lambda\setminus \Delta\ni (\omega_1,\omega_2)\mapsto \eta_{\omega_1,\omega_2}\in\mathcal{P}_{\xi(p)}(\Omega)
\end{equation*}
is a Borel map such that there exists symmetric Borel measurable $\epsilon\colon\Lambda\setminus \Delta\to (0,\infty)$ such that $\eta_{\omega_1,\omega_2}\leq \epsilon(\omega_1,\omega_2)\lambda(\omega_1,\cdot)$ and such that 
\begin{equation*}
\int_{\Omega}a\, d\eta_{\omega_1,\omega_2}=\int_{\Omega}a\, d\eta_{\omega_2,\omega_1}\text{ for all }a\in\mathcal{A}.
\end{equation*}
Then there exists $\pi'\in\Gamma_{\mathcal{F}}(\mu,\nu)$ and $\lambda'\in T(\pi')$ such that if $\theta$ is any disintegration of $\gamma$ with respect to the projection $p_1$ on the first co-ordinate, then the measure
\begin{equation*}
\xi_{\omega}=\frac1{c(\omega)}\int_{\Omega}\eta_{\omega_2,\omega}\frac{\mathbf{1}_{\Lambda\setminus \Delta}(\omega,\omega_2)}{\epsilon(\omega,\omega_2)}\, d\theta(\omega,\omega_2)
\end{equation*}
belongs to $\Xi(\lambda',\omega)$. 
Here
\begin{equation*}
c(\omega)=\int_{\Omega}\frac{\mathbf{1}_{\Lambda\setminus \Delta}(\omega,\omega_2)}{\epsilon(\omega,\omega_2)}\, d\theta(\omega,\omega_2).
\end{equation*}
Moreover, for any $f\in\mathcal{C}_{\xi(p),0}(\Omega)$, the map 
\begin{equation*}
\Omega\ni \omega\mapsto \int_{\Omega}f\,d \xi_{\omega}\in\mathbb{R}
\end{equation*}
is Borel measurable.
\end{lemma}
\begin{proof}
Let us take a map $\eta$ as in the formulation of the lemma.
Define measure $\rho$ via the formula
\begin{equation*}
\rho(A)=\int_{\Omega\times\Omega}\Big(\int_{\Omega}\mathbf{1}_A(\omega_1,\cdot)\, d\eta_{\omega_2,\omega_1}-\int_{\Omega}\mathbf{1}_A(\omega_1,\cdot)\, d\eta_{\omega_1,\omega_2}\Big)\,\frac{\mathbf{1}_{\Lambda\setminus \Delta}(\omega_1,\omega_2)}{\epsilon(\omega_1,\omega_2)} d\gamma(\omega_1,\omega_2)
\end{equation*}
for Borel sets $A\subset\Omega\times\Omega$.
Note that, by the property of $\epsilon$ we see that
\begin{equation*}
\int_{\Omega\times\Omega\times\Omega}\mathbf{1}_A(\omega_1,\cdot)\, d\eta_{\omega_1,\omega_2} \frac{\mathbf{1}_{\Lambda\setminus \Delta}(\omega,\omega_2)}{\epsilon(\omega_1,\omega_2)}\, d\gamma(\omega_1,\omega_2)\leq \int_{\Omega\times\Omega\times\Omega}\mathbf{1}_A(\omega_1,\cdot)\, d\lambda(\omega_1,\cdot)\, d\gamma(\omega_1,\omega_2).
\end{equation*}
As $\gamma\in\Gamma(\mu,\mu)$, and $\lambda\in T(\pi)$, the bound from above is equal to $\pi(A)$. Similarly, the first summand is bounded by $\pi(\Omega)$. 
Thus, by $\rho$ is a well-defined Radon measure, such that $\pi+\rho\geq 0$.

Now, since $\gamma$ is symmetric, the second marginal of $\rho$ is zero. The first marginal is zero as $\mathcal{A}$ contains constants. For any $a\in\mathcal{A}$ and all Borel sets $B\subset\Omega$ we have
\begin{equation*}
\int_{B\times\Omega}a(\omega_2)\, d\rho(\omega_1,\omega_2)=0,
\end{equation*}
by the choice of $\eta$.
Proposition \ref{pro:mod} shows that $\pi'=\pi+\rho\in\Gamma_{\mathcal{F}}(\mu,\nu)$. 
Let us take a disintegration $\theta$ of $\gamma$ with respect to the projection $p_1$ on the first co-ordinate, see Theorem \ref{thm:dis}. Then setting for $\omega\in\Omega$ and all Borel sets $A\subset\Omega$
\begin{equation*}
\lambda'(\omega,A)=\lambda(\omega,A)+\int_{\Omega}\big(\eta_{\omega_2,\omega}(A)-\eta_{\omega,\omega_2}(A)\big)\frac{\mathbf{1}_{\Lambda\setminus \Delta}(\omega,\omega_2)}{\epsilon(\omega,\omega_2)}\, d\theta(\omega,\omega_2)
\end{equation*}
we see that $\lambda'$ is a disintegration of $\pi'$, i.e., $\lambda'\in T(\pi')$.
We see that $\xi_{\omega}\in \Xi(\lambda',\omega)$.

Since $\eta$ and $\mathbf{1}_{\Lambda\setminus \Delta}$ are jointly measurable, the claim about measurability follows by the Fubini's theorem.
\end{proof}

\begin{example}
Let us consider the following simple, yet elucidating example. Let $\mu=\frac12(\delta_{\omega_1}+\delta_{\omega_2})$ and let $\nu=\frac12(\nu_1+\nu_2)$ be such that 
\begin{equation*}
a(\omega_i)=\int_{\Omega}a\, d\nu_i\text{ for }i=1,2\text{ and all }a\in\mathcal{A}.
\end{equation*}
Then there exists a trivial $\mathcal{F}$-transport $\pi\in\Gamma_{\mathcal{F}}(\mu,\nu)$,
\begin{equation*}
\pi=\frac12(\delta_{\omega_1}\otimes\nu_1+\delta_{\omega_2}\otimes\nu_2).
\end{equation*}
One may ask whether there exists another $\mathcal{F}$-transport between $\mu$ and $\nu$.
The above lemma builds upon the following simple idea.
Suppose that there exist probability measures $\eta_1,\eta_2$ such that 
\begin{equation*}
\int_{\Omega}a\, d\eta_1=\int_{\Omega}a\, d\eta_2\text{ for all }a\in\mathcal{A}\text{ and }\eta_i\leq C\nu_i\text{ for }i=1,2\text{ and some }C>0.
\end{equation*}
Then taking
\begin{equation*}
\lambda(\omega_1,\cdot)=\nu_1-\frac1C\eta_1+\frac1C\eta_2,\lambda(\omega_2,\cdot)=\nu_2-\frac1C\eta_2+\frac1C\eta_1
\end{equation*}
yields a disintegration of a measure $\pi'\in\Gamma_{\mathcal{F}}(\mu,\nu)$. If moreover $\eta_1,\eta_2$ are such that $\eta_i\geq c\nu_i$ for $i=1,2$ and some $c>0$, then the supports of any maximal disintegration are equal to $\mathrm{supp}\nu$.

\end{example}

Let us invoke the following theorem of Kellerer, see \cite[Proposition 3.3., p. 424 and Proposition 3.5., p. 425]{Kellerer1984}, also  see \cite[Proposition 2.1.]{Beiglbock2009}. Below $p_1,p_2\colon\Omega\times\Omega\to\Omega$ denote projections on the first and on the second co-ordinate respectively.

\begin{theorem}\label{thm:kellerer}
Let $\mu,\nu$ be Radon probability measures on $\Omega$.  Let $\Sigma\subset\Omega\times\Omega$ be Borel. Then for any $\gamma\in \Gamma(\mu,\nu)$ we have $\gamma(\Sigma)=0$ if and only if there exist Borel sets $B_1,B_2 \subset\Omega$ such that 
\begin{equation*}
\Sigma\subset p_1^{-1}(B_1)\cup p_2^{-1}(B_2)\text{ and }\mu(B_1)=\nu(B_2)=0.
\end{equation*}
\end{theorem} 

\begin{corollary}\label{col:symmetric}
Let $\mu$ be a Radon probability measure on $\Omega$.  Let $\Sigma\subset\Omega\times\Omega$ be Borel. Then for any symmetric $\gamma\in \Gamma(\mu,\mu)$ we have $\gamma(\Sigma)=0$ if and only if there exists a Borel set $B\subset\Omega$ such that 
\begin{equation*}
\Sigma\subset p_1^{-1}(B)\cup p_2^{-1}(B)\text{ and }\mu(B)=0.
\end{equation*}
\end{corollary}
\begin{proof}
If  $\gamma(\Sigma)=0$ for any symmetric $\gamma\in\Gamma(\mu,\mu)$, then the same equality holds for any $\gamma\in\Gamma(\mu,\mu)$. Indeed, if $\gamma\in\Gamma(\mu,\mu)$, then 
\begin{equation*}
\hat{\gamma}=\frac12(\gamma+s_{\#}\gamma))
\end{equation*}
is symmetric and belongs to $\Gamma(\mu,\mu)$. Now, Theorem \ref{thm:kellerer} implies that there exist Borel sets $B_1,B_2\subset\Omega$ such that 
\begin{equation*}
\Sigma\subset p_1^{-1}(B_1)\cup p_2^{-1}(B_2)\text{ and }\mu(B_1)=\mu(B_2)=0.
\end{equation*}
Taking $B=B_1\cup B_2$ completes the proof.
\end{proof}

We shall not consider the joint support of the family $\Gamma_{\mathcal{F}}(\mu,\nu)$, due to the reason explained in Example \ref{exa:jointbad}. For $\omega\in \Omega$ we set $S_{\omega}$ to be the support of $\lambda(\omega,\cdot)$, where $\lambda$ is a maximal disintegration of $\Gamma_{\mathcal{F}}(\mu,\nu)$ with respect to $p_1$, see Proposition \ref{pro:maxsup}.

\begin{lemma}\label{lem:support}
Let $\pi\in\Gamma_{\mathcal{F}}(\mu,\nu)$ and $\lambda\in T(\pi)$. Then there exists a Borel measurable set $B\subset\Omega$, with $\mu(B)=1$, such that whenever $\omega_1,\omega_2\in B$ then any measure in 
\begin{equation*}
\zeta_{\lambda,\omega_2}^{-1}(\Theta(\lambda,\omega_1)\cap \Theta(\lambda,\omega_2))\text{ is supported in } S_{\omega_1}.
\end{equation*}
\end{lemma}
\begin{proof}
Let us fix $\pi\in\Gamma_{\mathcal{F}}(\mu,\nu)$ and $\lambda\in T(\pi)$. 
Pick a sequence of measurable maps $\eta^k\colon \Lambda\setminus \Delta \to \mathcal{P}_{\xi(p)}(\Omega)$ as in Lemma \ref{lem:selection}. By Lemma \ref{lem:modification} for any symmetric $\gamma\in \Gamma(\mu,\mu)$ and any it's disintegration $\theta$, we may find $\pi'\in\Gamma_{\mathcal{F}}(\mu,\nu)$, $\lambda'\in T(\pi')$ and a Borel measurable function
\begin{equation*}
\epsilon^k\colon\Lambda\setminus\Delta\to (0,\infty)
\end{equation*}  
such that for any $\omega\in\Omega$
\begin{equation}\label{eqn:cosik}
\xi^k_{\omega}=\frac1{c(\omega)}\int_{\Omega}\eta^k_{\omega_2,\omega}\frac{\mathbf{1}_{\Lambda\setminus \Delta}(\omega,\omega_2)}{\epsilon^k(\omega,\omega_2)}\, d\theta(\omega,\omega_2)
\end{equation}
belongs to $\Xi(\lambda',\omega)$.
By the definition of $S_{\omega}$, for $\mu$-almost every $\omega\in\Omega$,
\begin{equation}\label{eqn:somega}
\lambda'(\omega, S_{\omega}^c)=0\text{ and therefore also }\xi^k_{\omega}(S_{\omega}^c)=0.
\end{equation}
Therefore, integrating (\ref{eqn:cosik}), we get that  
\begin{equation*}
\int_{\Omega}\eta^k_{\omega_2,\omega}(S_{\omega}^c)\mathbf{1}_{\Lambda\setminus \Delta}(\omega,\omega_2)\, d\gamma(\omega,\omega_2)=0.
\end{equation*}
It follows by Corollary \ref{col:symmetric} that there is a Borel set $B_k\subset\Omega$, with $\mu(B_k)=1$ such that
\begin{equation*}
\eta^k_{\omega_2,\omega}(S_{\omega}^c)\mathbf{1}_{\Lambda\setminus \Delta}(\omega,\omega_2)=0\text{ for all }\omega,\omega_2\in B_k.
\end{equation*}
Let $B$ denote the intersection of the sets $(B_k)_{k=1}^{\infty}$. Then $B$ is Borel and $\mu(B)=1$.
Now, by density, as stipulated in Lemma \ref{lem:selection}, it follows that whenever $\omega_1,\omega_2\in B$, then for any $\eta_1\in \Xi(\lambda,\omega_1)$ for which there exists  $\eta_2\in\ \Xi(\lambda,\omega_2)$ such that 
\begin{equation*}
\zeta_{\lambda,\omega_2}(\eta_1)=\zeta_{\lambda,\omega_2}(\eta_2)
\end{equation*}
we have
\begin{equation*}
\eta_2(S_{\omega_1}^c)=0\text{ for all }\omega_1\in B.
\end{equation*}
Here we use a result \cite[Theorem 8.2.3., p. 184]{Bogachev20072} on weak convergence of measures together with Lemma \ref{lem:regular}.
\end{proof}

\begin{remark}\label{rem:essentially}
As we suppose that $\Omega$ is separable and $\sigma$-compact, Remark \ref{rem:standard} shows that the Borel $\sigma$-algebra on $\Omega$ is countably generated. Therefore Theorem \ref{thm:dis} shows that any disintegration of $\pi\in\Gamma_{\mathcal{F}}(\mu,\nu)$ is essentially unique. It follows readily that if Lemma \ref{lem:support} holds true for some $\lambda\in T(\pi)$ then it holds true for all such disintegrations.
\end{remark}

\begin{remark}
One can show that, for $\mu$-almost every $\omega\in\Omega$, the weak closure of $\Xi(\omega)$ is equal to the set of measures $\eta\in\mathcal{P}_{\xi(p)}(\Omega)$ supported on $S_{\omega}$. This can yield a slightly alternative proof of the above lemma.

Phrased in other words, the proof relies on saying that the set of measures with support in $S_{\omega}$ is a weakly closed face in the space of Borel probability measures. It can be readily shown that any weakly closed face has to be of this form.
\end{remark}

\begin{example}
Let $\mathcal{A}$ be the space of affine functions on $\mathbb{R}^2$. Let $\nu_i$ be the uniform measure on the unit circle $S_i$ centred at $\omega_i\in\mathbb{R}^2$, $i=1,2$. Let 
\begin{equation*}
\mu=\frac12(\delta_{\omega_1}+\delta_{\omega_2}), \nu=\frac12(\nu_1+\nu_2).
\end{equation*}
Then $\mu\prec_{\mathcal{F}}\nu$.

For $i=1,2$ consider the set $\Theta_i$ of all functionals $\theta$ on $\mathcal{A}$ of the form
\begin{equation*}
\mathcal{A}\ni a\mapsto\int_{\Omega}a\, d\eta\text{ for some }\eta\in\mathcal{P}(\Omega), \eta\leq C\nu_i\text{ for some }C>0.
\end{equation*}
Then $\Theta_i$ can be identified with the relative interior of the convex hull of the unit circle $S_i$. 
Now, if $\omega_1,\omega_2$ are in distance at least two, then $\Theta_1\cap\Theta_2=\emptyset$, yet if the distance is equal to two, then the unit circles $S_1$ and $S_2$ do intersect.
If now the distance is less than two, then $\Theta_1\cap\Theta_2\neq\emptyset$. We see by Lemma \ref{lem:support} that in this case the supports of a maximal disintegration $\lambda$ are both equal to the union of the two circles.
\end{example}

\section{Notions of interia}\label{s:interia}

Before we proceed to the proofs of the main results we shall discuss various notions of interia a convex set.

\begin{definition}
Let $K$ be a convex set in a topological vector space. 

Let $\mathrm{rint}K$ denote the relative interior of $K$, i.e., the interior of $K$ viewed as a subset of its affine hull $\mathrm{Aff}K$.

Let $\mathrm{aint}K$ denote the algebraic relative interior of $K$, i.e., the set of all $k\in K$ for which the conical hull $\mathrm{Cone}(k-K)$ is a linear subspace. 

Let $\mathrm{qint}K$ denote the quasi-relative interior of $K$, i.e., the set of all $k\in K$ for which the closed conical hull $\mathrm{clCone}(k-K)$ is a linear subspace. 
\end{definition}

We shall need the following result, see e.g. \cite[Theorem 1.1.2, p. 4]{Zalinescu2002}.

\begin{theorem}\label{thm:interior}
Let $K$ be a convex subset of a topological vector space. If $\mathrm{rint}K\neq \emptyset$, then $\mathrm{rintcl}K=\mathrm{rint}K$ and $\mathrm{clrintK}=\mathrm{cl}K$.
\end{theorem}

Note that if $K$ is a subset of a finite-dimensional topological vector space, then it has non-empty relative interior.

\subsection{Gleason parts}\label{s:gleason}

If $K$ is a convex set in an infinite-dimensional topological vector space, it might happen that it has empty relative interior.
The following example is \cite[Example 5-3., p. 18]{Bear1970}.

\begin{example}\label{exa:interior}
There exist convex sets with empty algebraic interior. As an example consider the space of all Borel probability measures $\mathcal{P}([0,1])$ on $[0,1]$. For any $t\in [0,1]$ let $\delta_{t}$ denote the probability measure concentrated at the point $t$. Suppose that $\mu$ belongs to the algebraic interior of $\mathcal{P}([0,1])$.  For $t\in [0,1]$, there exists $r>0$ such that 
\begin{equation*}
(1+r)\mu-r\delta_t\in\mathcal{P}([0,1]).
\end{equation*}
That implies that $\mu(\{t\})\geq \frac{r}{1+r}>0$. This inequality cannot hold true for uncountably many $t\in [0,1]$. Thus $\mathrm{aint}\mathcal{P}([0,1])=\emptyset$.
\end{example}

The example shows that the relative interior of a convex set is not always a useful in the infinite-dimensional context.

\begin{definition}\label{def:gleason}
Let $K$ be a convex subset of a vector space. Let $k_1,k_2\in K$. We shall write $k_1\sim k_2$ whenever there exists a line segment 
\begin{equation*}
(l_1,l_2)=\{tl_1+(1-t)l_2\mid t\in (0,1)\}
\end{equation*}
containing $k_1,k_2$ and contained in $K$. An equivalence class $G(k,K)$ of a point $k\in K$ we shall call the \emph{Gleason part} of $k$ in $K$.
\end{definition}

We refer the reader to \cite[p. 122]{Alfsen1971} and \cite{Bear1970} for an account on the topic of Gleason parts.

\begin{proposition}\label{pro:charac}
Let $\mathcal{A}\subset\mathcal{C}_{\xi(p),0}(\Omega)$ be a linear subspace of continuous functions that  contains constants. Let $S\subset \Omega$ and let
\begin{equation*}
K=\mathrm{clConv}_H\Phi(S).
\end{equation*}
Then 
\begin{equation*}
K=\{a^*\in\mathcal{A}^*\mid a^*(a)\leq \sup a(S)\text{ for all }a\in\mathcal{A}\}
\end{equation*}
and
\begin{equation*}
K=\{a^*\in\mathcal{A}^*\mid a^*(1)=1, a^*(a)\geq 0\text{ for all }a\in\mathcal{A}\text{ non-negative on }S\}.
\end{equation*}
\end{proposition}
\begin{proof}
Clearly, 
\begin{equation*}
\Phi(S)\subset \{a^*\in\mathcal{A}^*\mid a^*(a)\leq \sup a(S)\text{ for all }a\in\mathcal{A}\},
\end{equation*}
so also $K$ is contained in  the considered  set.
If there was $a^*$ in the set but  not in  $K$, then by the Hahn--Banach theorem we could find $\epsilon>0$ and $a\in\mathcal{A}$ for which
\begin{equation*}
\Phi(\omega)(a)+\epsilon\leq a^*(a)\text{ for all }\omega\in S.
\end{equation*}
It follows  that
\begin{equation*}
\sup a(S)+\epsilon\leq a^*(a),
\end{equation*}
contrary  to the assumption on $a^*$. 
The second equality follows from the first one trivially, since $\mathcal{A}$ contains constants.
\end{proof}

\begin{lemma}\label{lem:harnack}
Let $\mathcal{A}\subset\mathcal{C}_{\xi(p),0}(\Omega)$ be a linear subspace that contains constants. Let $S\subset\Omega$. Let $a^*\in\mathrm{clConv}_H\Phi(S)$. Then the Gleason part of $a^*$ in $\mathrm{clConv}_H\Phi(S)$ is the set of all functionals $b^*\in \mathrm{clConv}_H\Phi(S)$ for which there exist constants $C>c>0$ such that
\begin{equation}\label{eqn:boundd}
ca^*(a)\leq b^*(a)\leq C a^*(a)\text{ whenever }a\in\mathcal{A}\text{ is non-negative on }S.
\end{equation}
\end{lemma}
\begin{proof}
Proposition \ref{pro:charac} tells us that $K=\mathrm{clConv}_H\Phi(S)$ is equal to the set of all $a^*\in\mathcal{A}^*$ such that $a^*(1)=1$ and $a^*(a)\geq 0$ whenever $a\in\mathcal{A}$ is non-negative on $S$.

Let $b^*\in K$. Then $a^*\sim b^*$ is equivalent to existence of $\epsilon>0$ such that 
\begin{equation}\label{eqn:ext}
(b^*+\epsilon (b^*-a^*)\big)(a)\geq 0\text{ and }(a^*+\epsilon (a^*-b^*)\big)(a)\geq 0\text{ if }a\in\mathcal{A}\text{ is non-negative on }S.
\end{equation}
This shows that $b^*$ satisfies (\ref{eqn:boundd}) with $c=\frac{\epsilon}{1+\epsilon}$ and $C=\frac{1+\epsilon}{\epsilon}$.
Conversely, if (\ref{eqn:boundd}) is satisfied, then (\ref{eqn:ext}) is also satisfied with 
\begin{equation*}
\epsilon=\frac1{\max\{C,\frac1c\}-1}.
\end{equation*}
Thus, $a^*\sim b^*$.
\end{proof}

\section{$\mathcal{F}$-convex sets}\label{s:fconvex}

Throughout this section we shall consider a set $\Omega$ and a lattice cone of functions $\mathcal{F}$ on $\Omega$.

In Section \ref{s:existence} we have proven a general version of the Strassen theorem. We assumed that any function in $\mathcal{G}$ is of $p$-growth for some proper, non-negative function $p$ in the complete lattice cone $\mathcal{H}$ generated by $\mathcal{G}$. In \cite{Ciosmak2023}
we provide two theorems that characterise the spaces $\Omega$, equipped with a convex cone of functions, that admit such $p$. In particular, we concern convex cones that satisfy the maximum principle, and lattice cones generated by a linear subspace.

\begin{definition}
Let $\mathcal{F}$ be a cone of functions on topological space $\Omega$. Let $S\subset\Omega$. We define closed $\mathcal{F}$-convex hull of $S$ by the formula
\begin{equation*}
\mathrm{clConv}_{\mathcal{F}}S=\{\omega\in\Omega\mid f(\omega)\leq \sup f(S)\text{ for all }f\in\mathcal{F}\}.
\end{equation*}
\end{definition}

The lemma  and the proposition below are proven in \cite{Ciosmak2023}.

\begin{lemma}\label{lem:convex}
Suppose that $\mathcal{H}$ is a complete lattice cone generated by a convex cone $\mathcal{G}$. Let $S\subset\Omega$. Then
\begin{equation*}
\mathrm{clConv}_{\mathcal{H}}S=\{\omega\in\Omega\mid g(\omega)\leq \sup g(S)\text{ for all }g\in\mathcal{G}\}.
\end{equation*}
\end{lemma}

\begin{proposition}\label{pro:convex}
Let $\mathcal{H}$ be  the complete lattice cone generated by a linear space $\mathcal{A}$ of continuous functions on topological space $\Omega$. Let $S\subset\Omega$. Then
\begin{equation*}
\mathrm{clConv}_{\mathcal{H}}S=\Phi^{-1}(\mathrm{clConv}_H\Phi(S)),
\end{equation*}
where $G$ is the cone of convex functions on $\mathcal{A}^*$, lower semi-continuous with respect to $\sigma(\tau(\mathcal{A}))$-topology.
\end{proposition}

The above lemma and proposition allowed in \cite{Ciosmak2023} for a general definition of $\mathcal{A}$-convex set.

\begin{definition}\label{def:conv}
A set $S\subset\Omega$ is said to be $\mathcal{A}$-convex whenever $S= \Phi^{-1}(\mathrm{Conv}_H\Phi(S))$.
\end{definition}

\begin{remark}\label{rem:gleasonf}
In Section \ref{s:gleason} we have defined Gleason parts of a convex set. However,  in \cite{Bear1970} the Gleason parts are defined in the general setting of function spaces. 

Suppose that $\mathcal{A}$ is a linear space of continuous functions on $\Omega$. Let $\omega\in\Omega$ and $S\subset\Omega$ be a closed, $\mathcal{A}$-convex set. Then the Gleason part $G'(\omega,S)$ of $\omega$ in $S$ is set of all $\omega'\in S$ such that there exist $C>c>0$ such that
\begin{equation*}
ca(\omega)\leq a(\omega')\leq Ca(\omega)\text{ for all }a\in\mathcal{A}\text{ non-negative on }S.
\end{equation*}
This is to say, the Gleason part of $\omega$ is the set of all points in $S$ for which the Harnack inequality is satisfied.
Let us observe that
\begin{equation*}
    G'(\omega,S)=\Phi^{-1}(G(\Phi(\omega),\mathrm{clConv}_H\Phi(S))).
\end{equation*}
In particular, Gleason parts are $\mathcal{A}$-convex.
\end{remark}

\begin{remark}\label{rem:fgaconvexity}
    Suppose that $\mathcal{G}$ generates $\mathcal{F}$. Let $\mathcal{A}=\mathcal{G}-\mathcal{G}$. Then a closed set that is $\mathcal{A}$-convex is $\mathcal{F}$-convex as well. Indeed, if $S$ is $\mathcal{A}$-convex, then it is $\mathcal{G}$-convex as well since $\mathcal{G}\subset\mathcal{A}$.  Lemma \ref{lem:convex} shows that it is $\mathcal{F}$-convex. 
\end{remark}

\section{Irreducible components}\label{s:irreducible}

The irreducible components we shall define as the relative interia of the respective closed convex hulls of the images of the supports of $\lambda$, where $\lambda\in\Lambda_{\mathcal{F}}(\mu,\nu)$ is a maximal disintegration of $\Gamma_{\mathcal{F}}(\mu,\nu)$, cf. Section \ref{s:maximaldis}. If $\mathcal{G}$ is a linear subspace, we may equivalently define them as the Gleason parts of the evaluation functional in the corresponding closed convex set.

Let us recall the Gelfand map $\Phi\colon\Omega\to\mathcal{A}^*$ cf. Section \ref{s:gelfand}, which is given by the formula
\begin{equation*}
\Phi(\omega)(a)=a(\omega)\text{ for }\omega\in\Omega\text{ and }a\in\mathcal{A}.
\end{equation*}

\begin{definition}\label{def:ircfin}
Suppose that $\mathcal{G}$ is a convex cone of functions on $\Omega$. Let $\mathcal{F}$ be the lattice cone generated by $\mathcal{G}$ and let $\mathcal{A}=\mathcal{G}-\mathcal{G}$. Let $\mu,\nu\in\mathcal{P}_{\xi(p)}(\Omega)$ be such that $\mu\prec_{\mathcal{F}}\nu$. Let $\lambda\in\Lambda_{\mathcal{F}}(\mu,\nu)$ be a maximal disintegration of $\Gamma_{\mathcal{F}}(\mu,\nu)$. For $\omega\in\Omega$ we define its irreducible component as 
\begin{equation*}
    \mathrm{rint}\mathrm{clConv}_{H}\Phi(\mathrm{supp}\lambda(\omega,\cdot)).
\end{equation*}
We shall denote it by $\mathrm{irc}_{\mathcal{A}}(\mu,\nu)(\omega)$.
\end{definition}

\begin{remark}
If $\mathcal{G}$ is a linear subspace, and $\mathrm{clConv}_{H}\Phi(\mathrm{supp}\lambda(\omega,\cdot))$ is finite-dimensional, we shall show in Lemma \ref{lem:dense} that the irreducible component can be defined as the Gleason part
    \begin{equation*}
G(\Phi(\omega),\mathrm{clConv}_{H}\Phi(\mathrm{supp}\lambda(\omega,\cdot)))\subset\mathcal{A}^*.
\end{equation*}
Under the same assumptions, Lemma  \ref{lem:harnack} shows that the irreducible component $\mathrm{irc}_{\mathcal{A}}(\mu,\nu)(\omega)$ is the set of the functionals $a^*\in\mathcal{A}^*$ such that there exists $C>1$ such that for all $a\in\mathcal{A}$, non-negative on $\mathrm{supp}\lambda(\omega,\cdot)$ we have
\begin{equation*}
\frac1Ca(\omega)\leq a^*(a)\leq Ca(\omega),
\end{equation*}
i.e., the set of those functionals for which the Harnack inequality holds true.
\end{remark}

\begin{remark}
Note that if $\lambda'\in\Lambda_{\mathcal{F}}(\mu,\nu)$ is another maximal disintegration of $\Gamma_{\mathcal{F}}(\mu,\nu)$ then for $\mu$-almost every $\omega\in\Omega$ we see that
\begin{equation*}
\mathrm{supp}\lambda(\omega,\cdot)=\mathrm{supp}\lambda'(\omega,\cdot).
\end{equation*}
Therefore the irreducible components are well-defined $\mu$-almost everywhere.
\end{remark}

\begin{remark}
We  want to show that any two irreducible components are either equal or disjoint. To this aim, we shall apply Lemma \ref{lem:support}. To make use  of it, we  shall show that for $\mu$-almost every $\omega\in\Omega$ the irreducible component $\mathrm{irc}_{\mathcal{A}}(\mu,\nu)(\omega)$ is equal to the set of all elements $a^*$ of $\mathcal{A}^*$ for which there exists a function $g\in L^{\infty}(\Omega,\lambda(\omega,\cdot))$, with $g\geq c$ for some $c>0$, such that
\begin{equation*}
a^*(a)=\int_{\Omega}ag\, d\lambda(\omega,\cdot),
\end{equation*}
provided that $\mathrm{clConv}_H\Phi(\mathrm{supp}\lambda(\omega,\cdot))$ is finite-dimensional.
\end{remark}

\begin{lemma}\label{lem:dense}
  Suppose that $\mathrm{clConv}_H\Phi(\mathrm{supp}\lambda(\omega,\cdot))$ is finite-dimensional. Then 
  \begin{equation*}
   \mathrm{irc}_{\mathcal{A}}(\mu,\nu)(\omega)=\zeta_{\lambda,\omega}(\mathrm{aint}\Xi(\lambda,\omega)).   
  \end{equation*}
  If $\mathcal{G}$ is a linear subspace, then
   \begin{equation*}
   \mathrm{irc}_{\mathcal{A}}(\mu,\nu(\omega)=G\big(\Phi(\omega),\mathrm{clConv}_H\Phi\big((\mathrm{supp}\lambda(\omega,\cdot))\big)\big).
   \end{equation*}
\end{lemma}
\begin{proof}
Theorem \ref{thm:interior} shows that the relative interior of $\mathrm{clConv}_H\Phi(\mathrm{supp}\lambda(\omega,\cdot))$ and the relative interior of any of its convex and dense subsets coincide. 

Let us show that
\begin{equation}\label{eqn:closures}
\mathrm{cl}\zeta_{\lambda,\omega}(\mathrm{aint}\Xi(\lambda,\omega))= \mathrm{clConv}_H\Phi(\mathrm{supp}\lambda(\omega,\cdot)).
\end{equation}
Note that the set on the right hand-side of the above equality consists of all functionals in $\mathcal{A}^*$ that are given by integration with respect to a measure in $\mathcal{P}_{\xi(p)}(\Omega)$ that is supported in $\lambda(\omega,\cdot)$. Thus clearly
\begin{equation*}
\mathrm{cl}\zeta_{\lambda,\omega}(\Xi(\lambda,\omega))\subset \mathrm{clConv}_H\Phi(\mathrm{supp}\lambda(\omega,\cdot)).
\end{equation*}
If we had an element $a^*\in \mathrm{clConv}_H\Phi(\mathrm{supp}\lambda(\omega,\cdot))\setminus \mathrm{cl}\zeta_{\lambda,\omega}(\mathrm{aint}\Xi(\lambda,\omega))$, then by the Hahn--Banach theorem there would exist a measure $\rho\in\mathcal{P}_{\xi(p)}(\Omega)$ supported in $\lambda(\omega,\cdot)$, $a\in\mathcal{A}$ and $\epsilon>0$ such that for all $\eta\in \Xi(\lambda,\omega)$
\begin{equation}\label{eqn:hahnik}
\int_{\Omega}a\, d\eta>\epsilon+\int_{\Omega}a \, d\rho. 
\end{equation}
Let $\epsilon'>0$. Suppose that 
\begin{equation}\label{eqn:inf}
\mathrm{inf}a(\mathrm{supp}(\lambda(\omega,\cdot))>-\infty. 
\end{equation}
Since $a$ is continuous, this is equivalent to assuming that the essential infimum of $a$ with respect to $\lambda(\omega,\cdot)$ is finite. Let $\eta$ be a measure, absolutely continuous with respect to $\lambda(\omega,\cdot)$, with density
\begin{equation}
\frac{\mathbf{1}_{\{\omega'\in\Omega\mid a(\omega')<\mathrm{inf}a(\mathrm{supp}(\lambda(\omega,\cdot))+\epsilon'\}}}{\lambda(\omega, {\{\omega'\in\Omega\mid a(\omega')<\mathrm{inf}a(\mathrm{supp}(\lambda(\omega,\cdot))+\epsilon'\})}}.
\end{equation}
Assumption (\ref{eqn:inf}) ensures that the denominator is positive.
We see that (\ref{eqn:hahnik}) implies that
\begin{equation*}
\epsilon'+\mathrm{inf}a(\mathrm{supp}(\lambda(\omega,\cdot))>\epsilon'+\int_{\Omega}a\, d\rho.
\end{equation*}
Letting $\epsilon'$ tend to zero, we get a contradiction. A similar argument works when the assumption (\ref{eqn:inf}) is not satisfied. For  if 
\begin{equation}\label{eqn:infinf}
\mathrm{inf}a(\mathrm{supp}(\lambda(\omega,\cdot))=-\infty,
\end{equation}
then for any $M>0$ we can take $\eta\in\mathcal{P}_{\xi(p)}(\Omega)$ to be a measure, absolutely continuous with respect to $\lambda(\omega,\cdot)$, with density
\begin{equation}
\frac{\mathbf{1}_{\{\omega'\in\Omega\mid a(\omega')<-M\}}}{\lambda(\omega, {\{\omega'\in\Omega\mid a(\omega')<-M\})}}.
\end{equation}
Then (\ref{eqn:hahnik}) shows that for any $M>0$
\begin{equation*}
-M>\epsilon+\int_{\Omega}a\, d\rho.
\end{equation*}
Letting $M$ tend to infinity we get a contradiction.
This shows that (\ref{eqn:closures}) holds true and in consequence completes the proof of the first claim of the lemma.

For the proof of the second claim, observe that if $\mathcal{G}$ is a linear subspace, then
\begin{equation*}
    \Phi(\omega)(a)=\int_{\Omega}a\, d\lambda(\omega,\cdot)\text{ for all }a\in\mathcal{A}.
\end{equation*}
Therefore the Gleason part  $G(\Phi(\omega),\mathrm{clConv}_H\Phi(\mathrm{supp}\lambda(\omega,\cdot)))$ contains the image $\zeta_{\lambda,\omega}(\mathrm{aint}\Xi(\lambda,\omega))$, which is equal to the irreducible component. Since the Gleason part is algebraically open, and the considered sets are finite-dimensional, it coincides with the component.
\end{proof}

\begin{theorem}\label{thm:finite}
Suppose that $\mathrm{clConv}_H\Phi(\mathrm{supp}\lambda(\omega,\cdot))$ is finite-dimensional subset of $\mathcal{A}^*$. Then $\mathrm{irc}_{\mathcal{A}}(\mu,\nu)(\omega)$ is the set of all functionals $a^*$ in $\theta(\lambda,\omega)$ for which there exists $\eta\in\Xi(\lambda,\omega)$ and constant $c>0$ such that
\begin{equation*}
\eta\geq c\lambda(\omega,\cdot)\text{ and }a^*=\zeta_{\lambda,\omega}(\eta).
\end{equation*}
\end{theorem}
\begin{proof}
Let us observe that the set of measures $\eta\in\mathcal{P}_{\xi(p)}(\Omega)$ for which there exist constants $C>c>0$ such that
\begin{equation*}
C\lambda(\omega,\cdot)\geq \eta\geq c\lambda(\omega,\cdot)
\end{equation*}
is precisely the algebraic relative interior of $\Xi(\lambda,\omega)$. We need to show that 
\begin{equation*}
\mathrm{irc}_{\mathcal{F}}(\mu,\nu)(\omega)=\zeta_{\lambda,\omega}(\mathrm{aint}\Xi(\lambda,\omega)).
\end{equation*}
Let $Z\subset \mathcal{A}^*$ be a finite-dimensional linear space that contains $\mathrm{clConv}_H\Phi(\mathrm{supp}\lambda(\omega,\cdot)$.
Clearly, $\zeta_{\lambda,\omega}\colon\Xi(\lambda,\omega)\to Z$ is a continuous linear map. 
Therefore by \cite[Proposition 1.2.7, p. 15]{Zalinescu2002}, as $Z$ is finite-dimensional,
\begin{equation*}
\zeta_{\lambda,\omega}(\mathrm{aint}\Xi(\lambda,\omega))=\mathrm{rint}\zeta_{\lambda,\omega}(\Xi(\lambda,\omega)).
\end{equation*}
The theorem follows now by  Theorem \ref{thm:interior} and by Lemma \ref{lem:dense}, (\ref{eqn:closures}).
\end{proof}

\begin{remark}\label{rem:arbsup}
Let us note that the same reasoning, as the one presented in the proof of Lemma \ref{lem:dense}, shows that for any measure $\rho\in\mathcal{P}_{\xi(p)}(\Omega)$ such that
\begin{equation*}
\mathrm{clConv}_H\Phi(\mathrm{supp}\rho)=\mathrm{clConv}_H\Phi(\mathrm{supp}\lambda(\omega,\cdot)),
\end{equation*}
and any $a^*\in \mathrm{irc}_{\mathcal{A}}(\mu,\nu)(\omega)$ one can find a measure $\eta\in\mathcal{P}_{\xi(p)}(\Omega)$ and a constant $C>1$ such that
\begin{equation*}
a^*(a)=\int_{\Omega}a\, d\eta\text{  for all }a\in\mathcal{A},
\end{equation*}
and
\begin{equation*}
\frac1C\rho\leq\eta\leq C\rho.
\end{equation*}
\end{remark}

\subsection{Partitioning}

\begin{definition}\label{def:local}
    We shall say that an $\mathcal{F}$-transport $\pi\in\Gamma_{\mathcal{F}}(\mu,\nu)$ is local whenever for $\mu$-almost every $\omega\in\Omega$ the set
    \begin{equation*}
        \mathrm{clConv}_H\Phi(\mathrm{supp}\lambda(\omega,\cdot))\subset \mathcal{A}^*
    \end{equation*}
    is a finite-dimensional convex set, where $\lambda\in T(\pi)$. 
\end{definition}

\begin{remark}\label{rem:anylocal}
    It is immediate to observe that any $\mathcal{F}$-transport between $\mu$ and $\nu$ is local if and only if there exists a maximal $\mathcal{F}$-transport between $\mu$ and $\nu$ that is local. Moreover, if $\mathcal{G}$ is a linear subspace, then this holds if and only if for $\mu$-almost every $\omega\in\Omega$ the irreducible component $\mathrm{irc}_{\mathcal{A}}(\mu,\nu)(\omega)$, defined as the appropriate Gleason part, is a finite-dimensional subset of $\mathcal{A}^*$.
\end{remark}

\begin{theorem}\label{thm:partition}
Suppose that $\mu\prec_{\mathcal{F}}\nu$. Assume that $\Omega$ is separable in $\tau(\mathcal{G})$. Then there exists a Borel measurable set $B\subset\Omega$ with $\mu(B)=1$ such that whenever $\omega_1,\omega_2\in B$ then 
\begin{equation*}
\mathrm{irc}_{\mathcal{A}}(\mu,\nu)(\omega_1)\cap\mathrm{irc}_{\mathcal{A}}(\mu,\nu)(\omega_2)=\emptyset
\end{equation*}
or
\begin{equation*}
\mathrm{irc}_{\mathcal{A}}(\mu,\nu)(\omega_1)=\mathrm{irc}_{\mathcal{A}}(\mu,\nu)(\omega_2)\text{ provided that these sets are finite-dimensional.}
\end{equation*}
Moreover, if $\mathrm{clConv}_H\Phi(\mathrm{supp}\lambda(\omega,\cdot))$ is finite-dimensional, then
\begin{equation}\label{eqn:closed}
    \Phi(\mathrm{supp}(\lambda(\omega,\cdot))\subset\mathrm{cl}\big(\mathrm{irc}_{\mathcal{A}}(\mu,\nu)(\omega)\big)\text{ for }\mu\text{-almost every  }\omega\in\Omega
\end{equation}
The sets $\mathrm{irc}_{\mathcal{A}}(\mu,\nu)(\cdot)$ are the smallest convex sets that satisfy  the above condition.

Furthermore, if $\mathcal{B}=\mathcal{G}\cap(-\mathcal{G})$, then
 \begin{equation*}
\Phi( \omega)\in R_{\mathcal{B}}\mathrm{irc}_{\mathcal{A}}(\mu,\nu)(\omega)\text{ for }\mu\text{-almost every }\omega\in\Omega.
 \end{equation*}
In particular, if $\mathcal{G}$ is a linear subspace, then 
 \begin{equation*}
    \Phi( \omega)\in\mathrm{irc}_{\mathcal{A}}(\mu,\nu)(\omega)\text{  for }\mu\text{-almost every }\omega\in\Omega.
 \end{equation*}
\end{theorem}
\begin{proof}
Let $\lambda\in\Lambda_{\mathcal{F}}(\mu,\nu)$ be a maximal disintegration of $\Gamma_{\mathcal{F}}(\mu,\nu)$. Then, thanks to Lemma \ref{lem:support}, there exists a Borel set $B\subset\Omega$ such that whenever $\omega_1,\omega_2\in B$, then any measure in $\zeta_{\lambda,\omega_2}^{-1}(\theta(\lambda,\omega_1)\cap\theta(\lambda,\omega_2))$ is supported in $\mathrm{supp}\lambda(\omega_1,\cdot)$. Suppose now that $\omega_1,\omega_2\in B$ are such that 
\begin{equation}
\mathrm{irc}_{\mathcal{A}}(\mu,\nu)(\omega_1)\cap\mathrm{irc}_{\mathcal{A}}(\mu,\nu)(\omega_2)\neq\emptyset,
\end{equation}
and such that these sets are finite-dimensional.
By Theorem \ref{thm:finite} there exist measures $\eta_i\in\Xi(\lambda,\omega_i)$, for $i=1,2$, with $\mathrm{supp}\eta_i=\mathrm{supp}\lambda(\omega_i,\cdot)$, such that $\zeta_{\lambda,\omega_1}(\eta_1)=\zeta_{\lambda,\omega_1}(\eta_2)$. By the property of the set $B$ we see that $\mathrm{supp}\lambda(\omega_1,\cdot)=\mathrm{supp}\lambda(\omega_2,\cdot)$. This immediately implies that
\begin{equation*}
\mathrm{irc}_{\mathcal{A}}(\mu,\nu)(\omega_1)=\mathrm{irc}_{\mathcal{A}}(\mu,\nu)(\omega_2).
\end{equation*}
The claim (\ref{eqn:closed}) follows immediately by Lemma \ref{lem:dense}.
It implies that the irreducible components are the smallest convex sets that satisfy the required conditions, as the components are defined as the relative interia of the closed convex hulls of $\Phi(\mathrm{supp}\lambda(\omega,\cdot))$ for $\omega\in\Omega$ and the relative interior of the closed convex hull of a convex set is contained in this set, by Theorem \ref{thm:interior}.

As $\mu\prec_{\mathcal{F}}\nu$, there exists a maximal disintegration $\lambda\in\Lambda_{\mathcal{F}}(\mu,\nu)$. Let $\mathcal{B}=\mathcal{G}\cap(-\mathcal{G})$. Then, by Theorem \ref{thm:finite},  for any measure $\eta\in\Xi(\lambda,\omega)$ such that for some $c>0$, 
\begin{equation*}
    \eta\geq c\lambda(\omega,\cdot),
\end{equation*}
the functional
\begin{equation*}
    \mathcal{A}\ni a\mapsto\int_{\Omega}a \,d\eta\in\mathbb{R},
\end{equation*}
belongs $\mathrm{irc}_{\mathcal{A}}(\mu,\nu)(\omega)$. For $b\in\mathcal{B}$, we see however that
\begin{equation*}
    \Phi_{\mathcal{B}}(\omega)(b)=b(\omega)=\int_{\Omega}bh\, d\lambda(\omega,\cdot),
\end{equation*}
which proves the assertion.

If we assume that $\mathcal{G}$ is a linear subspace, then  $\mathcal{A}=\mathcal{G}$, and $\mathcal{A}$ separates points of $\Omega$ by the assumption.

The proof is complete.
\end{proof}

The following theorem is a refinement of the previous one and provides a description of the intersections of the faces of the sets $\mathrm{clConv}_H\Phi(\mathrm{supp}\lambda(\omega,\cdot))$, for $\omega\in\Omega$,  such that there exists $\eta\in\mathcal{P}_{\xi(p)}(\Omega)$ and $C>0$ such that $\eta\leq C\lambda(\omega,\cdot)$ and
\begin{equation*}
a^*(a)=\int_{\Omega}a\, d\eta\text{ for all }a\in\mathcal{A}.
\end{equation*}

\begin{theorem}\label{thm:finestrucfin}
Suppose that $\Omega$ is separable in $\tau(\mathcal{G})$. There exists a Borel measurable set $B\subset\Omega$ with $\mu(B)=1$ such that whenever $\omega_1,\omega_2\in B$  are such that $\mathrm{clConv}_H\Phi(\mathrm{supp}\lambda(\omega_1,\cdot))$ and $\mathrm{clConv}_H\Phi(\mathrm{supp}\lambda(\omega_2,\cdot))$ are finite-dimensional and if 
\begin{equation*}
a^*\in\mathrm{clConv}_H\Phi(\mathrm{supp}\lambda(\omega_1,\cdot))\cap \mathrm{clConv}_H\Phi(\mathrm{supp}\lambda(\omega_2,\cdot)),
\end{equation*}
is such that for $i=1,2$ there exists $\eta_i\in\mathcal{P}_{\xi(p)}(\Omega)$ and $C_i>0$ such that $\eta_i\leq C_i\lambda(\omega_i,\cdot)$ and
\begin{equation}\label{eqn:intrep}
a^*(a)=\int_{\Omega}a\, d\eta_i\text{ for all }a\in\mathcal{A}.
\end{equation}
then the Gleason parts of $a^*$ in
\begin{equation*}
\mathrm{clConv}_H\Phi(\mathrm{supp}\lambda(\omega_i,\cdot))
\end{equation*}
for $i=1,2$ are equal.
\end{theorem}
\begin{proof}
Let $\lambda\in\Lambda_{\mathcal{F}}(\mu,\nu)$ be a maximal disintegration of $\Gamma_{\mathcal{F}}(\mu,\nu)$. Then, thanks to Lemma \ref{lem:support}, there exists a Borel set $B\subset\Omega$ such that whenever $\omega_1,\omega_2\in B$, then any measure in $\zeta_{\lambda,\omega_2}^{-1}(\theta(\lambda,\omega_1)\cap\theta(\lambda,\omega_2))$ is supported in $\mathrm{supp}\lambda(\omega_1,\cdot)$. 
Let $\omega_1,\omega_2\in B$ be such that there is 
\begin{equation*}
a^*\in\mathrm{clConv}_H\Phi(\mathrm{supp}\lambda(\omega_1,\cdot))\cap \mathrm{clConv}_H\Phi(\mathrm{supp}\lambda(\omega_2,\cdot)),
\end{equation*}
and such that (\ref{eqn:intrep}) is satisfied.

 By the property of the set $B$ we see that 
 \begin{equation}\label{eqn:containfin}
 \mathrm{supp}\eta_i\subset\mathrm{supp}\lambda(\omega_j,\cdot)\text{ for }i\neq j, i,j=1,2.
 \end{equation}
As in Lemma \ref{lem:dense} we see that 
the Gleason part of $a^*$ in $\mathrm{clConv}_H\Phi(\mathrm{supp}\lambda(\omega_i,\cdot))$ is equal to the image via $\zeta_{\lambda,\omega_i}$ of the Gleason part of $\eta_i$, so that it is equal to the Gleason part of $a^*$ in $\mathrm{clConv}_H\Phi(\mathrm{supp}\eta_i)$. 
Now, for $i,j=1,2$, (\ref{eqn:containfin}) shows that the Gleason part of $a^*$ in $\mathrm{clConv}_H\Phi(\mathrm{supp}\lambda(\omega_i,\cdot))$ is contained $\mathrm{clConv}_H\Phi(\mathrm{supp}\lambda(\omega_j,\cdot))$.

The proof is complete.
\end{proof}

\begin{remark}
We cannot hope to get a result that would read that if $a^*$ belongs to the intersection of two sets $\mathrm{clConv}_H\Phi(\mathrm{supp}\lambda(\omega_i,\cdot))$, $i=1,2$, then its Gleason parts in the corresponding sets are equal.  For example, let us consider 
\begin{equation*}
\mu=\frac12(\delta_{(0,0)}+\delta_{(0,1)})\text{ and }\nu=\frac14(\delta_{(0,1)}+\delta_{(0,-1)})+\frac12\lambda_S,
\end{equation*}
where $S$ is the circle of radius $1$ centred at $(0,1)$, and $\lambda_S$ is the uniform measure on $S$. Then if $\omega_1=(0,0)$, $\omega_2=(0,1)$ is is readily shown that $\lambda(\omega_1,\cdot)=\frac12(\delta_{(0,1)}+\delta_{(0,-1)})$ and that $\lambda(\omega_2,\cdot)=\lambda_S$. Therefore, $\mathrm{clConv}_H\Phi(\lambda(\omega_2,\cdot))$ intersects the irreducible component $\mathrm{irc}_{\mathcal{A}}(\mu,\nu)(\omega_1)$, yet, it does not contain it.
\end{remark}

\begin{remark}
Let us  consider the set $\mathrm{clConv}_H\Phi(\mathrm{supp}\lambda(\omega,\cdot))$. Assume that it is finite-dimensional. Let $C\subset\Omega$ be a Borel set such that $\lambda(\omega,C)>0$.  Let us consider the set of all $a^*\in\mathcal{A}^*$ such that there exists a measure $\eta\in\mathcal{P}_{\xi(p)}(\Omega)$ such that
\begin{equation*}
a^*(a)=\int_{\Omega}a\, d\eta\text{ for all }a\in\mathcal{A}\text{ and }\lambda(\omega,C\cap\cdot)\leq\eta\leq C\lambda(\omega,C\cap\cdot)\text{ for some }C>0.
\end{equation*}
This is a convex set, which is relatively open, since $\mathrm{clConv}_H\Phi(\mathrm{supp}\lambda(\omega,\cdot))$ is finite-dimensional, cf. Lemma \ref{lem:dense}. 
\end{remark}

\subsection{Disintegration}\label{s:dispartition}

Let $X$ be a topological space. We shall consider multifunctions with values in the space $F(X)$ of closed subsets of the space $X$, equipped with Effros Borel $\sigma$-algebra generated by the sets
\begin{equation}
\{F\in F(X)\mid F\cap U\neq\emptyset\}\text{ with }U\text{ open in }X.
\end{equation}
We refer the reader to \cite{Kechris1995}[Chapter II, 12. C, p. 72] for an account on Effros Borel $\sigma$-algebra.

\begin{remark}\label{rem:omega}
Let us note that by a result of Beer, if $X$ is a Polish space,  then so is $F(X)\setminus \{\emptyset\}$, see \cite{Beer1991}.
Therefore, one may easily show that if $\Omega$ is a Souslin space, then so is the space of its closed subsets $F(\Omega)\setminus\{\emptyset\}$. Indeed, $\Omega$ is a continuous image of a Polish space.
\end{remark}

\begin{remark}\label{rem:prelim}
We shall consider the set  $I$ of all preimages via $\Phi$ of irreducible components, and the map $s\colon\Omega\to I$ defined by the formula
\begin{equation*}
s(\omega)=\Phi^{-1}(\mathrm{irc}_{\mathcal{A}}(\mu,\nu)(\omega))\text{ for }\omega\in\Omega.
\end{equation*}
Let $r\colon\Omega\to F(\Omega)\setminus\{\emptyset\}$ be defined by the formula
\begin{equation*}
r(\omega)=\mathrm{supp}\lambda(\omega,\cdot)\text{ for }\omega\in\Omega.
\end{equation*}
 Assume that for $\mu$-almost every $\omega\in\Omega$, $\mathrm{irc}_{\mathcal{A}}(\mu,\nu)(\omega)$ is finite-dimensional. 
By Lemma \ref{lem:support} and Theorem \ref{thm:finite} it follows that there exists a Borel set $B\subset\Omega$, $\mu(B)=1$ such that for $\omega_1,\omega_2\in B$
\begin{equation*}
\mathrm{irc}_{\mathcal{A}}(\mu,\nu)(\omega_1)=\mathrm{irc}_{\mathcal{A}}(\mu,\nu)(\omega_2)
\end{equation*}
if and only if
\begin{equation*}
r(\omega_1)=r(\omega_2).
\end{equation*}
Indeed, if $r(\omega_1)=r(\omega_2)$, then the corresponding relative interia of closed convex hulls coincide. Conversely, if the irreducible components intersect, then Lemma \ref{lem:support} shows that the corresponding supports of maximal disintegration are equal.
This shows that there exists an injective map $t\colon I\to F(\Omega)\setminus\{\emptyset\}$ such that $r=t\circ s$, that maps $\Phi^{-1}(\mathrm{irc}_{\mathcal{A}}(\mu,\nu)(\omega))$ to $r(\omega)$, for $\omega\in B$. 
This map, as it is a bijection onto its range, allows us to introduce topology on $I$, which also results in the corresponding Borel $\sigma$-algebra on $I$
\end{remark}

\begin{theorem}\label{thm:partdis}
Assume that $\Omega$ is separable in $\tau(\mathcal{G})$. Let $\mu\prec_{\mathcal{F}}\nu$. Suppose that any $\mathcal{F}$-transport between $\mu$ and $\nu$ is local.

Then there exists a $\sigma$-algebra $\mathcal{E}$ of subsets of $I$, a probability measure $\theta$ on $\mathcal{E}$, and a map 
\begin{equation*}
\mu(\cdot,\cdot)\colon I\times \sigma(\tau(\mathcal{G}))\to\mathbb{R}
\end{equation*} 
such that
\begin{enumerate}
\item\label{i:proba} for every $i\in I$, $\mu(i,\cdot)$ is a probability measure,
\item\label{i:measure} for every $B\in\sigma(\tau(\mathcal{G}))$, the map $I\ni i\mapsto \mu(i,B)\in\mathbb{R}$ is Borel measurable,
\item for all $B\in\sigma(\tau(\mathcal{G}))$ and $E\in\mathcal{E}$
\begin{equation*}
\mu(B\cap s^{-1}(E))=\int_{E} \mu(i,B)\,d\theta(i),
\end{equation*}
\item  $\mu(i,\cdot)$ is concentrated on $s^{-1}(i)$ for $\theta$-almost every $i\in I$.
\end{enumerate}

\end{theorem}
\begin{proof}
Let  $\lambda\in\Lambda_{\mathcal{F}}(\mu,\nu)$ be a maximal disintegration.

Each compact subset of $\Omega$ is a Polish space, see Remark \ref{rem:standard}. Since $\Omega$ is a countable union of Polish spaces, it is a Souslin space, see \cite[Theorem 6.6.6., p. 21]{Bogachev20072}, and we may apply Theorem \ref{thm:dis}. 
By Remark \ref{rem:omega}, $F(\Omega)\setminus\{\emptyset\}$ is a Souslin space.

We claim that $r$ is a measurable map with respect to the Effros Borel $\sigma$-algebra on $F(\Omega)$ and Borel $\sigma$ algebra on $\Omega$. To this end, we need to verify that for each open set $U\in \tau(\mathcal{G})$ the set
\begin{equation*}
\{\omega\in \Omega\mid \mathrm{supp}\lambda(\omega,\cdot)\cap U\neq \emptyset\}
\end{equation*}
is Borel measurable.
Note  that the set is equal to 
\begin{equation*}
\{\omega\in \Omega\mid \lambda(\omega,U)>0\neq \emptyset\},
\end{equation*}
which is indeed Borel measurable, thanks to the measurability properties of $\lambda$.

Theorem \ref{thm:dis} yields a map $\tilde{\mu}(\cdot,\cdot)\colon F(\Omega)\times\sigma(\tau(\mathcal{G}))\to\mathbb{R}$ such that 
\begin{enumerate}
\item for every $D\in  F(\Omega)$, $\tilde{\mu}(D,\cdot)$ is a probability measure,
\item for every $B\in\sigma(\tau(\mathcal{G}))$, the map $F(\Omega)\ni D\mapsto \tilde{\mu}(D,B)\in\mathbb{R}$ is Borel measurable,
\item for all $B\in \sigma(\tau(\mathcal{G}))$ and any Borel measurable $E\subset F(\Omega)$
\begin{equation*}
\tilde{\mu}(B\cap r^{-1}(E))=\int_{E} \tilde{\mu}(D,B)\,dr_{\#}\tilde{\mu}(D),
\end{equation*}
\item  $\tilde{\mu}(D,\cdot)$ is concentrated on $r^{-1}(\{D\})$ for $r_{\#}\mu$-almost every $D\in F(\Omega)$.
\end{enumerate}
Let $\mathcal{E}$ be the $\sigma$-algebra of preimages by $t$ of Effros Borel sets in $F(\Omega)$, i.e., the Borel $\sigma$-algebra  on $I$, see Remark \ref{rem:prelim}.
We define for $E\in\mathcal{E}$,
\begin{equation*}
\theta(E)=s_{\#}\mu(E).
\end{equation*}
We define for $i\in I$ and $B\in\sigma(\tau(\mathcal{G}))$
\begin{equation*}
\mu(i,B)=\tilde{\mu}(t(i),B), 
\end{equation*} 
cf. Remark \ref{rem:prelim}.
It is then readily verifiable that these satisfy the requirements of the theorem.
\end{proof}

\begin{remark}
    Let us note that if $\mathcal{G}$ is a linear space, then $s(\omega)=i$ if and only if $\omega\in i$. 
\end{remark}

\section{Polar sets}\label{s:polarity}

\begin{definition}\label{def:polar}
Let $\mu,\nu$ be two Radon probability measures on $\Omega$ in $\mathcal{F}$-order, $\mu\prec_{\mathcal{F}}\nu$. We shall say that a Borel set $A\subset\Omega\times\Omega$ is a polar set for $\mathcal{F}$-transports between $\mu$ and $\nu$, whenever $\pi(A)=0$ for any $\pi\in\Gamma_{\mathcal{F}}(\mu,\nu)$. We shall also say in this situation that $A$ is $\Gamma_{\mathcal{F}}(\mu,\nu)$-polar.
\end{definition}

We shall be interested in characterising polar sets for $\mathcal{F}$-transports between $\mu$ and $\nu$. Let us first present the following two simple propositions.

\begin{proposition}
 Let $S$ be the joint support of $\Gamma_{\mathcal{F}}(\mu,\nu)$. Suppose that $U\subset\Omega\times\Omega$ is an open set. Then $U$ is $\Gamma_{\mathcal{F}}(\mu,\nu)$-polar if and only if $U\subset S^c$.
\end{proposition}
\begin{proof}
The proof follows immediately from Definition \ref{def:joint}, of the joint support.
\end{proof}

\begin{remark}
Example \ref{exa:jointbad} shows that that there exist measures $\mu$ and $\nu$ for which there are Borel polar sets with respect to all transports whose complement strictly contains the complement of the joint support $S$ of $\Gamma_{\mathcal{F}}(\mu,\nu)$. If however, the sections of $U$ are open, the situation is well-understood thanks to the existence of maximal disintegrations.
\end{remark}

For $U\subset\Omega\times\Omega$ and $\omega\in\Omega$ we set
\begin{equation*}
U_{\omega}=\{\omega'\in\Omega\mid (\omega,\omega')\in U\}.
\end{equation*}

\begin{proposition}
Suppose that $\Omega$ is separable in  $\tau(\mathcal{G})$. Let $\lambda\in\Lambda_{\mathcal{F}}(\mu,\nu)$ be a maximal disintegration of $\Gamma_{\mathcal{F}}(\mu,\nu)$. Suppose that $U\subset\Omega\times\Omega$ is a Borel set such that for $\mu$-almost every $\omega\in \Omega$ its section $U_{\omega}$ is open. Then $U$ is $\Gamma_{\mathcal{F}}(\mu,\nu)$-polar if and only if 
 \begin{equation*}
 U_{\omega}\subset \mathrm{supp}\lambda(\omega,\cdot)^c\text{ for }\mu\text{-almost every }\omega\in \Omega.
 \end{equation*}
\end{proposition}
\begin{proof}
The proof follows immediately from Definition \ref{def:maxdis}, of maximal disintegration.
\end{proof}

The situation when the sections of $U$ are not open is more intricate. 

\begin{theorem}\label{thm:polarfin}
Suppose that $\Omega$ is separable in $\tau(\mathcal{G})$ and that any  $\mathcal{F}$-transport between $\mu$  and $\nu$ is local.

Suppose that $U\subset\Omega\times\Omega$ is a Borel set such that for $\mu $-almost every $\omega\in \Omega$
\begin{equation}\label{eqn:relat}
\Phi(U_{\omega})\subset\mathrm{irc}_{\mathcal{A}}(\mu,\nu)(\omega).
\end{equation}
Then $U$ is $\Gamma_{\mathcal{F}}(\mu,\nu)$-polar set if and only if there exist Borel sets $N_1,N_2\subset\Omega$ with
\begin{equation}\label{eqn:null1}
\mu(N_1)=0,\nu(N_2)=0
\end{equation}
and 
\begin{equation}\label{eqn:null}
U\subset (N_1\times\Omega)\cup (\Omega\times N_2).
\end{equation} 
\end{theorem}

Let us note that a Borel set $U$ is $\Gamma(\mu,\nu)$-polar set if and only if there exist Borel sets $N_1,N_2$ such that (\ref{eqn:null1}) and (\ref{eqn:null}) are satisfied, thanks to Theorem \ref{thm:kellerer}. Thus, the above theorem shows that if $U$ satisfies also (\ref{eqn:relat}), i.e., if the sections are contained in appropriate irreducible components, then this condition is necessary and sufficient for $U$ to be $\Gamma_{\mathcal{F}}(\mu,\nu)$-polar.

\begin{lemma}\label{lem:faces}
Suppose  that $\Omega$ is separable in $\tau(\mathcal{G})$. Suppose that $\mu\prec_{\mathcal{F}}\nu$. Suppose that for $\mu$-almost every $\omega\in\Omega$, $\mathrm{irc}_{\mathcal{A}}(\mu,\nu)(\omega)$ is finite-dimensional. Let $\lambda\in\Lambda_{\mathcal{F}}(\mu,\nu)$ be a maximal disintegration. Let $f\in\mathcal{F}$ be such that 
\begin{equation*}
\int_{\Omega}f\,d\mu=\int_{\Omega}f\,d\nu<\infty.
\end{equation*}
Then, for each $i\in I$,  
there exists $g_i\in\mathcal{G}$ such that
$f=g_i$ on $i$ and $\mu(i,\cdot)$-almost everywhere and $\nu(i,\cdot)$-almost everywhere,
where for $i\in I$
\begin{equation*}
\nu(i,\cdot)=\int_{\Omega}\lambda(\omega,\cdot)\, d\mu(i,\cdot),
\end{equation*}
and $\mu(\cdot,\cdot)$ is defined in Theorem \ref{thm:partdis}.
\end{lemma}
\begin{proof}
For any non-negative, bounded $h$ on $\Omega$, $\sigma(\tau(\mathcal{G}))$-measurable we have
\begin{equation}\label{eqn:gie}
\int_{\Omega}hf\, \, d\mu\leq\int_{\Omega}\Bigg(h(\omega)\int_{\Omega}h(\omega)f\, d\lambda(\omega,\cdot)\Bigg)\, d\mu(\omega).
\end{equation}
However, by the assumption,
\begin{equation}\label{eqn:ekwal}
\int_{\Omega}f\, d\mu=\int_{\Omega}\int_{\Omega}f\, d\lambda(\omega,\cdot)\, d\mu(\omega).
\end{equation}
If we had a strict inequality in (\ref{eqn:gie}) for some non-zero $h$, then we would have a strict inequality in (\ref{eqn:ekwal}) as well. This shows that in (\ref{eqn:gie}) there is an equality for all $h$. 
Therefore
\begin{equation}\label{eqn:lamb}
f(\omega)=\int_{\Omega}f\, d\lambda(\omega,\cdot)
\end{equation}
for $\mu$-almost every $\omega\in\Omega$. Note also that for any $g\in\mathcal{G}$
\begin{equation}\label{eqn:lamba}
g(\omega)\leq\int_{\Omega}g\, d\lambda(\omega,\cdot),
\end{equation}
for $\mu$-almost every $\omega\in\Omega$.
For $\omega\in\Omega$ pick $g\in\mathcal{G}$ be such that $g\leq f$ on $s(\omega)$ and $f(\omega)=g(\omega)$.  Such $g\in\mathcal{G}$ exists, since $\Phi(s(\omega))$ is finite-dimensional and so  is $\mathrm{Conv}_H \big(\Phi(s(\omega))\cup\{\omega\}\big)$. 

Then (\ref{eqn:lamb}) and (\ref{eqn:lamba}) show that 
\begin{equation*}
f(\omega)=g(\omega)\text{ and } f=g\text{ for }\lambda(\omega,\cdot)\text{-almost everywhere, for }\mu\text{-almost every }\omega\in\Omega.
 \end{equation*}
 
 It follows that $f=g$, $(\mu(i,\cdot)+\nu(i,\cdot))$-almost everywhere. 
 By Theorem \ref{thm:finite} it follows that the equality holds true also on $i$.
\end{proof}

\begin{lemma}\label{lem:selectionbis}
There exists an analytic measurable map $I\ni i\mapsto \omega(i)\in\Omega$ such that for all $i\in I$, $\omega(i)\in s^{-1}(i)$.
Moreover, if $N\subset\Omega$ is a Borel set of $\mu$-measure zero, we may arrange that $\omega(i)\notin N$ for $\theta$-almost all $i\in I$.
\end{lemma}
\begin{proof}
By the definition, $I$ is isomorphic to the range of $r$ in $F(\mathcal{A}^*)$. The latter is a Polish space, whence $I$ is a Souslin space. 
Moreover 
\begin{equation*}
\Omega\ni\omega\mapsto s(\omega)\in I
\end{equation*}
is Borel measurable, as $r$ is Borel measurable.
The space $\Omega$, as a countable union of Polish spaces is Souslin, see \cite[Theorem 6.6.6., p. 21]{Bogachev20072}. The space $\Omega\setminus N$ is a Souslin space as well, thanks to \cite[Theorem 6.6.7., p. 22]{Bogachev20072}.
Now, the range of $s(\Omega\setminus I)$, by the definition of the measure $\theta$ is $I\setminus M$, where $M$ is of $\theta$-measure zero.
It is moreover a Souslin space, see \cite[Theorem 6.7.3., p. 24]{Bogachev20072}, as $s$ is a Borel map, thanks to the definition of topology on $I$. 
We may therefore employ the Jankoff selection theorem \cite[Theorem 6.9.1., p. 34]{Bogachev20072} to infer that there exists an analytically measurable map
\begin{equation*}
I\setminus M\ni i\mapsto\omega(i)\in \Omega\setminus N
\end{equation*}
such that $s(\omega(i))=i$ for all $i\in  I\setminus M$. The proof is complete.
\end{proof}

Let us recall that analytic sets are universally measurable, i.e., measurable with respect to any complete probability measure defined on Borel sets, see e.g. \cite[Theorem 29.7, p. 227]{Kechris1995}.

\begin{lemma}\label{lem:joint}
Suppose that $\Omega$ is separable in $\tau(\mathcal{G})$. Let $I,\theta$ be as in Theorem \ref{thm:partdis}. Suppose that
\begin{equation*}
I\times\sigma(\tau(\mathcal{A}))\ni (i,B)\mapsto\mu(i,B)\in\mathbb{R}
\end{equation*}
 is a function such that
\begin{enumerate}
\item\label{i:probaprim} for every $i\in I$, $\mu(i,\cdot)$ is a Borel probability measure,
\item\label{i:measureprim} for every $C\in\sigma(\tau)$, the map $I\ni i\mapsto \mu(i,C)\in\mathbb{R}$ is Borel measurable.
\end{enumerate}
Suppose that $\Omega\ni\omega\mapsto B(\omega)\in F(\Omega)$ is a Borel measurable multifunction, with values in closed subsets of $\Omega$. Let $I\ni i\mapsto \omega(i)\in\Omega$ be $\theta$-measurable. Then
\begin{equation*}
i\mapsto \mu(i,B(\omega(i)))
\end{equation*}
is $\theta$-measurable. Moreover, for any bounded, Borel measurable function $b\colon\Omega\to\mathbb{R}$, 
\begin{equation*}
i\mapsto\int_{\Omega}b\, d\mu_{B(\omega(i))}(i,\cdot)
\end{equation*}
is $\theta$-measurable. Here
\begin{equation*}
\mu_{B(\omega(i))}(i,C)=\mu(i,C\cap B(\omega(i)))\text{ for }C\in\sigma(\tau(\mathcal{A})).
\end{equation*}
\end{lemma}
\begin{proof}
Let $(K_i)_{i=1}^{\infty}$ be an increasing sequence of compact sets that exhaust $\Omega$. Then each $K_i$ is a Polish space. Let $K_0=\emptyset$. For $i=1,2,\dotsc$ consider the set 
\begin{equation}\label{eqn:ki}
\{(\omega_1,\omega_2)\in (K_i\setminus K_{i-1})\times\Omega\mid \omega_2\in B(\omega_1)\cap K_i\}
\end{equation}
From \cite[Theorem 28.8, p. 220]{Kechris1995} and the assumption that $B$ is Borel measurable multifunction it follows that it is Borel. 
By \cite[Theorem 28.7, p. 220]{Kechris1995} it follows that there exist Borel sets $(B_n^i)_{n=1}^{\infty}$ in $\sigma(\tau(\mathcal{G}))$ and closed sets $ (V_n^i)_{n=1}^{\infty}$ such that the set in (\ref{eqn:ki}) can be written as
\begin{equation*}
\bigcup_{n=1}^{\infty}B_n^i\times V_n^i.
\end{equation*}
Therefore
\begin{equation*}
\{(\omega_1,\omega_2)\in\Omega\times\Omega\mid \omega_2\in B(\omega_1)\}=\bigcup_{i=1}^{\infty}\bigcup_{n=1}^{\infty}B_n^i\times V_n^i.
\end{equation*}
We see that there exist Borel sets $(C_n)_{n=1}^{\infty}$ and pairwise disjoint $(W_n)_{n=1}^{\infty}$ such that 
\begin{equation*}
B(\omega)=\bigcup\{W_n\mid n\text{ is such  that }\omega\in C_n\}.
\end{equation*}
Then for any $C\in\sigma(\tau(\mathcal{A}))$ the map
\begin{equation*}
i\mapsto\mu(i,C\cap B(\omega(i)))=\sum_{n=1}^{\infty}\mathbf{1}_{C_n}(\omega(i))\mu(i,C\cap W_n)
\end{equation*}
is $\theta$-measurable.
This proves the first assertion. Also the second assertion is proven for $b=\mathbf{1}_C$. The case of general measurable function $b$ follows. Indeed, by  linearity, measurability holds for step functions. Any non-negative measurable function is a limit of an increasing sequence of step functions, so the claim follows by  the monotone convergence theorem. 
\end{proof}

\begin{proof}[Proof of Theorem \ref{thm:polarfin}]
If (\ref{eqn:null}) is satisfied then clearly $U$ is a polar set for $\mathcal{F}$-transports between $\mu$ and $\nu$.

Suppose that $U\subset\Omega\times\Omega$ is a polar set such that (\ref{eqn:relat}) holds true. We shall show that $U$ is also a polar set with respect to all transports between $\mu$ and $\nu$. Then the conclusion will follow by Theorem \ref{thm:kellerer}. 

\textbf{Step 1 Preparations}

Let $\lambda\in\Lambda_{\mathcal{F}}(\mu,\nu)$ be a maximal disintegration.

 Theorem \ref{thm:partdis} yields a set $I$, a measure $\theta$ on a $\sigma$-algebra $\mathcal{E}$ of subsets of $I$ and a map
 \begin{equation*}
 \mu(\cdot,\cdot)\colon I\times \sigma(\tau(\mathcal{G}))\to\mathbb{R},
 \end{equation*}
 whose properties are specified by the theorem.

  In particular, $\mu(i,\cdot)$ is concentrated on $s^{-1}(i)$ for a single preimage $i$ of an irreducible component.
  For $\theta$-almost every $i$, $\mathrm{clConv}_H\Phi(i)$ is finite-dimensional, by the assumption.
  
  Let us pick a compact set  $K\subset i$, which we shall specify later.
 
 Let $\gamma\in \Gamma(\mu,\nu)$ be any transport and let $\phi\in \Lambda(\mu,\nu)$ be its disintegration. 
 
For $\omega\in\Omega$ we define measures $\phi_K(\omega,\cdot)$ by the formula
\begin{equation*}
\phi_K(\omega,B)=\phi(\omega,B\cap K).
\end{equation*} 
Consider measure
\begin{equation*}
\psi_K(\omega,\cdot)=\phi_K(\omega,\cdot)+\lambda(\omega,\cdot).
\end{equation*}
  Then
\begin{equation*}
\mathrm{clConv}_H\Phi(\mathrm{supp}(\psi_K\omega,\cdot))=\mathrm{clConv}_H\Phi(\mathrm{supp}\lambda(\omega,\cdot)).
\end{equation*}
By Remark \ref{rem:arbsup}, there exists $\eta\in\mathcal{P}_{\xi(p)}(\Omega)$ and $D>0$ such that
\begin{equation}\label{eqn:bouns}
\frac1D\psi_K(\omega,\cdot)\leq \eta\leq D\psi_K(\omega,\cdot)
\end{equation}
and for all $a\in\mathcal{A}$
\begin{equation}\label{eqn:mom}
\int_{\Omega}a\, d\eta=\int_{\Omega}a\, d\lambda(\omega,\cdot).
\end{equation}
In particular, for $g\in\mathcal{G}$,
\begin{equation}\label{eqn:momo}
g(\omega)\leq \int_{\Omega}g\, d\eta.
\end{equation}
For a fixed $D>0$, let $\tilde{\Omega}_D$ be a set of $\omega\in\Omega$ for which there exists $\eta\in\mathcal{P}_{\xi(p)}(\Omega)$ such that (\ref{eqn:bouns}), (\ref{eqn:mom}) and (\ref{eqn:momo}) hold  true.

By the Radon--Nikodym theorem, for each $\omega\in\tilde{\Omega}_D$, we can find two measures $\tilde{\eta}_{\omega}, \bar{\eta}_{\omega}$ such that 
\begin{equation}\label{eqn:rn}
\frac1D\phi_K(\omega,\cdot)\leq\tilde{\eta}_{\omega}\leq D\phi_K(\omega,\cdot), \frac1D\lambda(\omega,\cdot)\leq\tilde{\eta}_{\omega}\leq D\lambda(\omega,\cdot)
\end{equation}
and such that $\eta_\omega=\tilde{\eta}_{\omega}+\bar{\eta}_{\omega}$.
Reasoning in an analogous way to the way in the proof of Lemma \ref{lem:selection} and employing Theorem \ref{thm:kuratowski}, we may pick a Borel measurable selections 
\begin{equation*}
\tilde{\Omega}_D\ni\omega\mapsto \tilde{\eta}_{\omega}\in\mathcal{P}_{\xi(p)}(\Omega),\tilde{\Omega}_D\ni\omega\mapsto \bar{\eta}_{\omega}\in\mathcal{P}_{\xi(p)}(\Omega).
\end{equation*}
Let us define measures 
\begin{equation}\label{eqn:defmu}
\mu_{D,K}(i,\cdot)=\mu(i,\cdot\cap\tilde{\Omega}_D\cap s^{-1}( K)),
\end{equation}
\begin{equation}\label{eqn:defnu}
\nu_{D,K}(i,\cdot)=\int_{\tilde{\Omega}_D\cap s^{-1}( K)}\lambda(\omega,\cdot)\, d\mu(i,\omega),
\end{equation}
\begin{equation}\label{eqn:defnutilde}
\tilde{\nu}_{D,K}(i,\cdot)=\int_{\tilde{\Omega}_D\cap s^{-1}( K)}\tilde{\eta}_{\omega}\, d\mu(i,\omega),
\end{equation}
and
\begin{equation}\label{eqn:defnubar}
\bar{\nu}_{D,K}(i,\cdot)=\int_{\tilde{\Omega}_D\cap s^{-1}( K)}\bar{\eta}_{\omega}\, d\mu(i,\omega),
\end{equation}

\textbf{Step 2 $\mathcal{F}$-ordering for some positive $\epsilon$}

We shall show that if $\epsilon>0$ is sufficiently small then\footnote{There is no guarantee that the measure on the right hand-side of the inequality is non-negative.}
\begin{equation}\label{eqn:orders}
(1-\epsilon)\int_{\Omega}f\, d \mu_{D,K}(i,\cdot)\leq\int_{\Omega}f\,d \big(\nu_{D,K}(i,\cdot)-\epsilon(\tilde{\nu}_{D,K}(i,\cdot)+\bar{\nu}_{D,K}(i,\cdot)\big)
 \end{equation}
for any $f\in\mathcal{F}$. 

Note that for each $i\in I$, all measures $\nu_{D,K}(i,\cdot), \tilde{\nu}_{D,K}(i,\cdot),\bar{\nu}_{D,K}(i,\cdot)$ are supported on the closure of $i$, which, by the assumption, is a preimage by $\Phi$ of a  finite-dimensional convex set. Measure $\mu_{D,K}(i,\cdot)$ is concentrated on $s^{-1}(K)$. 

Employing the Gelfand transform, cf. Theorem \ref{thm:embed}, we may assume that $\Omega$ is a convex subset of a finite-dimensional linear space and that $\mathcal{A}$ consists of affine functions on that space. 
Let $\omega_0\in i$ be a fixed point. 
Let $S\subset i$ be such that $\mathrm{Conv}_H\Phi(S)$ is a convex polygon that contains $K$ in its relative interior, see e.g. \cite{Ciosmak2024}.

Let $V\subset S$ denote the finite set of vertices of $S$. Let $\mathcal{B}$ denote the set of all non-negative functions $b\in\mathcal{F}$ that vanish  at $\omega_0$ and  such that
\begin{equation*}
\max b(V)=1.
\end{equation*}
Since $V$ is the set of extreme points of $S$, it follows that any function $b\in\mathcal{B}$ is bounded on $S$ by one. 

In order to prove (\ref{eqn:orders}) let us first show that
\begin{equation}\label{eqn:in}
\inf\Bigg\{ \frac{\int_{\Omega}b \, d(\nu_{D,K}(i,\cdot)-\mu_{D,K}(i,\cdot))}{\int_{\Omega}b \, d(\tilde{\nu}_{D,K}(i,\cdot)+\bar{\nu}_{D,K}(i,\cdot)-\mu_{D,K}(i,\cdot))}\mid b\in\mathcal{B}\Bigg\}>0.
\end{equation}
Note that the quotient in (\ref{eqn:in}) is non-negative, since
\begin{equation}\label{eqn:ordering}
\mu_D(i,\cdot)\prec_{\mathcal{F}}\nu_{D,K}(i,\cdot)\text{ and }\mu_{D,K}(i,\cdot)\prec_{\mathcal{F}}\tilde{\nu}_{D,K}(i,\cdot)+\bar{\nu}_{D,K}(i,\cdot),
\end{equation}
by  the definitions (\ref{eqn:defmu}), (\ref{eqn:defnu}), (\ref{eqn:defnutilde}), (\ref{eqn:defnubar}), and by (\ref{eqn:rn}),  (\ref{eqn:mom}), (\ref{eqn:momo}).
By (\ref{eqn:ordering}) and Lemma \ref{lem:faces} we see that if for some $b\in\mathcal{B}$ we had
\begin{equation*}
\int_{\Omega}b\, d(\tilde{\nu}_{D,K}(i,\cdot)+\bar{\nu}_{D,K}(i,\cdot)-\mu_{D,K}(i,\cdot))=0,
\end{equation*}
then there would exist a function  $g\in\mathcal{G}$ such that 
\begin{equation*}
b=g\text{ on }i\text{ and }(\tilde{\nu}_{D,K}(i,\cdot)+\bar{\nu}_{D,K}(i,\cdot)+\mu_{D,K}(i,\cdot))\text{-almost everywhere. }
\end{equation*}
This would imply that $g=0$. Indeed, $g\circ\Phi^{-1}$ would yield a non-negative, affine function on $\mathrm{Conv}_H\Phi(i)$, that vanishes at a point $\Phi(\omega_0)$. 
It follows that 
\begin{equation*}
b=0\quad (\tilde{\nu}_{D,K}(i,\cdot)+\bar{\nu}_{D,K}(i,\cdot)+\mu_{D,K}(i,\cdot))\text{-almost everywhere,}
\end{equation*}
and hence also $\nu_{D,K}(i,\cdot)$-almost everywhere. By Jensen's inequality $b=0$ on $i$. 
This contradicts the fact that $\max b(V)=1$.

Suppose that (\ref{eqn:in}) was false. 
Then for any $n=1,2,\dotsc$, we would be able to find an element $b_n\in\mathcal{B}$ such that 
\begin{equation}\label{eqn:onen}
\frac{\int_{\Omega}b_n \, d(\nu_{D,K}(i,\cdot)-\mu_{D,K}(i,\cdot))}{\int_{\Omega}b_n \, d(\tilde{\nu}_{D,K}(i,\cdot)+\bar{\nu}_{D,K}(i,\cdot)-\mu_{D,K}(i,\cdot))}<\frac1n.
\end{equation}
Let us define
  \begin{equation*}
  b=\limsup_{n\to\infty}b_n.
  \end{equation*}
  Then $b\in\mathcal{B}$. Let us note here that to show that $\max b(V)=1$, we might need to pass to a subsequence of elements that attain the maximum at a fixed vertex of $S$.

Observe now that (\ref{eqn:onen}) would mean that
\begin{equation}\label{eqn:onetwo}
\int_{\Omega}b_n\, d\bigg(\nu_{D,K}(i,\cdot)-\frac1{n}\bar{\nu}_{D,K}(i,\cdot)\bigg)<\int_{\Omega}b_n\, d\bigg(\frac1{n}\tilde{\nu}_{D,K}(i,\cdot)+\Big(1-\frac1n\Big)\mu_{D,K}(i,\cdot)\bigg).
\end{equation}
If $n$ is large enough, $n>2D$, then 
\begin{equation}\label{eqn:lowerbound}
\nu_{D,K}(i,\cdot)-\frac1{n}\bar{\nu}_{D,K}(i,\cdot)\geq \frac12\nu_{D,K}(i,\cdot).
\end{equation}
Note also that $\tilde{\nu}_{D,K}(i,\cdot)$ and $\mu_{D,K}(i,\cdot)$ are concentrated in $K$. Therefore $(b_n)_{n=1}^{\infty}$ is bounded $\tilde{\nu}_{D,K}(i,\cdot)$- and $\mu_{D,K}(i,\cdot)$-almost everywhere by one. 
Therefore, the right hand-side of (\ref{eqn:onetwo}) is uniformly bounded for all $n=1,2,\dotsc$.
Then, by (\ref{eqn:lowerbound}) we see that there exists a number $N>0$ such that for all $n=1,2,\dotsc$
\begin{equation*}
\int_{\Omega}b_n\, d\nu_{D,K}(i,\cdot)\leq N.
\end{equation*}
Employing the Koml\'os theorem, we may assume that the sequence converges $\nu_{D,K}(i,\cdot)$-almost everywhere and in $L^1(\Omega,\nu_{D,K}(i,\cdot))$ to a $\nu_{D,K}(i,\cdot)$-integrable function, which coincides $\nu_{D,K}(i,\cdot)$-almost everywhere with $b$. 
Passing to the limit in (\ref{eqn:onetwo}) yields now that 
  \begin{equation*}
\int_{\Omega}b\, d\nu_{D,K}(i,\cdot)\leq \int_{\Omega}b\, d\mu_{D,K}(i,\cdot).
  \end{equation*}
  However, by (\ref{eqn:ordering}), we also have the reverse inequality.
  An equality in the above inequality, together with Lemma \ref{lem:faces}, show that $b=0$ on $i$ and $\nu_{D,K}(i,\cdot)$-almost everywhere. Thus, also $b=0$ on $V$. This stands in contradiction with $\max b(V)=1$.
  
Let $\epsilon_0$ denote the positive infimum in (\ref{eqn:in}). Then we immediately see that (\ref{eqn:orders}) is satisfied for all $\epsilon<\epsilon_0$ for all $b\in\mathcal{B}$. Let us take $\epsilon<\frac1D$.

We shall show that it is valid for all $f\in\mathcal{F}$. Any such function is Lipschitz, so in particular it is finite on $V$.
We may pick $g\in\mathcal{G}$ such that $f\geq g$ on $i$  and $f(\omega_0)=g(\omega_0)$. If  $\max(f-g)(V)>0$, then we set
\begin{equation*}
b=\frac{f-g}{\max (f-g)(V)}.
\end{equation*}
Then  $b\in\mathcal{B}$. Since for this function (\ref{eqn:orders}) is satisfied, so it is for $f-g$ and also for $f$, since $g\in\mathcal{A}$ and we have (\ref{eqn:mom}).
 If $\max(f-g)(V)=0$, then $f$ does coincide with a function from $\mathcal{A}$ on $S$.
As $\mu_{D,K}(i,\cdot),\tilde{\nu}_{D,K}(i,\cdot)$ are supported on $S$ we see that (\ref{eqn:orders}) is  equivalent to 
\begin{equation}\label{eqn:ep}
0\leq \int_{\Omega}(f-g)\, d\big(\nu_{D,K}(i,\cdot)-\epsilon\bar{\nu}_{D,K}(i,\cdot)\big).
\end{equation}
As $\epsilon<\frac1D$, we see by (\ref{eqn:rn}) that 
\begin{equation*}
\nu_{D,K}(i,\cdot)-\epsilon\bar{\nu}_{D,K}(i,\cdot)\geq 0.
\end{equation*} 
As $f-g\geq 0$, (\ref{eqn:ep}) holds true.

Thus, (\ref{eqn:orders}) is proven for all $f\in\mathcal{F}$.

\textbf{Step 3 Specification of $K$  }

Let us now specify the set $K\subset i$ that we would like to take. Let $\omega\in\Omega$. Theorem \ref{thm:finite} shows that for $\mu$-almost every $\omega\in\Omega$, $s(\omega)$ consists of these $\omega'\in\Omega$ for which there exist $C>1$ and a measure $\eta\in\mathcal{P}_{\xi(p)}(\Omega)$ such that
\begin{equation}\label{eqn:condi}
\frac1C\lambda(\omega,\cdot)\leq\eta\leq C\lambda(\omega,\cdot)\text{ and }a(\omega')=\int_{\Omega}a\, d\eta\text{ for all }a\in\mathcal{A}.
\end{equation}
For $C>1$ let $\Omega_C(\omega)$ denote the set of $\omega'\in\Omega$ for which there is $\eta \in\mathcal{P}_{\xi(p)}(\Omega)$ such that (\ref{eqn:condi}) holds true.
Again, reasoning as in the proof of Lemma \ref{lem:selection} one can readily show that $\omega\mapsto\Omega_C(\omega)$ is a measurable multifunction.

Let us take a measurable selection, $I\ni i\mapsto \omega(i)\in\Omega$ such that for $\theta$-almost all $i\in I$, $\omega(i)\in i$ and moreover
\begin{equation*}
\bigcup\{\Omega_C(\omega(i))\mid C>1\}=i.
\end{equation*}
Such a selection exists thanks to Lemma \ref{lem:selectionbis} and thanks to the above observation that for $\mu$-almost every $\omega\in\Omega$ 
\begin{equation*}
\bigcup\{\Omega_C(\omega)\mid C>1\}=s(\omega).
\end{equation*} 
Moreover, the set of $\omega\in\Omega$ such that above the equality holds true, is a Borel set, since we have Borel measurable multifunctions $s$ and $\bigcup\{\Omega_C\mid C>1\}$ and we may  employ \cite[Exercise 12.11, ii), p. 76]{Kechris1995}.

For $i\in I$ we take $K=\Omega_C(\omega(i))$.

\textbf{Step 4 Measurability of $\epsilon$}

For $I\in I$, let $\epsilon(i)$ be given as minimum of $\frac1{2D}$ and the positive minimum of (\ref{eqn:in}) corresponding to the set $K=\Omega_C(\omega(i))$.

Since $\omega\mapsto\Omega_C(\omega)$ is a measurable multifunction, and $i\mapsto\omega(i)$ is $\theta$-measurable, Lemma \ref{lem:joint} shows that for each $D,C>1$ and $f\in\mathcal{F}$ the function
\begin{equation*}
i\mapsto \frac{\int_{\Omega}f \, d(\nu_{D,\Omega_C(\omega(i))}(i,\cdot)-\mu_{D,\Omega_C(\omega(i))}(i,\cdot))}{\int_{\Omega}f \, d(\tilde{\nu}_{D,\Omega_C(\omega(i))}(i,\cdot)+\bar{\nu}_{D,\Omega_C(\omega(i))}(i,\cdot)-\mu_{D,\Omega_C(\omega(i))}(i,\cdot))}
\end{equation*}
is $\theta$-measurable. Now, as a function of $f\in\mathcal{F}$, it is continuous with respect to topology of $\mathcal{D}_{p}(\Omega)$. Since $\mathcal{G}$ is separable, we see that there exists a countable subset $\mathcal{C}$ of $\mathcal{F}$ such that
\begin{equation*}
\inf\Bigg\{\frac{\int_{\Omega}f \, d(\nu_{D,\Omega_C(\omega(i))}(i,\cdot)-\mu_{D,\Omega_C(\omega(i))}(i,\cdot))}{\int_{\Omega}f \, d(\tilde{\nu}_{D,\Omega_C(\omega(i))}(i,\cdot)+\bar{\nu}_{D,\Omega_C(\omega(i))}(i,\cdot)-\mu_{D,\Omega_C(\omega(i))}(i,\cdot))}\mid f\in\mathcal{F}\Bigg\},
\end{equation*}
coincides with the respective infimum taken with respect to $\mathcal{C}$. We see, therefore, that the above infimum is measurable as a function of $i$.

Let  us denote the above  infimum by $I\ni i\mapsto \epsilon(i)\in\mathbb{R}$. We have shown that it is $\theta$-measurable.

\textbf{Step 5 Integration, construction of an $\mathcal{F}$-transport}

Let $f\in\mathcal{F}$. Taking into account Lemma \ref{lem:joint} and integrating (\ref{eqn:orders}), with the above choice of $\epsilon$ and $K$, with respect to $\theta$, yields 
\begin{equation*}
\int_{\Omega}f \, d\big(\tilde{\nu}_{D,C}+\bar{\nu}_{D,C}-\bar{\mu}_{D,C}\big)\leq\int_{\Omega}f\, d\big(\nu_{D,C}-\mu_{D,C}\big),
\end{equation*}
where 
\begin{align*}
& \mu_{D,C}=\int_I \mu_{D,\Omega_C(\omega(i))}\, d\theta(i)\\
& \bar{\mu}_{D,C}=\int_I\epsilon(i) \mu_{D,\Omega_C(\omega(i))}\, d\theta(i)\\
& \nu_{D,C}=\int_I \nu_{D,\Omega_C(\omega(i))}\, d\theta(i)\\
& \tilde{\nu}_{D,C}=\int_I\epsilon(i) \tilde{\nu}_{D,\Omega_C(\omega(i))}\, d\theta(i)\\
& \bar{\nu}_{D,C}=\int_I\epsilon(i) \bar{\nu}_{D,\Omega_C(\omega(i))}\, d\theta(i).
\end{align*}
 Note that
\begin{equation*}
\bar{\mu}_{D,C}\leq \mu_{D,C}\leq \mu
\end{equation*}
and since $\epsilon(i)\leq \frac1{2D}$, we see that  by (\ref{eqn:rn}),
\begin{equation*}
\tilde{\nu}_{D,C}+\bar{\nu}_{D,C}\leq\nu\text{ and } \nu_{D,C}\leq \nu.
\end{equation*}
Now, by Remark \ref{rem:finite},
\begin{equation*}
\mu-\bar{\mu}_{D,C}\prec_{\mathcal{F}}\nu-\big(\tilde{\nu}_{D,C}+\bar{\nu}_{D,C}\big).
\end{equation*}
Theorem \ref{thm:mart} proves that there exists an  $\mathcal{F}$-transport 
\begin{equation*}
\pi\in\Gamma_{\mathcal{F}}\Big(\mu-\bar{\mu}_{D,C},\nu-\big(\tilde{\nu}_{D,C}+\bar{\nu}_{D,C}\big)\Big).
\end{equation*}
Let us now explicitly define an $\mathcal{F}$-transport $\Phi_{D,C}\in\Gamma_{\mathcal{F}}(\bar{\mu}_{D,C},\tilde{\nu}_{D,C}+\bar{\nu}_{D,C})$ by setting for a Borel set $\Sigma\subset\Omega\times\Omega$ 
\begin{equation*}
\Phi_{D,C}(\Sigma)=\int_I \epsilon(i)\bigg(\int_{\Omega} \big(\tilde{\eta}_{\omega}+\bar{\eta}_{\omega}\big)(\Sigma_{\omega})\,d\mu_{D,\Omega_C(\omega(i))}(i,\omega) \bigg)\, d\theta(i).
\end{equation*}
Since $\eta_{\omega}=\bar{\eta}_{\omega}+\tilde{\eta}_{\omega}$ is a probability measure, we see that the first marginal of $\Phi$ is equal to $\bar{\mu}_{D,C}$. It is also readily visible that the second marginal is equal to $\tilde{\nu}_{D,C}+\bar{\nu}_{D,C}$. Since $\eta_{\omega}$ satisfies (\ref{eqn:mom}), $\Phi_{D,C}$ is indeed an $\mathcal{F}$-transport.

Now, $\pi+\Phi_{D,C}\in\Gamma_{\mathcal{F}}(\mu,\nu)$.

\textbf{Step 6 Conclusion}

If $U\subset\Omega\times\Omega$ is a polar set with respect to all $\mathcal{F}$-transports between $\mu$ and $\nu$, then $\Phi_{D,C}(U)=0$. 
By (\ref{eqn:rn}), we see that this implies that 
\begin{equation}\label{eqn:polar}
\int_I \bigg(\int_{\Omega} \phi(\omega,U_{\omega}\cap \Omega_C(\omega(i)))\,d\mu_{D,\Omega_C(\omega(i))}(i,\omega) \bigg)\, d\theta(i)=0.
\end{equation}
Note that for $\theta$-almost every $i\in I$, by Theorem \ref{thm:finite},
\begin{equation*}
\bigcup_{C=1}^{\infty}\Omega_C(\omega(i))=i\text{ and }
\bigcup_{D=1}^{\infty}\tilde{\Omega}_D=\Omega.
\end{equation*}
Therefore (\ref{eqn:polar}) and the assumption (\ref{eqn:relat}) show that 
\begin{equation}
\gamma(U)=\int_I \bigg(\int_{\Omega} \phi(\omega,U_{\omega})\,d\mu(i,\omega) \bigg)\, d\theta(i)=0.
\end{equation}
Thus, $U$ is a polar set with respect to all transports between $\mu$ and $\nu$. Theorem \ref{thm:kellerer} completes the proof.
\end{proof}

\begin{corollary}
Suppose that $\Omega$ is separable in $\tau(\mathcal{G})$, $\mu\prec_{\mathcal{F}}\nu$. The following are equivalent:
\begin{enumerate}
    \item\label{i:single} all $\mathcal{F}$-transports between $\mu$ and $\nu$ are local and there is a Borel set $B$, $\mu(B)=1$ such that
    \begin{equation*}
    \mathrm{irc}_{\mathcal{A}}(\mu,\nu)(\omega_1)=\mathrm{irc}_{\mathcal{A}}(\mu,\nu)(\omega_2)\text{ for all }\omega_1,\omega_2\in B,
    \end{equation*}
    \item\label{i:g} the set $\mathrm{clConv}_H\Phi(\mathrm{supp}\nu)$ is finite-dimensional and whenever $f\in\mathcal{F}$ is such that 
    \begin{equation*}
        \int_{\Omega}f\, d\mu=\int_{\Omega}f\, d\nu,
    \end{equation*}
    then there exists $g\in\mathcal{G}$ such that $f=g$, $(\mu+\nu)$-almost everywhere.
\end{enumerate}
\end{corollary}
\begin{proof}
That \ref{i:g} follows from \ref{i:single} is a consequence of Lemma \ref{lem:faces} and the assumption that all $\mathcal{F}$-transports are local.

The proof of Theorem \ref{thm:polarfin} shows that if \ref{i:g} holds true, then there exists an $\mathcal{F}$ transport $\pi\in\Gamma_{\mathcal{F}}(\mu,\nu)$ such that $\pi(U)>0$ if $U\subset\Omega\times\Omega$ is a Borel set, such that its respective projections have positive $\mu$ and $\nu$-measures. By Theorem \ref{thm:finite}, \ref{i:single} follows immediately.
\end{proof}

\begin{remark}\label{rem:extension}
Our method is robust enough to deal with a more general case. Namely, suppose that a polar set $U\subset\Omega\times\Omega$ for $\Gamma_{\mathcal{F}}(\mu,\nu)$ is such that for $\mu$-almost every $\omega\in\Omega$
\begin{equation*}
\Phi(U_{\omega})\subset \bigcap\Big\{\bigcup\{\theta(\lambda,\omega')\mid \omega'\in  N\}\bigg)\mid N\subset s(\omega), \mu(s(\omega),N)>0\Big\}.
\end{equation*}
We conjecture that it follows that it follows that $U\subset N_1\times\Omega\cup \Omega\times N_2$ for some sets $N_1,N_2$ such that $\mu(N_1)=0$ and $\nu(N_2)=0$.
Let us note that 
\begin{equation*}
\mathrm{irc}_{\mathcal{A}}(\mu,\nu)\subset \bigcap\{\theta(\lambda,\omega')\mid \omega'\in s(\omega)\},
\end{equation*}
so this would indeed be an extension of the previous result.
The strategy of the proof is the same. Instead of choosing in Step 1., a compact set in $i$, we pick a compact set in the preimage of a Gleason part of the above-considered set. Then Step 2., works as well as it did previously. However, it is unclear to us, how to pick a measurable choice of such compact sets. We do not pursue this direction, due to involved technicalities  and due to its little relevance to our interests.
Moreover, we conjecture that a similar approach allows to fully characterise polar sets for $\Gamma_{\mathcal{F}}(\mu,\nu)$, under the assumption that for $\theta$-almost all $i\in I$, $\mu(i,\cdot)$ or $\nu(i,\cdot)$ has at most countable support.
In that case, see Section \ref{s:abs}, one could pick a maximal disintegration such that any other disintegration would be absolutely continuous with respect to that one, on the preimages of the irreducible components.
Then, we conjecture, that the set $U\subset\Omega\times\Omega$ is polar for all $\mathcal{F}$-transports in $\Gamma_{\mathcal{F}}(\mu,\nu)$ if and only if
\begin{equation*}
U\subset N_1\times\Omega\cup \Omega\times N_2\cup \{(\omega,\omega')\in\Omega\times\Omega\mid \omega'\in\Phi^{-1}(\theta(\lambda,\omega))\}^c.
\end{equation*}
Note that when we may find a disintegration with respect to which all are absolutely continuous, then the set $\theta(\lambda,\omega)$ is independent on the choice of the maximal disintegration and on the point $\omega$ in the preimage of the irreducible component, cf.  proof of Lemma \ref{lem:support}.
\end{remark}

\begin{remark}
 The closure of the set $\theta(\lambda,\omega)$ contains $\Phi(\mathrm{supp}\lambda(\omega,\cdot))$ for $\mu$-almost every $\omega\in\Omega$. Indeed, its closure is equal to $\mathrm{clConv}_H\Phi(\mathrm{supp}\lambda(\omega,\cdot))$, see Remark \ref{lem:dense}.

Let us note however that it might happen that $\Phi^{-1}(\theta(\lambda,\omega))\cap\mathrm{supp}\lambda(\omega,\cdot)=\emptyset$. As an example, consider a uniform measure on a circle. Then $\Phi^{-1}(\theta(\lambda,\omega))$ is the open unit disc, hence it is disjoint from the circle, which is equal to the support of the considered measure.

Clearly, such situation occurs whenever $\mathrm{cclConv}_H\Phi(\mathrm{supp}\lambda(\omega,\cdot))$ is strictly convex and we $\lambda(\omega,\cdot)$ is non-atomic and supported on the boundary of this set.
\end{remark}

\begin{remark}
Let us exhibit a simple example of measures $\mu,\nu$, that shows that not all of polar sets with respect to all $\mathcal{F}$-transports between $\mu$ and $\nu$ can be written as a union 
\begin{equation}\label{eqn:unio}
N_1\times\Omega\cup\Omega\times N_2\cup\{(\omega,\omega')\in\Omega\times\Omega\mid \omega'\notin \Phi^{-1}(\mathrm{clConv}_H\mathrm{supp}\lambda(\omega,\cdot)\}.
\end{equation}
Let 
\begin{equation*}
\mu=\frac12(\delta_{0}+\delta_2),\nu=\frac14(\delta_{-1}+\delta_1+\lambda_{(1,3)}),
\end{equation*}
where $\lambda_{(1,3)}$ is the Lebesgue measure on the interval $(1,3)$.
Then $\mu\prec_{\mathcal{F}}\nu$ are in convex order. Observe that\begin{equation*}
\mathrm{irc}_{\mathcal{F}}(\mu,\nu)(0)=(-1,1),\mathrm{irc}_{\mathcal{F}}(\mu,\nu)(2)=(1,3),
\end{equation*}
but the set
\begin{equation*}
U=\{(2,1)\}
\end{equation*}
is a polar set, yet it is  not of the form (\ref{eqn:unio}). Let us note that this example is however captured by the reasoning presented in Remark \ref{rem:extension}.
\end{remark}

\section{Irreducible components in the non-local infinite-dimensional case}\label{s:issues}

Let us discuss the issues one encounters when one wants to extend the above results to the non-local infinite-dimensional setting.

Let $\lambda\in\Lambda_{\mathcal{F}}(\mu,\nu)$ be  a maximal disintegration. As we have seen in the proof of Theorem \ref{thm:partition} the methods developed in Section \ref{s:mixing} allow us to infer that if 
$\theta(\lambda,\omega_1)\cap\theta(\lambda,\omega_2)\neq\emptyset$ for some $\omega_1,\omega_2\in\Omega$ then the supports of measures 
$\eta_1,\eta_2\in\mathcal{P}_{\xi(p)}(\Omega)$ such that for some $C_1,C_2>0$
\begin{equation*}
    \eta_i\leq C_i\lambda(\omega_i,\cdot)\text{ for }i=1,2,
\end{equation*}
such that for all $a\in\mathcal{A}$
\begin{equation}\label{eqn:etai}
\int_{\Omega}a\, d\eta_1=\int_{\Omega}a\,d\eta_2,
\end{equation}
satisfy $\mathrm{supp}\eta_i\subset\mathrm{supp}\lambda(\omega_j,\cdot)$ for $i,j=1,2$, $i\neq j$. 
Therefore any notion of irreducible component that would allow for a use of results of Section \ref{s:mixing}, has to depend solely on the support of a maximal disintegration.
That is, the irreducible components $\mathrm{irc}_{\mathcal{A}}(\mu,\nu)(\cdot)$ should be defined so that
\begin{equation}\label{eqn:con}
\mathrm{irc}_{\mathcal{A}}(\mu,\nu)(\omega)\subset \theta(\lambda,\omega) \text{ for }\mu\text{-almost every }\omega\in\Omega,
\end{equation} 
and the definition depends only on the supports of a maximal disintegration $\lambda\in\Lambda_{\mathcal{F}}(\mu,\nu)$ and on $\omega\in \Omega$. Moreover, we would want the component to be dense in $\mathrm{clConv}_H\Phi(\mathrm{supp}\lambda(\omega,\cdot))$.

\subsection{Gleason parts}

Suppose that $\mathcal{G}$ is a linear subspace. We shall show in Theorem \ref{thm:nonresult} that we may not hope to define $\mathrm{irc}_{\mathcal{A}}(\mu,\nu)(\omega)$ as the Gleason part of the evaluation functional $\Phi(\omega)$ in $\mathrm{clConv}_H\Phi(\mathrm{supp}\lambda(\omega,\cdot))$,  unless the convex set $\mathrm{clConv}_H\Phi(\mathrm{supp}\lambda(\omega,\cdot))$ is finite-dimensional.

Like in the finite-dimensional setting, see Theorem \ref{thm:finite},
in the infinite-dimensional setting we would like to identify the image via $\zeta_{\lambda,\omega}$ of the Gleason part of $\lambda(\omega,\cdot)$ in $\Xi(\lambda,\omega)$ with the Gleason part of the evaluation functional $\Phi(\omega)\in\mathcal{A}^*$ in the set $\mathrm{clConv}_H\Phi(\mathrm{supp}\lambda(\omega,\cdot))$, where $\lambda\in\Lambda_{\mathcal{F}}(\mu,\nu)$ is a maximal disintegration.

\begin{theorem}\label{thm:nonresult}
Suppose that $\Omega$ is a separable and reflexive normed space, with norm $p$. Let $\xi\colon\mathbb{R}\to [1,\infty)$ be an increasing, convex function of superlinear growth. Let $\mathcal{A}$ be the space of all continuous affine functionals on $\Omega$. Then the following conditions are equivalent:
\begin{enumerate}
\item\label{i:gleason} for any measure $\lambda\in\mathcal{P}_{\xi(p)}(\Omega)$ the Gleason part $G(\Phi(\omega),\mathrm{clConv}_H\Phi(\mathrm{supp}\lambda))$  coincides  with the set of all functionals of the form 
\begin{equation*}
a\mapsto \int_{\Omega}a\,d\eta
\end{equation*}
with $\eta\in\mathcal{P}_{\xi(p)}(\Omega)$ being a Radon probability measure such that
\begin{equation*}
    c\lambda\leq \eta\leq C\lambda\text{ for some }C>c>0,
\end{equation*}
\item\label{i:closedL} $\mathcal{A}\subset L^1(\Omega,\lambda)$ is a closed subspace for any $\lambda\in \mathcal{P}_{\xi(p)}(\Omega)$ with $\mathrm{supp}\lambda=\Omega$,
\item \label{i:closedLone} $\mathcal{A}\subset L^1(\Omega,\lambda)$ is a closed subspace for some $\lambda\in \mathcal{P}_{\xi(p)}(\Omega)$ with $\mathrm{supp}\lambda=\Omega$,
\item\label{i:finited} $\mathcal{A}$ is finite-dimensional.
\end{enumerate}
\end{theorem}
\begin{proof}
Clearly, \ref{i:finited} implies \ref{i:closedL}. Now \ref{i:closedL} implies \ref{i:closedLone}. For if  $(\omega_i)_{i=1}^{\infty}$ is  a dense subset of $\Omega$, we may put
\begin{equation*}
\lambda=\frac{1}{\sum_{i=1}^{\infty}c_i}\sum_{i=1}^{\infty}c_i \delta_{\omega_i}
\end{equation*}
for some positive $(c_i)_{i=1}^{\infty}$ such that 
\begin{equation*}
    \sum_{i=1}^{\infty}c_i<\infty\text{ and }  \sum_{i=1}^{\infty}c_i\xi(p(\omega_i))<\infty,
\end{equation*}
e.g., $c_i=(2^i\xi(p(\omega_i)))^{-1}$ for $i=1,2,\dotsc$.
In this way, we get a measure $\lambda\in\mathcal{P}_{\xi(p)}(\Omega)$ such that $\mathrm{supp}\lambda=\Omega$.

To see that \ref{i:closedLone} implies \ref{i:finited}, let us pick $\lambda\in\mathcal{P}_{\xi(p)}(\Omega)$ with $\mathrm{supp}\lambda=\Omega$.
Let $\lambda_{\xi(p)}$ be a finite measure with density  $\xi(p)$ with respect to $\lambda$. Let $\mathcal{B}$ be the linear space of functions of the form $\frac{a}{\xi(p)}$ for $a\in\mathcal{A}$.
Then 
\begin{equation*}
\mathcal{B}\subset L^1(\Omega,\lambda_{\xi(p)})\text{ is a closed subspace and }\mathcal{B}\subset L^{\infty}(\Omega,\lambda_{\xi(p)}).
\end{equation*}
It follows by a theorem of Grothendieck \cite[Theorem 5.2 , p. 117]{Rudin1991}, that $\mathcal{B}$ is finite-dimensional in $L^1(\Omega,\lambda_{\xi(p)})$. Thus, $\mathcal{A}$ is finite-dimensional as well. Indeed, if we had an infinite set of linearly independent functions in $\mathcal{A}$, then the corresponding functions in $\mathcal{B}$ would be linearly independent in $L^1(\Omega,\lambda_{\xi(p)})$ as well, thanks to continuity of functions in $\mathcal{A}$ and the assumption that $\mathrm{supp}\lambda=\Omega$.

Let us assume that  \ref{i:gleason} is satisfied. 
Let $\lambda\in\mathcal{P}_{\xi(p)}(\Omega)$ be such that $\mathrm{supp}\lambda=\Omega$. Lemma \ref{lem:harnack} shows that the Gleason part $G(\Phi(\omega),\mathrm{clConv}_H\Phi(\mathrm{supp}\lambda))$ consists of all $a^*\in\mathcal{A}^*$ such that $a^*(1)=1$, since there are no functions $a\in\mathcal{A}$ such  that $a\geq 0$ on $\mathrm{supp}\lambda$.
Since the condition $a^*(1)=1$ is merely a normalisation, the assumption that \ref{i:gleason} is satisfied, implies that any functional $a^*\in\mathcal{A}^*$ is continuous, with respect to the norm on $ L^1(\Omega,\lambda)$. Therefore, any sequence of elements in $\mathcal{A}$ that converges weakly in $L^1(\Omega,\lambda)$ would also weakly converge in $\mathcal{A}$. As $\mathcal{A}$ is reflexive, the limit of the sequence also belongs to $\mathcal{A}$. This is to say, \ref{i:closedL} holds true.

Suppose now that \ref{i:finited} is satisfied. Then Theorem \ref{thm:finite} shows that \ref{i:gleason} holds true.
\end{proof}

\begin{remark}
For the definition of the irreducible components, we also know that $\lambda(\omega,\cdot)$ is a representing measure, i.e., for any $a\in\mathcal{A}$
\begin{equation*}
    a(\omega)=\int_{\Omega}a\, d\lambda(\omega,\cdot).
\end{equation*}
Note however that on a linear space, any measure with finite first moments represents its barycentre. Consequently, the assumption that a measure is a representing measure does not impose a substantial constraint. 
\end{remark}

Let us recall the extension theorem, due to M. Riesz \cite{Riesz1923}. 

\begin{theorem}\label{thm:riesz}
Let $K$ be a convex cone in a linear space  $X$. Let $Y\subset X$ be a subspace such that  $X=Y+K$. Let $\phi$ be a linear functional on $Y$ such that $\phi(y)\geq 0$ for all $y\in K\cap Y$.  Then there exists an extension of $\phi$ to a linear functional $\psi$ on $X$ such that $\psi(x)\geq 0$ for all $x\in K$.
\end{theorem}

\begin{proposition}\label{pro:extende}
Let $\mathcal{A}$ be a linear subspace of continuous affine functionals on a symmetric subset $\Omega$ of a normed space $X$, with norm $p$. Let $\xi\colon \mathbb{R}\to [1,\infty)$ be an increasing, convex function with superlinear growth.
 Suppose that $a^*\in\mathcal{A}^*$ is a continuous linear functional, such that
\begin{equation*}
a^*(a)\geq 0\text{ if }a\geq 0.
\end{equation*}
Then there exists a continuous functional $b^*$ on $\mathcal{A}+\mathbb{R}\xi(p)$ that extends $a^*$ and such that 
\begin{equation*}
b^*(b)\geq 0\text{ if }b\geq 0.
\end{equation*}
\end{proposition}
\begin{proof}
We put 
\begin{equation*}
\phi(\xi(p))=\sup \{a^*(a)\mid a\leq \xi(p)\}
\end{equation*}
and extend $\phi$ linearly to the space $\mathcal{A}+\mathbb{R}\xi(p)$. 
We shall show that 
\begin{equation*}
\sup \{a^*(a)\mid a\leq \xi(p)\}<\infty.
\end{equation*}
Since $a^*$ is continuous, there exists $D>0$ such that for $a\in\mathcal{A}$
\begin{equation*}
\abs{a^*(a)}\leq D\sup\bigg\{\frac{\abs{a(\omega)}}{ p(\omega)}\mid \omega\in \Omega\bigg\}<\infty,
\end{equation*}
as the supremum above is the norm of $a^*$.
The assumption that $\Omega$ is a symmetric subset implies that
\begin{equation*}
\sup \{a^*(a)\mid a\leq \xi(p)\}=\sup \{a^*(a)\mid \abs{a}\leq \xi(p)\}
\end{equation*}
If $\xi(p)\geq a$, then by the definition $\phi(\xi(p))\geq a^*(a)$, so
\begin{equation*}
\phi(\xi(p)-a)\geq 0.
\end{equation*}
If there exists $a\in\mathcal{A}$ such that  $a-\xi(p)\geq 0$, then $\phi(a-\xi(p))\geq 0$, since trivially
\begin{equation*}
\sup \{a^*(a)\mid a\leq \xi(p)\}\leq \inf \{a^*(a)\mid  \xi(p)\leq a\}.
\end{equation*}
This shows that 
\begin{equation*}
\phi(a+t\xi(p))\geq 0\text{ if } a+t\xi(p)\geq 0 \text{ for some }t\in\mathbb{R}.
\end{equation*}
\end{proof}

\begin{remark}\label{rem:automat}
We could prove that \ref{i:closedL} implies \ref{i:gleason} of Theorem \ref{thm:nonresult} in an alternative way.
Suppose that $\Omega$ is a symmetric subset of a normed space, with norm $p$.
Let $\lambda\in\mathcal{P}_{\xi(p)}(\Omega)$. Let $a^*\in\mathcal{A}^*$ be the functional
\begin{equation*}
\mathcal{A}\ni a\mapsto \int_{\Omega}a\,d\lambda\in\mathbb{R}.
\end{equation*} 
By Proposition \ref{pro:charac}, the Gleason part of $a^*$ in $\mathrm{clConv}_H\Phi(\mathrm{supp}\lambda)$ consists of all functionals $b^*\in\mathcal{A}^*$ such that 
\begin{equation*}
0\leq \frac1C a^*(a)\leq b^*(a)\leq Ca^*(a)\text{ whenever }a\in\mathcal{A}\text{ is non-negative on }\mathrm{supp}\lambda.
\end{equation*}
Let $c^*=b^*-\frac1Ca^*$. 

 We first extend the functional $c^*$ to a non-negative functional on $\mathcal{A}+\mathbb{R}\xi(p)$ employing Proposition \ref{pro:extende}.

Note that $\mathcal{A}$ is closed if and only if $\mathcal{B}=\mathcal{A}+\mathbb{R}p$ is closed in $L^1(\Omega,\lambda)$. Moreover, if we let $\mathcal{B}_+$ denote the set of all non-negative functions in $\mathcal{B}$, then  $\mathcal{B}=\mathcal{B}_+-\mathcal{B}_+$. Then  by \cite[Theorem, Chapter V, 5.5, p. 228 ]{Schaefer1971}, it  follows that any non-negative functional on $\mathcal{B}$ is automatically continuous.

Now, if $a^*$ as non-negative and continuous with  respect to  $L^1(\Omega,\lambda)$ norm on $\mathcal{A}$, $c^*$ is non-negative and continuous  with  respect to  $L^1(\Omega,\lambda)$ norm on $\mathcal{B}$. Thanks to \cite[Theorem, 5.4, p. 227 ]{Schaefer1971}, it can be extended to a non-negative, continuous functional on $L^1(\Omega,\lambda)$.

We are thus able to find a non-negative $h\in L^{\infty}(\Omega,\lambda)$ such that 
\begin{equation*}
b^*(a)=c^*(a)+\frac1Ca^*(a)=\int_{\Omega}a\Big(h+\frac1C\Big)\, d\lambda\text{ for all }a\in\mathcal{A}.
\end{equation*}
 Therefore \ref{i:gleason} holds true.

 \end{remark}
 
\begin{remark}
Let $\lambda\in\mathcal{P}_{\xi(p)}(\Omega)$ be such that $\mathrm{clConv}_H(\Phi(\mathrm{supp}\lambda))$ is a linear subspace. Then the characterisation of $G(\Phi(\omega),\mathrm{clConv}_H\Phi(\mathrm{supp}\lambda))$, as in \ref{i:gleason}, can be shown to imply that $\mathcal{A}\subset L^1(\Omega,\lambda)$ is closed and, in turn, finite-dimensional.
\end{remark}

\begin{remark}
    Let us mention two natural approaches to the definition of the irreducible components. One would be to define the components as the corresponding relative interia of the sets
    \begin{equation*}
\mathrm{clConv}_H\Phi(\mathrm{supp}\lambda(\omega,\cdot))\subset\mathcal{A}^*,\omega\in\Omega.
\end{equation*}
However, as shown in Example \ref{exa:interior}, the relative interior of a set in infinite-dimensional space might be empty.
A way to overcome this difficulty is to consider the quasi-relative interia of these sets. 
The advantage is that, under mild assumptions, the quasi-relative interia are non-empty, see \cite[Proposition 1.2.9, p. 18]{Zalinescu2002}. Moreover, \cite[Proposition 1.2.7, p. 15]{Zalinescu2002} shows that the image via $\zeta_{\lambda,\omega}$ of the quasi-relative interior of $\Xi(\lambda,\omega)$ is contained in the quasi-relative interior  $\mathrm{qintclConv}_G\Phi(\mathrm{supp}\lambda(\omega,\cdot))$.

Furthermore, the quasi-relative interior of $\Xi(\lambda,\omega)$ is the set of non-negative measures with bounded and almost everywhere positive density with respect to $\lambda(\omega,\cdot)$. That would be exactly the case that we can handle using our method of mixing. Note also that \cite[Proposition 1.2.7, p. 15]{Zalinescu2002}  shows that the above-defined component coincides with the previous definition in the finite-dimensional setting.

However, the same reasoning as presented in the proof of Theorem \ref{thm:nonresult} shows that a characterisation of the quasi-relative interior as the image via $\xi$ of the quasi-relative interior of  $\Xi(\lambda,\omega)$ is not generally possible in the infinite-dimensional case.
\end{remark}

\section{Further research}

\subsection{Maximal disintegrations with respect to absolute continuity}\label{s:abs}

If one of the measures $\mu$ or $\nu$ has at most countable support, then one may prove the existence of a disintegration with respect to which all other disintegrations are absolutely continuous $\mu$-almost surely. If $\nu$ has at most countable support, it suffices to take a maximal disintegration, since the inclusion of the supports implies readily the absolute continuity condition. When $\mu$ has at most countable support, then one may take also a maximal disintegration. Indeed, we see that for all disintegrations $\lambda\in\Lambda_{\mathcal{F}}(\mu,\nu)$ and for any $\omega$ in the support of $\mu$ there is 
\begin{equation*}
\lambda(\omega,\cdot)\ll \nu.
\end{equation*} 
Now, \cite{Berger1951} tells us that for a subset $M$ of Borel probability measures the following two conditions are equivalent:
\begin{enumerate}
\item\label{i:abs} there exists a probability measure $\rho$ such that for each $\eta\in M$ there is $\eta\ll\rho$,
\item\label{i:total} the set $M$ is separable with respect to the metric 
\begin{equation*}
d(\eta_1,\eta_2)=\sup\{\abs{\eta_1(A)-\eta_2(A)}\mid A\text{ is a Borel set}\}.
\end{equation*}
\end{enumerate}
Moreover, if \ref{i:total} is satisfied, $\rho=\sum_{i=1}^{\infty}\frac1{2^i}\eta_i$, where $(\eta_i)_{i=1}^{\infty}$ is a dense subset of $M$, will satisfy the requirements of \ref{i:abs}.

Now, we know by the above equivalence that the set 
\begin{equation*}
M=\{\lambda(\omega,\cdot)\mid \omega\in\Omega, \lambda\in\Lambda_{\mathcal{F}}(\mu,\nu)\}
\end{equation*}
is separable.
Let 
\begin{equation*}
\lambda'=\sum_{i=1}^{\infty}\frac1{2^i}\lambda_i,
\end{equation*}
where $(\lambda_i)_{i=1}^{\infty}\subset\Lambda_{\mathcal{F}}(\mu,\nu)$ is such that for each $\omega\in \Omega$, $(\lambda_i(\omega,\cdot))_{i=1}^{\infty}$ is a dense subset of 
\begin{equation*}
M(\omega)=\{\lambda(\omega,\cdot)\mid  \lambda\in\Lambda_{\mathcal{F}}(\mu,\nu)\}.
\end{equation*}
Then, again by the above equivalence, for all $\lambda\in\Lambda_{\mathcal{F}}(\mu,\nu)$,
\begin{equation}\label{eqn:maximal}
\lambda(\omega,\cdot)\ll \lambda'(\omega,\cdot)\text{ for all }\omega\in\Omega.
\end{equation}
This proves the asserted claim.

In the case of the existence of $\lambda\in\Lambda_{\mathcal{F}}(\mu,\nu)$ such that (\ref{eqn:maximal}) holds true, the theory that we developed is simpler.

 With the existence of such a kernel, we may define irreducible components and prove their pairwise disjointness, without the assumption of locality of all $\mathcal{F}$-transports.

For $\omega\in\Omega$, the irreducible components would be defined as the set of all functionals $a^*\in\mathcal{A}^*$ that are given as integration with respect to a probability measure equivalent to $\lambda(\omega,\cdot)$. By the maximality property, this definition is independent on the choice of the kernel $\lambda$.  

Moreover, in this way, we could also provide a more precise characterisation of the fine structure of intersections, akin to Theorem \ref{thm:finestrucfin}, and a full characterisation of polar sets with  respect to all $\mathcal{F}$-transports, cf. Remark \ref{rem:extension}.

\bibliographystyle{amsplain}
\bibliography{references2}

\end{document}